
\input amstex
\documentstyle{amsppt}

\magnification=1200
\parindent 20 pt
\def\nologo{\let\log@\empty}
\nologo
\NoBlackBoxes
\hsize 6.3 true in
\define\tz{\tilde z}
\define\cd{\cdot}
\define\ep{\endproclaim}
\define \a{\alpha}
\define \be{\beta}

\define \g{\gamma}
\define \lm{\lambda}
\define \s{\sigma}
\define \fa{\forall}

\define \hL{\hat{L}}
\define \iy{\infty}

\define \r{\rho}

\define \tg{\tilde{\gamma}}
\define \tP{\tilde{\Pi}}
\define \la{\langle}
\define \ra{\rangle}

\define \G{\Gamma}

\define \Ce{\varepsilon}

\define \CP{\Bbb C\Bbb P}

\define \CPt{\Bbb C\Bbb P^2}

\define \BR{\Bbb R}
\define \F{\Bbb F}

\define \Z{\Bbb Z}
\define \ve{\varepsilon}
\define \vp{\varphi}

\define \Dl{\Delta}
\define \dl{\delta}
\define \C{\Bbb C}

\define \1{^{-1}}\define \2{^{-2}}

\define \ri{\rightarrow} \define \Ri{\Rightarrow}

\define \Aff{\operatorname{Aff}}
\define \Gal{\operatorname{Gal}}

\define \un{\underline}
\define \ov{\overline}
\define \ub{\underbar}
\define \df{\dsize\frac}

\define \edm{\enddemo}
\define \Ss{S^{(6)}}
\define \sk{\smallskip}
\define \mk{\medskip}
\define \bk{\bigskip}

\topmatter

\heading{\bf Braid Group Technique in Complex Geometry V:\\
The Fundamental Group of a Complement of a Branch Curve of a Veronese
Generic Projection}\rm\endheading

\vskip 1.0in

\centerline{by}

\vskip 1.0in

$$\alignat 2
&\text{Boris Moishezon} \qquad \qquad \qquad \qquad
&& \text{Mina Teicher}\\
&\text{Dept. of Math.}
&&\text{Dept. of Math. \& Comp. Sci.}\\
&\text{Columbia University}
&& \text{Bar-Ilan University}\\
&\text{New York, N.Y. 10027}
&&\text{Ramat-Gan  52900}\\
& \text{U.S.A.}
&& \text{Israel}\endalignat$$
\baselineskip 20pt

\newpage

In this paper we compute $\pi_1(\C^2-S_3)$, \ $S_3$  the
branch curve of a generic projection of the Veronese surface $V_3$ on $\CPt. $
Throughout this paper we use $G$ for $\pi_1(\C^2-S_3).$
We also have a similar result concerning $\ov{G} =  \pi_1
(\CPt - S_3).$

Fundamental groups of complements of curves are very important invariants but
it is very difficult to compute them.  We obtained $G$ and $\ov G$
using our braid group techniques.

This paper is a continuation of ``Braid Group Techniques in Complex Geometry,
I,
II, III, IV.'' (BGT I, BGT II, BGT III, BGT IV for short.)
In BGT I we laid the foundation of our braid monodromy techniques
and applied them to line arrangements.
In BGT II we dealt with the braid monodromy of almost real curves and showed
how
to regenerate these formulas to cuspidal curves.
In BGT III we presented a series of generic projected degenerations of the
Veronese surface $V_3$ and the branch curve $S_3$ of its generic projection to
$\CPt$ to a union of 9 planes, and a branch curve $S^{(6)}$ which is a union of
lines.
In BGT IV we computed the braid monodromy of $S_3\cap\Bbb C^2$ using the braid
monodromy of $S^{(6)}$ and the regeneration rules proved in BGT II.
We obtained the factorized expression for the braid monodromy denoted by
$\varepsilon (18).$  In this paper we use $\varepsilon
(18)$ and the Van Kampen Theorem to compute $G$ and $\ov G.$
If the reader is only interested in the final results he can go directly to
Chapter VI.

\heading{{\bf CHAPTER 0.\quad Definition of Braid Group and Presentation of the
Van Kampen Theorem}}\endheading

We need certain terminology in order to formulate  the Van Kampen Theorem.

We first recall from BGT I (Section II) the definition of a good geometric base
of $\pi_1(D-K, *)$ for $K$  a finite set in a disc $D.$

\smallskip
\demo{Definition}\ $\underline{\ell(\g)}$.

Let $D$ be a disc.
Let $w_i$\ $i=1,\dots,n$ small discs in $\operatorname{Int}(D)$ s.t. $w_i\cap
w_j=\emptyset$ for $i\ne j$.
Let $u\in\partial D.$
 Let $\g$ be a simple path connecting $u$
with one of the $w_i$'s\,, say $w_{i_{0}}$\,, which does not meet any other
$w_j\,,\,j\ne i_0$\,.
We assign to $\g$ a loop (actually an element of $\pi_1(D-K,u)$) as follows.
Let $c$ be a simple loop equal to the (oriented) boundary of a small
neighborhood $ V$  of $w_{i_0}$ chosen such that $\g^{\prime}=\g-V\cap \g$ is a
simple path.  Then $\ell(\g)=\g^{\prime}\cup c\cup \g^{\prime\,-1}.$
We also use the same notation $\ell(\g)$   for the element of
$\pi_1(D-K,u)$ corresponding to $\ell(\g)$ (see Fig. 0.1).\enddemo

\midspace{1.50in}\caption{Fig. 0.1}

\demo{Definition} \ $\underline{g\text{-base (good geometric base)}}$

Let $D$ be a disk, $K\subseteq D,$ $\# K<\iy.$
Let $u\in D-K.$
Let $\{\g_i\}$ be a bush in $(D,K,u),$ i.e.,\ $\forall i,j$\
$\g_i\cap\g_j=u,$ \ $\forall i$\ $\g_i\cap K$ = one point, and
$\g_i$ are ordered counterclockwise around $u.$
Let $\G_i=\ell (\g_i) \in \pi_1 (D-K,u)$ be the loop around $K\cap\g_i$
determined by $\g_i.$
$\{\G_i\}$ is a $g$-base of $\pi_1(D-K, *).$\edm
\demo{Definition} $\un{\text{Braid group}\ B_n=B_n[D,K]}$

Let $D$ be
a closed disc in $\Bbb R^2,$ \ $K\subset D,$ \ $K$ finite.
Let $B$ be the group of
all diffeomorphisms $\beta$ of $D$ such that $\beta(K) = K\,,\, \beta
|_{\partial D} = \text{Id}_{\partial D}$\,.
For $\beta_1 ,\beta_2\in B$\,, we
say that $\beta_1$ is equivalent to $\beta_2$ if $\beta_1$ and $\beta_2$ induce
the same automorphism of $\pi_1(D-K,u)$\,.
The quotient of $B$ by this
equivalence relation is called the braid group $B_n[D,K]$ ($n= \#K$).
We sometimes denote by $\overline\beta$ the braid represented by $\beta.$
The elements of $B_n[D,K]$ are called braids.
\enddemo
\demo{Definition}\ \underbar{$H(\sigma)$, half-twist defined by $\sigma$}

Let
$D,K$ be as above.
Let $a,b\in K\,,\, K_{a,b}=K-a-b$ and $\sigma$ be a simple
path in $D-\partial D$ connecting $a$ with $b$ s.t. $\sigma\cap
K=\{a,b\}.$
Choose a small regular neighborhood $U$ of $\sigma$ and an
orientation preserving diffeomorphism $f:{\Bbb R}^2 \rightarrow {\Bbb C}^1$\
(${\Bbb C}^1$ is taken with the usual ``complex'' orientation) such that
$f(\sigma)=[-1,1]\,,\,$ \ $ f(U)=\{z\in{\Bbb C}^1 \,|\,|z|<2\}$\,.
Let $\alpha(r),r\ge 0$\,,
be a real smooth monotone function such that $
\alpha(r) = 1$ for $r\in [0,\tsize{3\over 2}]$ and
                $\alpha(r) =   0$ for $ r\ge 2.$

Define a diffeomorphism $h:{\Bbb C}^1 \longrightarrow {\Bbb C}^1$ as follows.
For $z\in {\Bbb C}^1\,,\, z= re^{i\varphi},$ let
$h(z) =re^{i(\varphi +\alpha(r))}$\,.
It is clear that on
$\{z\in{\Bbb C}^1\,|\,|z|\leq\tsize{3\over 2}\}$,\ $h(z)$ is the
positive rotation by $180^{\tsize{\circ}}$ and that
$h(z)=\text{Identity on }\{z\in{\Bbb C}^1\,|\,|z|\ge 2\}$\,, in
particular, on ${\Bbb C}^1 -f(U)$\,.
Considering $(f\circ h\circ f^{-1})|_{D}$ (we always take composition from
left to right), we get a diffeomorphism of $D$ which switches $a$ and $b$ and
is
the identity on $D-U$\,.
Thus it defines an element of $B_n[D,K],$ called the
half-twist defined by $\sigma$ and denoted $H(\sigma).$ \edm
 \demo{Definition}\  $\underline{\text{Frame of}\ B_n[D,K]}$

Let $D$ be a disc in $\BR^2.$
Let $K=\{a_1,\ldots ,a_n\},$\ $K\subset D.$
 Let $\sigma_1,\ldots ,\sigma_{n-1}$
be a system of simple paths in $D-\partial D$ such that each $\sigma_i$
connects
$a_i$ with $a_{i+1}$ and for
$$
i,j\in\{1,\ldots ,n-1\}\ ,\ i<j\quad ,\quad
\sigma_i\cap\sigma_j =
    \cases \emptyset \ \ &\text{if } |i-j|\ge 2\\
           a_{i+1} \ \ &\text{if } j=i+1\,.
    \endcases
$$
Let $H_i = H(\sigma_i)$\,.
We call the ordered system of (positive) half-twists $(H_1,\ldots ,H_{n-1})$ a
frame of $B_n[D,K]$ defined by $(\sigma_1,\ldots ,\sigma_{n-1})$\,, or a frame
of $B_n[D,K]$ for short. \enddemo

\demo{Notation}

$[A, B] = ABA^{-1} B^{-1}$.

$\la A, B\ra = ABAB^{-1}A^{-1}B^{-1}$.

$(A)_B= B^{-1} AB$.\edm

\proclaim\nofrills {Theorem.\ (E.~Artin's braid group presentation)}
\  $B_n$
is generated by the half-twists $H_i$ of a frame $\{H_i\}$ and all the
relations
between $H_1 ,\dots, H_{n-1}$ follow from $$
\align [H_i,H_j] &=1\qquad \text {if} \quad |i-j|>1,\\
\langle H_i, H_j \rangle &=1 \qquad \text {if} \quad |i-j|=1,
\\ &1 \leq i,j \leq n-1.
\endalign
$$
\endproclaim

\demo {Proof} See proof in \cite{MoTe4}, Chapter 4.\quad $\square$\edm
\smallskip
\demo{Proposition-Definition}\ $\underline{\Delta_n^2\ (\in B_n)}$

$\Delta_n^2=(H_1\dots H_n)^{n-1}$ for any frame $H_1\dots H_{n-1}$ of
$B_n.$\edm

We shall need the following definition:

\demo{Definition}\ $\underline{R(\Ce), N(\varepsilon), G(\varepsilon)\
\text{where}\ \varepsilon\ \text{is a factorized expression in}\ B_m}$

Let $D$ be a disk in $\C,$ $K\subseteq D,$ $\# K=m.$

Let $\F_m=\pi_1 (D-K).$
Consider the natural action of $B_m=B[D,K]$ on $\F_m,$ denoted by
$(\G) g_i.$

Let $\varepsilon=g_1\cdot\dots\cdot g_t$ be a factorized expression in $B_m.$

Let $\G_1\cdots\G_m$ be a good geometric base of $\F_m.$

Let $M(\varepsilon)$ be the subgroup of $B_m$ generated by $\{g_i\}^t_{i=1}.$

Let $R(\varepsilon)$ be the subgroup of $\F_m$ generated by
$\{(\G_j)g_i\circ\G_j^{-1}\}^{t\quad m}_{i=1, j=1}.$

Let $N(\varepsilon)$ be the normal subgroup of $\F_m$ generated by
$R(\varepsilon).$

Let $G(\varepsilon)=\df{\F_m}{N(\varepsilon)}.$
\edm

\proclaim {Lemma 0.1}
\roster\item"(i)" $N(\varepsilon) = \{(\a)\be\cdot\a^{-1} \bigm|
\a\in\F_m$,\quad $ \be\in M(\varepsilon)\}.$
\item"(ii)"  $G(\ve) = F_m$ with the relations induced from
$R(\varepsilon).$\endroster\endproclaim

\demo{Proof}\
Trivial. \qed\edm
\mk

We recall from BGT IV \cite{MoTe 7} the definitions of Hurwitz equivalent
factorizations and factorization invariant under $h\in B_m.$
\demo{Definition} \ $\underline{\text{Hurwitz move}}$

Let $g_i\cdot \dots \cdot g_k = h_i\cdot \dots \cdot h_k$ be two
factorized expressions of the same element in a group $G.$ \
We say that $g_i\cdot\dots \cdot g_k$ is obtained from
$h_i\cdot \dots \cdot h_k$ by a Hurwitz move if $\exists \ 1\leq p\leq k-1$
s.t. $g_i = h_i \ i\neq p, p+1,$ $g_p = h_p h_{p+1} h_p^{-1}$ and $g_{p+1} =
h_p$ or $g_p = h_{p+1}$ and $g_{p+1}=h_{p+1}^{-1} h_p h_{p+1}$
  \enddemo

\newpage

\demo{Definition}\ $\underline{\text{Hurwitz equivalence of factorized
expressions}}$

Let $g_i\cdot \dots \cdot g_k = h_i\cdot \dots \cdot h_k$ be two
factorized expressions of the same element in a group $G.$ \
We say that $g_i\cdot \dots \cdot g_k$ is a Hurwitz equivalent to
$h_i\cdot \dots \cdot h_k$ if $h_i\cdot \dots \cdot g_k$ is obtained
from $h_i\cdot \dots \cdot h_k$ by a finite number of Hurwitz moves.
We denote it by $g_1\cdot\dots\cdot g_k \underset{He}\to\simeq
h_i\cdot \dots \cdot h_k.$\enddemo

\demo{Definition} \ $\underline{\text{Factorized expression in}\ B_m\
\text{invariant under}\ h\in B_m}$

We say that a factorized expression $g_1\dots g_t$ is invariant under $h$ if
$(g_1)_h\cdot\dots\cdot (g_t)_h$ is Hurwitz equivalent to $g_1\dots g_t,$
i.e., can be obtained from $g_1\dots g_t$ by a finite number of Hurwitz moves
$((g)_h=h^{-1}gh).$  \edm

\proclaim{Lemma 0.2}\
If a factorized expression $\varepsilon=g_1\dots g_t$ in $B_m$ is invariant
under $h\in B_m,$ then $h$ induces an automorphism of
$G(\varepsilon).$\endproclaim

\demo{Proof}\
The group $B_m$ acts on $\pi_1 (D-K, u);$ thus there is a natural action of
$h\in B_m$ on $\pi_1 (D-k, u)= \F_m.$
Therefore, $h$ induces an automorphism of $\F_m.$
Since $\varepsilon$ is Hurwitz equivalent to $(\varepsilon)_h,$ we get that
$h^{-1} M(\varepsilon)
h=M(\varepsilon)$ and thus $(N(\varepsilon)) h\subset N(\varepsilon),$ and $h$
induces an automorphism of $\F/N(\varepsilon) = G(\varepsilon).$\quad
$\qed$\edm
\medskip
Certain factorized expressions of $\Delta_m^2$ in $B_m$ play an important role
in the computation of the fundamental group of complements of curves, as we
shall see in Theorem 0.3.

Let $S$ be a curve in $\Bbb C\Bbb P^2$ of degree $m.$

We refer the reader to  BGT I, Chapter VI, \cite{MoTe4} for  the definition of
a certain factorized expression  in $B_m$ related to $S$:\quad
 $\bold S${\bf-factorization of} $\bold{\Delta_m^2}$\quad   or
 $\bold S${\bf-factorization}\quad or \quad
  {\bf product form of} $\bold{\Delta_m^2}$ \quad or
 {\bf braid monodromy} {\bf factorization w.r.t.} $\bold S$ {\bf and}
$\bold u.$

\proclaim{Theorem 0.3.\ Zariski-Van Kampen Theorem} \rom{(see \cite{VK})}
 Let $\ov{S}$ be a   curve in $\CPt$ of $\deg m,$ s.t. $\overline S$
is transversal to the line in infinity.
 Let $S=\ov S\cap \C^2.$
 Let $\ve$ be the braid monodromy factorization w.r.t. to $S$ and $
u.$ Let $\C_u=u\times \C.$
 Let $\{\G_i\}$ be a $g$-base of $\pi_1(\C_u - S, u).$
Then:

$\pi_1(\C^2-S, *) = G(\varepsilon)$ \ and \
$\pi_1 (\CPt - \ov{S}, *) = G(\varepsilon)$ with extra relation\
$\prod\limits^m_{i=1} \G_i=1.$\endproclaim
\sk
\demo{Remark}\ We shall use this theorem for $S_3,$ the branch curve of a
generic projection of $V_3$ (the Veronese of order 3) to $\CPt.$
In BGT IV we computed the braid monodromy factorization related to $S_3.$
We denote it $\ve(18).$
We shall again present  $\ve(18)$ in the next chapter.
\edm
\sk

We are going to reformulate the Zariski-Van Kampen Theorem in a more precise
form for a cuspidal curve, i.e., for a curve with only nodes and cusps.

\proclaim{Theorem 0.4.\ (Zariski)}\ If $S$ is a cuspidal curve, then the
related
braid monodromy factorization $\ve$ is of the form $\prod\limits^p_{j=1}
V_j^{\nu_j},$ where $V_j$ is a half-twist and $\nu_j = 1$ or $2$ or $3$. \ep

\proclaim{Theorem 0.5.\ Zariski-Van Kampen (precise version)}\
Let $\ov{S}$ be a cuspidal curve in $\CPt.$
Let $S=\C^2\cap\ov{S}.$
Let $\ve$ be a braid monodromy factorization w.r.t. $S$ and $u.$
Let $\ve =
\prod\limits^p_{j=1} V_j^{\nu_j},$ where $V_j$ is a half-twist and $\nu_j =
1,2,3.$

For every $j=1\dots p$\ let $A_j, B_j\in\pi_1 (\C_u-S,u)$ be such that
 $A_j, B_j$ can be extended to a $g$-base of $\pi_1(\C_u-S,u)$
and $(A_j) V_j=B_j.$
Let $\{\G_i\}$ be a $g$-base of $\pi_1(\C_u-S, u).$
Then $\pi_1(\C^2-S,u)$ is generated by the images of $\{\G_i\}$ in
$\pi_1(\C^2 - S, u)$ and the only relations are those implied from $\left\{
V_j^{\nu_j}\right\},$ as follows:
$$\alignat 3
&A_j\cdot B_j^{-1} = 1 \qquad &&\text{if} \qquad &&\nu_j=1\\
&[A_j, B_j] = 1  &&\text{if}  &&\nu_j=2\\
&\la A_j, B_j\ra = 1  &&\text{if} &&\nu_j=3\endalignat$$

$\pi_1(\CPt-\ov{S}, *)$ is generated by $\{\G_i\}$ with the above relations and
one more relation $\prod\limits_i\G_i=1.$\endproclaim

\subheading{Remark 0.5$'$.\
How to determine $\bold{A_j}$  and $\bold{B_j}$   from
$\bold{V_j}$  or how to determine
 $\bold {A_V}$   and $\bold {B_V}$   from
$\bold{V=H(\sigma)}$  (see
formulation of Van Kampen Theorem)}

To be able to use the Zariski-Van Kampen Theorem, we must know how to compute
$A_j$ if $B_j$ for every $j=1\dots p.$
Assume, for simplicity, that $u_o$ is below real lines and $\{q_i\}=\C_u\cap
S$ are real points.
Assume that $\rho(\dl_j)=V_j^{\nu_j},$ where $V_j=H(\sigma)$, a half-twist
corresponding to a path $\s$ from $q_1$ to $q_2.$
Take a homotopically-equivalent path $\s'$ that passes through $u_0.$
Let $\s_1, \s_2$ be the part of $\s'$ from $u_0$ to $q_1, q_2$ respectively.
Let $A_j=\ell(\s_1)$ $B_j=\ell(\s_2)$ be the loops of $\pi_1(\C^2 - S,u_0)$
built from $\s_1,\s_2$ as in the definition in the beginning of the Chapter.
See Fig.  0.2 for an example how to determine $A_V$ and $B_V$ for
$V=H(\sigma).$

\midspace{2.00in}\caption{Fig. 0.2}

\proclaim{Proposition 0.6}\ If an $S$-factorization $\Dl^2=\prod g_i$
is invariant under $h\in B_m$ then $R((g_i)_h)$ is also a relation on
$\pi_1(\C^2-S,*).$\endproclaim

\demo{Proof}\   Zariski-Van Kampen Theorem and the fact that an invariant
factorization of a braid monodromy factorization, namely $\prod (g_i)_h,$ is
also a braid monodromy factorization (see Proposition  VI.4.2 from BGT I,
\cite{MoTe4}).\quad $\qed$\edm
\demo{Remark} Proposition 0.6 indicates  why it is important to prove
invariant properties of $S$-factorizations.
We use such properties to induce more relations on the fundamental
group.\enddemo \newpage

\heading{{\bf CHAPTER I.\ The braid monodromy related to $\bold S,$ the branch
curve of a $\bold{V_3}$-projection}}\endheading
\demo{Notation} \ $\un{V_3,\ f,\ S_3, S}$

Let $V_3$ be the Veronese surface of order 3,
i.e., the following embedding
of $\CPt$ into $\CP^N:$
$$  (x,y,z)\ri(\dots, x^iy^jz^k,\dots\ )_{i+j+k=3}.$$

Let $f=f_3$ be a generic projection:\ $V_3\overset{f}\to\ri \CPt.$

Let $S_3$ be its branch curve in $\CPt.$

Let $\C^2$ be a ``generic'' affine piece of $\CPt.$

Let $S=S_3\cap\C^2.$

$\deg S=\deg S_3=18.$

We are interested in $\pi_1(\CPt-S_3,*)$ and $\pi_1(\C^2-S, *).$

We constructed in BGT III \cite{MoTe7} a projective degeneration of $V_3
\overset{f}\to\ri \CPt$ into $Z^{(6)} \overset{f^{(6)}}\to\ri \CPt$ where
$Z^{(6)}$ is a union of $9$ planes $P_j,$\ $ j=1\dots 9,$ the ramification
curve
is a union of $9$ intersection lines $\hat{L}_i, i=1\dots 9,$ as in Fig. II.1.
(Each $\hat L_i$ is an intersection line of 2 $P_j$'s.)
$S^{(6)},$ the branch curve of $f^{(6)}$ in $\CPt,$ is a union of 9
lines $L_i, i=1\dots 9.$\ \
$(L_i=\pi^{(6)} (\hL_i)).$

$K^{(6)} = \C u\cap\Ss.$

$\#\ K^{(6)}=9.$

\edm

\demo{Definition} \ $\underline{u, \C_u, K.}$

Let us choose $u$ in the $x$-axis of $\C^2$ far away from the $x$-projection
singularities of the $x$-projection of $S_3$ and of $S^{(6)}.$

$\C_u = u\times\C.$

$K= \C_u\cap S.$
\edm

{}From the regeneration process it is obvious that for every point
$q_i$ that we had in $K^{(6)}$ we have   2 points $q_t, q_{t'}$ in $K$ which
are
close to each other. Recall that we used a ``real model'' of $(\C_u, K).$
Thus, we assume that $K=\{q_i, q_{i'}\}^9_{i=1},$\ $q_i, q_{i'}$ are real.

\demo{Remark}\ For arbitrary $n$ we would get that $V_n\ri\CPt$ degenerated
into
$Z \ri\CPt$ where $Z=$ a union of $n^2$ planes with a ramification curve which
consists of $\df{3}{2} n(n-1)$  intersection lines, and a branch curve
consisting of $\df{3}{2} n(n-1)$ lines.
Thus, $S_n,$  the branch curve of $V_n\ri\CPt,$ is of order $3 n(n-1).$
\edm

In BGT IV we computed a braid monodromy factorization of $S_3$ denoted
$\varepsilon(18).$

We also proved those invariance properties of $\varepsilon(18)$ and invariance
under complex conjugation.
We shall repeat  these results here.
For this we have to recall some notations.

\demo{Notations}
$$\alignat 3
&\un{z}_{ij}=z_{ij} &&= \text{a path below the real line from}\ q_i\ \text{to}
\ q_j. \qquad \qquad \qquad \ &&Z_{ij}=H(z_{ij})\\
&\ov{z}_{ij} &&= \text{a path above the real line from}\ q_i\ \text{to}
\ q_j.  &&\ov {Z}_{ij}=H(\ov{z}_{ij})\\
&\overset{(a)}\to
{\un{z}_{ij}}
&&= \text{a path above}\ a\ \text{and below the
real line elsewhere from}\ q_i\ \text{to} \ q_j.
&&\overset{(a)}\to{\un{Z}_{ij}}= H(\un{z}_{ij})\\
&\underset{(a)(b)}\to{\ov{z}_{ij}}
&&= \text{a path below}\ a\ \text{and}\ b\ \text{above the real line elsewhere
from}\ q_i\ \text{to} \ q_j.  &&\underset{(a)(b)}\to{\ov{Z}_{ij}}=
H \underset{(a)(b)}\to {(\ov{z}_{ij})}\\
&\rho_i&&=Z_{ii'} =H(z_{ii'}) = \text{half-twist corresponds to the shortest
line}&&\\
& &&\quad\text{ between}\ q_i\ \text{and}\ q_{i'}&&\\
&Y^{(2)}_{i,jj'}&&=\prod\limits_{m=0}^1 (Y_{ij}^2)\rho_j^m&&\\
&Y^{(2)}_{ii',jj'}&&=\prod\limits_{\ell=0}^1\prod\limits_{m=0}^1
(Y_{ij}^2)\rho_i^\ell\rho_j^m&&\\
&Y^{(3)}_{i,jj'}&&=\prod\limits_{m=-1}^{1} (Y_{ij}^3)\rho_j^m&&  \endalignat
$$\edm

\flushpar For $\be>\g$

$\tilde{\un{z}}_{\g\g'(\be)} =$

\midspace{.75in}

$\un{\tilde{z}}_{\be\be'(\g)}=$

\mk

With the above notations we recall the braid monodromy factorization of
$S=S^{(0)},$ denoted $\ve(18)$ and some invariance properties.

\proclaim{Theorem I.1}

$\varepsilon(18) =\prod\limits^1_{\nu=7}\ C_\nu\
H_\nu$ \qquad
\flushpar where $H_i,$  the  braid monodromy of $S $ around $v_i$ factors as
follows:

$H_1 = Z^{(3)}_{11',2} \cdot \tilde Z_{22'(1)}$

$H_2 = Z^{(3)}_{11',3} \cdot \tilde Z_{33'(1)}$

$H_3 = Z^{(3)}_{44',6} \cdot \tilde Z_{66'(4)}$

$H_5 = Z^{(3)}_{55',9} \cdot \tilde Z_{99'(5)}$

$H_6 = Z^{(3)}_{6', 77'} \cdot \tilde Z_{66'(7)}$

$H_7 = Z^{(3)}_{8', 99'} \cdot \tilde Z_{88'(9)}$

$H_4=Z_{2',33'}^{(3)}\tilde
Z_{88'}Z_{44',8'}^{(2)}\left(Z_{33',8}^{(2)}\right)^\bullet\overline
Z_{55',8}^{(3)} \left(Z_{44',8}^{(2)}\right)^\bullet
\left(Z_{33',8}^{(2)}\right)^\bullet\hat F_1(\hat F_1)_{\rho^{-1}}$

$\qquad \cdot Z_{77',8}^{(3)} \ Z_{2',i}^2, i=8,8,7',7,5',5\quad\overline
Z_{2',44'}^{(3)} \ Z_{2i}^2 \ \ i=8',8,7',7,5',5\quad\tilde
Z_{22'}$

\flushpar where
$\tilde Z_{22'},$\ $\tilde Z_{88'}$ correspond to the  paths $\tz_{22'},$
$\tz_{88'}$
described in Fig. \rom{I.0(a)}.

\midspace{1.00in}\caption{Fig. I.0\text{\rm(a)}}

$\bullet$ denotes conjugation by a braid $b^\bullet$ induced from the motion
described in Fig. \rom{I.0(b)}.

\midspace{1.00in}\caption{Fig. I.0\text{\rm(b)}}

(In fact,\ $b^\bullet = Z_{2'3'}^2Z_{2'3'}^2Z_{7'8}\2Z_{78}\2.)$
$$\align
&\hat F_1 = Z^{(3)}_{3,44'} Z^{(3)}_{55',7} \alpha^{(1)}
\overset {(4)}\to Z^2_{3'7} \overline Z^2_{3'7'} \\
&\hat F_2 = \left(Z^{(3)}_{3,44'}\right)\rho^{-1} \left(Z^{(3)}_{55',7}
\right)_{\rho^{-1}} \left(\alpha^{(1)} \right)_{\rho^{-1}}
\left(\overset {(4)}\to {\underline {Z}}^2_{3'7} \right)_{\rho^{-1}}
\left( \overline Z_{3'7'}\right)^2 \rho^{-1}.\endalign$$\mk

$\a^{(1)} = \a_1\a_2$

\flushpar where $\a_1$ and $\a_2$ are the curves described in Fig.
\rom{I.0(c)}, \rom{I.0(d)}, respectively.

\midspace{1.25in}\caption{Fig. I.0\text{\rm(c)}}

\midspace{1.25in}\caption{Fig. I.0\text{\rm(d)}}

$\rho=\rho_7\rho_3$

 $(\a^{(1)})\rho = (\a^{(1)})\rho^{-1}_3\rho^{-1}_7 = (\a_1)(\rho^{-1}_3
\rho^{-1}_7)\cdot(\a_2)(\rho_3\rho_7)^{-1}$

\flushpar and
$$\align
C_1 &= C_2 = I_d\\
C_3 &= Z^{(2)}_{11', 44'} \prod_{i=1,2,3,5} Z^{(2)}_{ii', 66'}\\
C_4 &= \ov{Z}^{(2)}_{11', 55'} \ \underset{(6)}\to{\ov{Z}}^{(2)}_{11', 77'}\
\ov{Z}^{(2)}_{11', 88'}\ \ov{Z}^{(2)}_{66', 88'}\\
C_5 &= \prod^7 \Sb i=1\\ i\neq 5\endSb Z^{(2)}_{ii', 99'}\\
C_6 &= C_7 = I_d.\endalign$$\ep

The following remark gives an explicit description of a half-twist conjugated
by some $\rho_j$ and will help us later to deduce relations from $\ve(18)$
using the Van Kampen method.
\subheading{Remark I.1}
\roster\item"(i)" $Z_{i,jj'}^{(3)}=\overset (j)\to{\un Z}_{ij'}^3\un
Z_{ij}^3\un Z_{ij'}^3$\quad (Fig. I.1(a))
\item"(ii)" $Z_{ii',jj'}^{(2)}=Z_{ij}^2Z_{ij'}^2Z_{i'j}^2Z_{i'j'}^2.$\quad
(Fig.
I.1(b)) \item"(iii)" Let $Y_{ij}=H(y_{ij})$ where $y_{ij}$ is a path connecting
$q_i$ or $q_{i'}$ with $q_j$ or $q_{j'}.$
The following graph (see Fig. I.1(c)) indicates the conjugation of
$Y_{ij}$ by $\rho_j,\rho_j\1,\rho_i,\rho_i\1$ for
different types of $y_{ij}.$
In the graph we only indicate the action of $\rho_j$ and $\rho_j\1$ on the
``head'' of $y_{ij}$ within a small circle around $q_j$ and $q_{j'}$ and the
action of $\rho_i$ and $\rho_i\1$ on the ``tail'' of $y_{ij}$ in a small
circle around $q_i$ and $q_{i'}.$
The ``body'' of $y_{ij}$ is not changing under $\rho_j^{\pm 1}$ and
$\rho_i^{\pm 1}.$ \endroster

\newpage

\midspace{8.00in}\caption{Fig. I.1}

\newpage
\proclaim{Theorem I.2} \ \text{\rm Invariance Theorem}\ \text{\rm(BGT IV,
Proposition 18, \cite{MoTe7})}

Let $\rho=\rho_{m_1\dots
m_4,m_6,m_9}=\rho_1^{m_1}(\rho_2\rho_8)^{m_2}\cdot(\rho_3\rho_7)^{m_3}\cdot
(\rho_4\rho_5)^{m_4}\cdot \rho_6^{m_6}\cdot \rho_9^{m_9}$
then, $\varepsilon(18)$ is invariant under $\rho$ for every $m_i\in \Bbb Z.$\
$\rho_i = Z_{ii'}).$   \endproclaim

\medskip

\proclaim{Theorem I.3} \text{\rm Complex Conjugation Theorem}\ \text{\rm (BGT
IV, Proposition 19, \cite{MoTe7})}

$\varepsilon(18)$ is invariant under complex conjugation.\endproclaim

\subheading{A finite set of generators for $\bold{\pi_1 (\C^2 - S, u_0)}$ and
$\bold{\pi_1 (\CPt - S_3, u_0)}$}\

Let us choose $u_0\in\C_u, u_0$ below the real line.
Let $\{\G_i, \G_{i'}\}$ be a $g$-base of\linebreak
$\pi_1(\C_u-S, u_o).$
When considered as elements of $\pi_1(\C^2 - S, u_0)$ and of $\pi_1(\CPt -
S_3, u_o),$ they generate (not freely) the groups.
Thus we have $\{\G_i,\G_{i'}\}_{i=1}^9,$ a set of generators for
$\pi_1(\C^2-S,u_0)$ and $\pi_1(\CP^2-S_3,u_0).$

We want to compute $G=\pi_1(\C^2-S, u_0)$ and $\ov{G}=\pi_1(\CPt - S_3, u_0).$

By the Zariski-Van Kampen Theorem, $G\simeq G(\varepsilon(18))$ and $\ov{G}
\cong \df{G(\varepsilon(18))}{\prod\limits^9_{i=1}\G_i\G_{i'}}.$

\proclaim{Corollary I.4}\
$G=\pi_1 (\C^2 - S, u_0)$ satisfies all the relations induced in
$R(\varepsilon(18)),$ all the relations induced from $(R(\varepsilon(18))
\rho_{m_{1},\dots,m_{4}, m_{6}, m_{9}}$   and all the relations induced from
the complex conjugation of $\varepsilon(18),$ where $\rho_{m_1\dots
m_4,m_6,m_9}=\rho_1^{m_1}(\rho_2\rho_8)^{m_2}\cdot(\rho_3\rho_7)^{m_3}\cdot
(\rho_4\rho_5)^{m_4}\cdot \rho_6^{m_6}\cdot \rho_9^{m_9}$
 and $\rho_i = Z_{ii'}.$
Moreover,
$$\align
&G\sim G (\varepsilon(18))\\
&\ov{G} = G/\prod_i \G_i\G_{i'}.\endalign$$
\endproclaim

\demo{Proof}\ Theorems I.1, I.2, I.3 and the Van Kampen Theorem.\ \ \
$\square$ \edm

\proclaim{Corollary I.5}\ Let $G=\pi_1 (\Bbb C^2 - S, u_0).$
If $R$ is any relation
in $G$  then $(R)\rho_
{m_{1},\dots,  m_4,m_6,m_9}$ is also a relation in $G,$ where
$(R)\rho_{m_{1},\dots, m_4,m_6, m_9}$ is the relation induced from $R$ by
replacing  $\G_i$ and $\G_{i'}$ with $(\G_i)\rho_i^{m_i}$ and $
(\G_{i'}) \rho_i^{m_{i}},\ i= 1,2,3,4,6,9,$ respectively and replacing
$\G_8, \G_{8'}, \G_7, \G_{7'}, \G_5, \G_{5'}$ with
$(\G_8)\rho_8^{m_2}, (\G_{8'})\rho_{8}^{m_2},
(\G_7)\rho_7^{m_3}, (\G_{7'})\rho_{7}^{m_3},
(\G_5)\rho_5^{m_4}, (\G_{5'})\rho_{5}^{m_4},$\linebreak
respectively.\endproclaim

\demo{Proof}\ Proposition 0.6 and Theorem I.2.\ \ \ $\square$\edm

\demo{Remark}\
In other words, $\rho_{m_1\dots m_4,m_6,m_9}$ defines an automorphism of
$G.$\edm

\demo{Notation}

$\G_{\un{i}} =$ any element of the set $\left\{
(\G_i)\rho_i^{m_i}\right\}_m\in\Bbb Z.$\edm
\mk
\proclaim{Corollary I.6}\ \text{\rm(Complete invariance in $3$-points)}

Let $(\a,\be)=(1,2)$ or $(1,3)$ or $(4,6)$ or $(5,9)$ or $(6,7$) or $(8,9).$
Any relation $R$ in $G$ that involves only $\G_\a$ and $\G_\be$ is true when
$\G_{\un{\a}}$ replaces $\G_\a$ amd $\G_{\un{\be}}$ replaces
$\G_\be.$\endproclaim

\demo{Proof}\ Without loss of generality assume $(\a, \be)=(1,2).$
Let $\G_{\un{1}}=(\G_1)\rho^{m_1}_1$ $\G_{\un{2}}=(\G_2)\rho^{m_2}_2.$
By Corollary I.5, $(R) \rho_{m_{1},m_{2},0,0,0,0}$ is also a relation on $R,$
and differs from $R$ by replacing
$\G_{\un{1}}$ by  $\G_1$ and $\G_{\un{2}}$ by $\G_2.$\ \ \ \qed\edm

\proclaim{Corollary I.7}\ Let $R$ be a relation in $G$ that involves at
most one index of each of the pairs $(2,8),$ $(3,7)$, $(4, 5).$ For every
index $i$ that appears in $R$ choose some $\G_{\un i}.$ Then $R$ is true when
we simultaneously replace $\G_i$  by  $\G_{\un{i}}$ (or $E_i$) and
$\G_i'$   by $(\G_{\un i})\rho_i$ (or $(E_{\un i})\rho_i,$
respectively), for $i\in\{$indices that appear in $R\}.$\endproclaim

\demo{Proof}\  We can assume w.l.o.g. that $\G_8,\G_{8'}, \G_7\,\G_{7'},
\G_5,\G_{5'}$ appear in $R$ and\newline
$\G_2,\G_{2'},\G_3,\G_{3'},\G_4,\G_{4'}$
do not appear in $R.$ Let $i$ be s.t. $\G_i$ or $\G_{i'}$ appears in
\newline $R$\ $(i\in 1,5,6,7,8,9).$ There exist
$m_i$  s.t. $\G_{\un{i}} = (\G_i)\rho_i^{m_i}.$ \newline Let
$\rho=(\prod_{i=1,5,6,8,9}\rho_i^{m_i})\rho_2^{m_8}\rho_3^{m_7}\rho_4^{m_5}.$
By Corollary I.5 $(R)\rho$ is also
a relation in $G.$
Since $\G_2,\G_{2'}, \G_3,\G_{3'},
\G_4,\G_{4'}$ do not appear in $R$ the relation $(R)\rho$ is actually equal to
$(R)\prod_{i=1,5,6,8,9}\rho_i^{m_i}$ and it   differs from $R$ by
replacing $\G_i$ with $(\G_i) \rho_i^{m_i}$ and $\G_{i'}$ by
$(\G_{i'})\rho_i^{m_i}$ for $i=1,5,6,7,8,9.$
But $(\G_i)\rho_i^{m_i}=\G_{\un i}$ and
$(\G_{i'})\rho_i^{m_i}=(\G_i)\rho_i\rho_i^{m_i}=(\G_i)\rho_i^{m_i}\rho_i=(\G_{\un
i})\rho_i,$ so we get the corollary.\ \ \ \qed\edm
\demo{Remark} We can replace every $\G_i$ that appear in $R$ by any $\G_{\un
i}=(\G_i)\rho_i^{m_i}$ (i.e., different $m_i$'s for different $\G_i$'s) since
both $\G_3$ and $\G_7$\ $(\G_2$ and $\G_8,$\ $\G_5$ and $\G_4,$ respectively)
do not appear in $R.$
If both $\G_3$ and $\G_7$ appear in $R,$ then we could only replace $\G_3$ by
$(\G_3)\rho_3^m$ and $\G_7$ by $(\G_7)\rho_7^m$ for the same $m.$\edm

\newpage

\heading{{\bf CHAPTER II.\ List  of relations in $\bold G$}}\endheading

We are going to describe $G$ using a different set of generators than those
introduced in Chapter I. We use all the notations from Chapter I, all the
relations induced from the braid monodromy factorization $\varepsilon(18)$
(Theorem I.1), the complex conjugation (Theorem ~I.3), and Corollaries I.5 and
I.6.
\subheading{Remark II.0}\ $\un{\text{First set of generators for}\ G}$

Let $\{\G_i,\G_{i'}\}$ be a $g$-base of $\pi_1(\Bbb C_u-S,u_0).$
Considered. as elements of \newline $G=\pi_1(\Bbb C^2-S,u_0),\
\{\G_i,\G_{i'}\}$
generates $G.$ \demo{Definition}
$$ E_i=\cases
\G_i \qquad &i\neq 2,7\\
\G_{i'} &i=2,7\endcases$$
$\qquad \qquad\qquad\qquad\qquad\qquad\quad\ \ E_{i'}=(E_i)\rho_i.$\edm

\demo{Notations}

$\rho_i=Z_{ii'}$ the half-twist corresponding to the shortest path between
$q_j$ and $q_{j'}.$

 $E_{\un{i}} = $ an element of $\{(E_i)\rho^m_i\}_m\in\Z.$

$\G_{\un{i}} =$ an element of $\{(\G_i)\rho^m_i\}_m\in\Z.$\edm

\medskip
In order to use the invariance theorem we need the following lemma.

\proclaim{Lemma II.0}
\roster
\item"(i)"\ $(\G_j)\rho_i=\G_j\quad (\G_{j'}) \rho_i = \G_{j'}$\ \text{for}\
$i\neq j.$
\item"(ii)"\ $
(\G_i)\rho_i=\G_{i'} \quad (\G_{i'})\rho_i=\G_{i'}\G_i
\G_{i'}^{-1} \quad (\G_i)\rho_i^{-1}
 = \G_i^{-1}\G_{i'}\G_i.$
\item"(iii)"\ Let $\rho = \rho_{m_1,\dots, m_4, m_6,
m_9}=\rho_1^{m_1}(\rho_2\rho_8)^{m_2}\cdot(\rho_3\rho_7)^{m_3}\cdot
(\rho_4\rho_5)^{m_4}\cdot \rho_6^{m_6}\cdot \rho_9^{m_9}.$ Then $(\G_i)\rho =
\G_{\un{i}}.$
\item"(iv)"\ $\G_{\un{i}}\in \la \G_i,  \G_{i'}\ra.$\endroster\endproclaim

\demo{Proof}\ Geometric observation (Fig. II.0(a) and (b)) or BGT I, Section
II,
\S 2. \qed\edm

\newpage

\midspace{8.00in}\caption{Fig. II.0}

\newpage

\proclaim{Lemma II.1}
\roster
\item"(i)"\ $E_{i'}=\cases
\G_{i'}  &i\neq 2,7\\
\G_{i'}\G_i \G_{{i}'}^{-1}
&i= 2,7\endcases$
\item""
\item"(ii)"\ $\G_{i'}\G_i= E_{i'} E_i$
\item"(iii)" $(E_i)\rho_i=E_{i'}$
\item"(iv)" $(E_{i'})\rho_i=E_{i'}E_iE_{i'}^{-1}$
\item""
\item"(v)"\ $\G_i=\cases
E_i \qquad   &i\neq 2,7\\
E_{i}^{-1} E_{i'} E_i &i=2,7\quad ( =
(E_{i'})\rho_i^{-2}=(E_i)\rho^{-1}_i).\endcases$\endroster\endproclaim
\demo{Proof} Trivial.\qed\edm
\proclaim{Lemma II.2}\ $\{E_i, E_{i'}\}$ generate $G.$\endproclaim

\demo{Proof}\ Trivial. \qed\edm
\smallskip
\subheading{Remark II.2}\ $\un{\text{A second set of generators for}\ G}$

We start with a set of
generators $\G_i,\G_{i'}$ and exchange it for a set of generators
$E_i,E_{i'}.$
 \demo{Definition}\
Let $\psi$ be the classical monodromy homomorphism  from $G$ to the
symmetric group of order 9 induced by the projection $V_3\ri\CPt.$\edm

\medskip

\proclaim{Lemma II.3}
$\psi(\G_i)=\psi(\G_{i'}) = \psi(E_i)=\psi(E_{i'})=(k_i, \ell_i)$ \ where
\
$\hat{L}_i= P_{k_i}\cap P_{\ell_i},$ and $\{L_i\}$ and $\{P_j\}$
are arranged as in Fig. \rom{II.1}.\endproclaim
\midspace{3.00in}\caption{Fig. II.1}

\demo{Proof}\ Let $\g_i$ be a path in $\Bbb C^2$ connecting $u_0$ to $q_i$
s.t. $\G_i=\ell(\g_i)$ (see the definition of $\ell(\g)$ in Chapter 0). One has
to consider the degeneration of $V_3\overset{f}\to\ri\CPt$ to $Z^{(6)}\overset
f^{(6)} \to\ri \CPt$ and of $S_3$ to  $S^{(6)},$
constructed in  BGT III, \cite{MoTe7}.
The surface $Z^{(6)}$ is a union of 9
planes, $P_1\cdots P_9.$
The configuration of the planes and their intersection lines $\hat
L_1,\dots,\hat L_9$ are as in Fig. II.1.
 Let $f^{(6)}$ be a generic projection
$Z^{(6)}\overset{f^{(6)}}\to\ri\CPt$ and $ S^{(6)}$ its branch curve in $
\C^2.$
We choose $V_3$ to be close to $Z^{(6)}.$ Let $p_i^{(6)}= P_i\cap
\pi^{(6)^{-1}}
(u_0).$
Let $p_i$ be a point in $\pi^{-1}(u_o)$ which is close to $p_i^{(6)}.$
Fix $i$ between 1 and 9.
It is clear that when we move along $\g_i$ from $u_0$ to $q_i,$ the lifted path
in $Z^{(6)}$ which starts in $p^{(6)}_{k_i}$ will lie in $P_{k_i}$ and will
end on a point in $\hat{L}_i$ above $q_i.$
The lifted path in $Z^{(6)}$ that starts in $p^{(6)}_{\ell_i}$ will lie in
$P_{\ell_i}$ and will end in the same point in $\hat{L}_i.$
Thus, in the regenerated case,  the lifted path of $\G_k$ that starts
in $p_i$ will end in $p_j.$
The lifted paths of $\g_i$ in $Z^{(6)}$ that start in $p^{(6)}_t,$\ $t\neq
k_i,\ell_i$ will be closed loops.
Thus, in the regenerated case, the path obtained from lifting $\G_k$ that
starts in  $p_t,$\ $t\neq k_i, \ell_i$ is a loop.
Thus, $\psi(\G_k)=$ the transposition $(k_i\ \ \ell_i)$ of the symmetric group
on 9 elements In the same way $\psi(\G_{i'})=(k_i\ \ \ell_i),$ and thus
$\psi(E_i)=\psi(E_{i'}) =(k_i\ \ \ell_i).$\quad\qed \edm
\proclaim{Corollary II.3}\
The transpositions $\psi(E_i)$ are as follows:

$\psi(E_1)= (1\quad 2)$

$\psi(E_2)=(2\quad 3)$

$\psi(E_3)=(2\quad 4)$

$\psi(E_4)= (3\quad 5)$

$\psi(E_5)=(4\quad 7)$

$\psi(E_6)=(5\quad 6)$

$\psi(E_7)= (5\quad 8)$

$\psi(E_8)=(7\quad 8)$

$\psi(E_9)=(7\quad 9)$

Moreover, $\psi(E_i)$ and $\psi(E_j)$ have one common index
$\Leftrightarrow \hat L_i$ and $\hat L_j$ are edges of some triangle in Fig.
\rom{II.1}. \endproclaim

\demo{Proof}\ Immediate from the previous Lemma.\qed\edm

\mk

Before continuing with $G,$ we want to prove some claims concerning an
arbitrary
group.

\proclaim{Claim II.4}

\flushpar \text{\rm(a)}\ If $\la A,B\ra = 1,$ i.e., $ABAB\1A\1B\1 = 1,$ then
$$\align
ABA &= BAB\\
A^{-1}B^k A &= BA^kB^{-1}\\
AB^kA^{-1} &= B^{-1}A^kB\\
A^{-1}B^{-1} A^k &= B^kA^{-1}B^{-1}.\endalign$$

\flushpar  \text{\rm(b)} \ $[A,B]=[A,C]=1\Rightarrow [A,D] = 1\ \text{for}\
D\in\la B,C\ra =\text{subgroup generated by}\  B,C.$

\flushpar  \text{\rm(c)} \ $[A_C, B_D] = 1,\ [C,B] =
[C,D] = [A,D] = 1 \Rightarrow [A,B] = 1.$

\flushpar  \text{\rm(d)} \ $[A,XYZ]=[A,X]\ [A,Y]_{X^{-1}}[A,Z]_{Y^{-1}X^{-1}}.$

\flushpar  \text{\rm(e)} \ $[XY,A]=[Y,A]_{X^{-1}}[X,A].$

\flushpar  \text{\rm(f)} \ $[A_C,B_C]=[A,B]_C.$

\flushpar  \text{\rm(g)} \ $[X,Z]=1\Rightarrow [X,Y_Z]=[X,Y]_Z.$

\flushpar  \text{\rm(h)} \ If $\la x,y\ra=1,$ then $\la
Ax,y\ra=1\Leftrightarrow A_{y\1x\1}=A\1A_{y\1}.$ A\endproclaim
\newpage

\demo{Proof} (a) -- (g) are easy to verify.

We shall only prove (h) here.

If $\la Ax,y\ra=1, $ then:
$$\alignat 4
& && \ \ &&1 = \la A x, y\ra &&= A x y x y\1 x\1 A\1 y\1  \\
& && && &&= A\cd A_{y\1 x\1} x y x y\1 x\1  A\1 y\1 y A\1 y\1\\
& && && &&= A\cd A_{y\1 x\1} \cd 1\cd A\1_{y\1}\\
& &&\Rightarrow && A_{y\1 x\1} &&= A\1 A_{y\1}.\endalignat$$
If $
A\1 A_{y\1} = A_{y\1x\1},$ then $$\alignat2
& &&A A_{y\1x\1} A\1_{y\1} = 1 \\
&\Rightarrow &&Axy \underbrace {Ay\1 x\1 y} A\1 y\1 = 1 \\
& \Rightarrow &&Axy\overbrace {Ax y\1 x\1}\ A\1 y\1 = 1  \\
&\Rightarrow &&Ax\cd y \cd Ax \cd y\1\cd (Ax)\1\cd y\1 = 1 \\
&\Rightarrow &&\la Ax, y\ra = 1.\endalignat$$
\edm\qed

\proclaim{Lemma II.5}

\rom{(a)}\ $[\G_i,\G_j]=[\G_{i'},\G_j] = 1 \Rightarrow [\G_{\un{i}},
\G_{\un{j}}]=1.$

\rom{(b)}\ $[\G_i,\G_j]=[\G_i^{-1} \G_{i'}\G_i, \G_j] = 1 \Rightarrow
[\G_{\un{i}}, \G_j]=1.$\endproclaim

\smallskip

\demo{Proof}\ We only prove (a);
(b) is the same argument.
Since $\G_{\un{j}}$ is a product of $\G_j$ and $\G_{j'},$ we can apply Claim
II.4 (b) to get $[\G_i, \G_{\un{j}}]=1.$
In particular, $[\G_i,(\G_j)\rho_j^{-1}]=1.$
We use invariance under $\rho_i\rho_j$ and Lemma 0.2 on it to get
$[\G_{i'},\G_j]=1.$
We use invariance under $\rho_i\rho_j$ and Lemma 0.2 on $[\G_i, \G_j]=1$ to
get    $[\G_{i'},\G_{j'}]=1.$
{}From $[\G_{i'}, \G_{j'}] = [\G_{i'}, \G_j] = 1$ we get, using Lemma II.4(b),
that $[\G_{i'}, \G_{\un{j}}] = 1.$
{}From $[\G_i, \G_{\un{j}}] = [\G_{i'}, \G_{\un j}] = 1$ we get, using Lemma
II.4(b), that
$[\G_{\un{i}}, \G_{\un{j}}] = 1.$\qquad \qed \edm

\newpage

\proclaim{Proposition II.6}\ The following relations hold in $G:$
\roster\item $\la E_{\un{i}},E_{\un{j}}\ra=1\ \forall i, j\ \text{s.t.}\
\psi(E_i)\  \text{and}\ \psi(E_j)$ have exactly one common
index.
\item $[ E_{\un{i}},E_{\un{j}}]=1\
\forall i, j\ \text{s.t.}\ \psi(E_i)$\ and $\psi(E_j)$\  have no
common index.
\item $1=(E_7E_5E_3^{-1}E_4^{-1}E_2E_4E_3E_5^{-1}E_7^{-1}E_8^{-1})
(\rho_3\rho_7)^i (\rho_4\rho_5)^j\quad \forall i, j\in\Z.$
\item
$E_{4'}E_4E_{3'}E_3E_{2'}E_3^{-1}E_{3'}^{-1}E_4^{-1}E_{4'}^{-1}=E_2.$
\item $E_5^{-1}E_{5'}^{-1}E_7^{-1}E_{7'}^{-1}E_8E_{7'}E_7E_{5'}E_5=E_{8'}.
$ \item
$$\align
					E_{\be'}=E_{\a}^{-1}E_{\a'}^{-1}E_\be E_{\a'}E_\a\qquad (\a,\be) &=(1,2)\\
	&= (1,3)\qquad\qquad\qquad\qquad\qquad\qquad\\
&= (4,6)\\
&= (5,9).\endalign$$
\item$$\align
					E_{\a}=E_{\be'}E_{\be} E_{\a'}E_{\be}^{-1}E_{\be'}^{-1}\quad\quad (\a,\be)
&=(6,7) \qquad\qquad\qquad\qquad\qquad\qquad \\
&=(8,9)\endalign$$ \endroster
Definitions of $\psi(E_i)$ as in Corollary \rom{II.3}.\endproclaim

\smallskip

\demo{Proof}
We divide the proof into the following 48 claims:\enddemo
\demo{Claim 0}\ $[\G_{\un i},\G_{\un j}]=1$  for $i,j$ s.t.
$\hat{L}_i\cap\hat{L}_j$ do not intersect, $j\neq 9,$  i.e., for
$(i,j)=(1,4),(1,5), (1,6), (1,7), (1,8), (2.6),
(3.6), (6.8).$ \edm
\demo{Claim 1}\ $\la\G_{\un i},\G_{\un j}\ra=1\ (i,j)=(1,2), (1,3), (4,6),
(5,9),
(6,7), (8,9)$\edm \demo{Claim 2}\ $\la\G_{2'},\G_{\un{3}}\ra=1.$\edm
\demo{Claim 3}\
$\G_{3'}^\bullet=\G_3\G_{2'}\G_{3}^{-1}=\G_{2'}^{\1}\G_3\G_{2'}.$\edm
\demo{Claim 4}\
$\G_{3}^\bullet=\G_{2'}^{\1}\G_{3'}^{\1}\G_{3'}\G_3\G_{2'}.$\edm
\demo{Claim 5}\ $\la\G_{\un{7}},\G_8\ra=1.$\edm
\demo{Claim 6}\
$\G_{7'}^\bullet=\G_8\G_{7'}\G_7\G_{7'}^{\1}\G_{8'}^{\1}.$\edm
\demo{Claim 7}\
$\G_{7}^\bullet=\G_8\G_{7'}\G_{8'}^{\1}=\G_{7'}^{\1}\G_8\G_{7'}.$\edm
\demo{Claim 8}\
$\G_{2'}^\bullet=\G_{3'}\G_3\G_{2'}\G_{3}^{\1}\G_{3'}^{\1}$

$\qquad\quad\
\G_{8}^\bullet=\G_{8}^{\1}\G_{7'}^{\1}\G_8\G_{7'}\G_{7}.$\edm
\demo{Claim 9}\ $\G_{2}^\bullet=\G_2, \G_{8'}^\bullet=\G_{8'},
\G_{4}^\bullet=\G_{4},\G_{4'}^\bullet=\G_{4'},\G_{5'}^\bullet=\G_{5'},
\G_{5}^\bullet=\G_{5}.$\edm
\demo{Claim 10}\ $[\G_{\un{2}},\G_{\un{i}}]=1\qquad i=5,7,8.$\edm
\demo{Claim 11}\ $[\G_{8'},\G_i]=1 \qquad i = 3,4.$\edm
\demo{Claim 12}\ $[\G_{\un{8}},\G_{\un{i}}]=1\qquad i=3,4.$\edm
$$ \alignat2 \text{{\it Claim 13.}}\quad
E_{\be'}=E_{\a}^{-1}E_{\a'}^{-1}E_\be E_{\a'}E_\a\qquad
\qquad &(\a,\be) &&=(1,2) \qquad\qquad\qquad\qquad\qquad\qquad\qquad\qquad\\
& &&= (1,3)\\
& &&= (4,6)\\
& &&= (5,9)\\
E_{\a}=E_{\be'}E_{\be} E_{\a'}E^{-1}_{\be}E_{\be'}^{-1}\qquad\qquad
&(\a,\be) &&=(6,7) \\
& &&=(8,9)\endalignat$$
\demo{Claim 14}\ $[\G_{3'}^\bullet,\G_{7'}^\bullet]=1.$\edm
\demo{Claim 15}\ $[\G_{3},\G_{7'}\G_7\G_{7'}^{-1}]=1.$\edm
\demo{Claim 16}\ $[\G_{3'},\G_{7'}\G_7\G_{7'}^{-1}]=1.$\edm
\demo{Claim 17}\ $[\G_{\un{3}},\G_{7'}\G_7\G_{7'}^{-1}]=1.$\edm
\demo{Claim 18}\ $[\G_{\un{3}},\G_{\un{7}}]=1.$\edm
\demo{Claim 19}\ $[\G_{\un{3}}^\bullet,\G_{\un{7}}^\bullet]=1.$\edm
\demo{Claim 20}\ $\la\G_{\un{3}}^\bullet,\G_{\un{4}}^\bullet\ra=1.$\edm
\demo{Claim 21}\ $\la\G_{\un{5}}^\bullet,\G_{\un{7}}^\bullet\ra=1.$\edm
\demo{Claim 22}\ $\G_7^\bullet \G_5^\bullet\G_{7}^{\bullet^{-1}}=
\G_{3'}^{\bullet^{-1}}\G_{4}^{\bullet}\G_{3'}^\bullet\ (\rho_3\rho_7)^i
(\rho_4\rho_5)^j.$\edm
\demo{Claim 23}\ $[\G_7^\bullet,
\G_{4'}^\bullet\G_4^\bullet\G_{3'}^\bullet\G_4^{\bullet^{-1}}\G_{4'}^{\bullet^{-1}}]=1.$\edm
\demo{Claim 24}\ $\la\G_{\un{3}}^\bullet,\G_{\un{5}}^\bullet\ra=1.$\edm
\demo{Claim 25}\ $\la\G_{\un{4}}^\bullet,\G_{\un{7}}^\bullet\ra=1.$\edm
\demo{Claim 26}\ $\la\G_{\un{3}},\G_{\un{5}}\ra=1.$\edm
\demo{Claim 27}\ $\la\G_{\un{7}},\G_{\un{4}}\ra=1.$\edm
\demo{Claim 28}\ $[\G_{4'},\G_5]=1.$\edm
\demo{Claim 29}\ $[\G_{4},\G_{5}]=1 $\edm
\demo{Claim 30}\ $[\G_{\un{4}},\G_{5}]=1.$\edm
\demo{Claim 31}\ $[\G_{\un{4}},\G_{\un{5}}]=1.$\edm
\demo{Claim 32}\ $[\G_{2'},\G_{7'}^\bullet]=1.$\edm
\demo{Claim 33}\ $[\G_{3},\G_{7}^\bullet]=1.$\edm
\demo{Claim 34}\ $[\G_{7}^\bullet,\G_{5}^2\G_3\G_5^{-2}]=1.$\edm
\demo{Claim 35}\ $[\G_{2'},\G_{3}^2\G_5\G_3^{-2}]=1.$\edm
\demo{Claim 36}\ $[\G_{\un{3}},\G_{\un{4}}]=1.$\edm
\demo{Claim 37}\
$[\G_{7'}^\bullet,\G_4^2\G_{3'}^\bullet\G_4^{-2}]=1.$\edm
\demo{Claim 38}\ $[\G_{7'}, \G_{3'}^\bullet]=1, \qquad
[\G_{8},\G_{3'}^\bullet]= 1.$\edm
\demo{Claim 39}\ $[\G_4,\G_{7'}^2\G_8\G_{7'}^{-2}]=1.$\edm
\demo{Claim 40}\ $[\G_5, \G_{7'}]=1.$\edm
\demo{Claim 41}\
$[\G_{\un{5}},\G_{\un{7}}]=1.$\edm \demo{Claim 42}\
$\la\G_{\un{4}},\G_{\un{2}}\ra=1.$\edm \demo{Claim 43}\
$\la\G_{\un{2}},\G_{\un{3}}\ra=1.$\edm

$\qquad\qquad \la\G_{\un{7}},\G_{\un{8}}\ra=1.$
\demo{Claim 44}\ $\la\G_{\un{5}},\G_{\un{8}}\ra=1.$\edm
\demo{Claim 45}\ $\G_8=\G_{7'}\G_5\G_3^{-1}\G_4^{-1}\G_{2'}
\G_4\G_3\G_5^{-1}\G_{7'}^{-1}.$\edm
\demo{Claim 46}\ $\G_{2'}=\G_{4'}\G_4\G_{3'}\G_3\G_{2'}
\G_2\G_{2'}^{-1}\G_4^{-1}\G_{3}^{-1}\G_{3'}^{-1}\G_{4}^{-1}\G_{4'}^{-1}.$
\edm
\demo{Claim 47}\ $\G_{8'}=\G_{5}^{-1}\G_{5'}^{-1}\G_{7}^{-1}\G_{7'}^{-1}
\G_8\G_{7'}\G_7\G_{5'}\G_5.$\edm
\demo{Claim 48}\ $[\G_{\un{i}},\G_{\un{9}}]=1 \qquad
i=1,2,3,4,6,7.$\edm

\demo{Proofs of the Claims}

We use the braids in $\Ce(18)$ (see Theorem I.1) to induce relations on $G$ via
the Van Kampen Theorem (Theorem 0.5).
For every factor $V^\nu $ in $\ve(18),$ we have to find $A_V$ and $B_V$ to get
a relation.
In Remark 0.4 one can find an algorithm how to determine $A_V$ and $B_V.$
In $\ve(18)$ we have sometimes used a compact notation for a
product of a few factors.  Then we use Remark I.1 to determine the factors
precisely. Sometimes, instead of using a factor $b$ in $\Ce(18),$ we shall use
its complex conjugation $\ov{b}.$  In that way, we get  $R(b)$ and/or $R(\ov
b)$ a relation on $G$ induced by the Van Kampen Theorem. We also use Corollary
I.5 to get other relations using $(R)\rho_{m_1,\dots,m_4,m_6,m_9}.$

\smallskip

\demo{Proof of Claim 0}\ Taking the factors $C'_i$ or the complex
conjugate of $C'_i\ \ i=1\dots 4$ and applying the Van Kampen
method on it to produce relations on $G,$ we get $[\G_{ii'},
\G_{jj'}]=1\ \forall i,j$ s.t. $\hat{L}_i\cap\hat{L}_j=\emptyset\
j\neq 7, 9.$
By Lemma II.5 we get $[\G_{\un{i}},\G_{\un{j}}]=1 \ \forall i,j$ s.t.
\linebreak $\hat{L}_i\cap\hat{L}_j=\emptyset\ j\neq 7, 9.$ For $j=7$ we
consider $C_4'.$ \bigskip
In $C'_4$ we have $\tilde{Z}_{17}=\underset{(6,6')}\to{\overline Z}_{17}: $
\bigskip
Its complex conjugate is $\overset {(6,6')}\to {\underline Z}_{17}:$
\bigskip
\flushpar which implies on $G$ the relation:
$$[\G_1,\G_6\1\G_{6'}\1\G_7\G_{6'}\G_6]=1.$$
Since we already know that $[\G_{\un{1}}, \G_{\un{6}}]=1$ we get
$[\G_1,\G_7]=1.$
Now, we are using the Invariance Theorem (Corollary I.6) and we get
$[\G_{\un{1}},  \G_{\un{7}}]=1.$ \enddemo

\demo{Proof of Claim 1}\ From $H'_i,$\ $i\neq 4$ and Corollary I.5.
\enddemo

\demo{Proof of Claim 2}\ By using $Z^{(3)}_{2', 33'}$ of $H'_4$ we get
$\la\G_{2'},\G_3\ra=\la\G_{2'},\G_{3'}\ra=\la\G_{2'},\G_{3'}\G_3\G_{3'}^{-1}\ra.$
We apply on it the invariance automorphism $(\rho_3\rho_7)^{m_3}$ for all
possible values of $m_3$ to get $\la\G_{2'},\G_{\un{3}}\ra = 1.$\enddemo

\demo{Proof of Claim 3} By definition of $\bullet$ (see Theorem I.1).
$$\align \G_{3'}^\bullet& =
\G_{3'}\G_3\G_{2'}\G_3^{-1}\G_{3'}\G_3\G_{2'}^{-1}\G_3^{-1}\G_3^{-1}\\
& =\G_{3'}\G_3
(\G_3^{-1}\G_{3'}\G_3)^{-1}\G_{2'}(\G_3^{-1}\G_{3'}\G_3)\cdot\G_3^{-1}
\G_3^{-1}\\
&\overset{\text{by Claims 2 and II.4(a) }}\to=\G_3\G_{2'}\G_3^{-1} \\
 & \overset{\text{by Claims 2 and
II.4(a)}}\to=\G_{2'}^{-1}\G_3\G_{2'}.\endalign$$\enddemo

\demo{Proof of Claim 4} By definition of $\bullet$ (see Theorem I.1).
$$\align \G_{3'}^\bullet &=
\G_{3'}\G_3\G_{2'}\underbrace{\G_3^{-1}\G_{3'}^{-1}\G_3\G_{3'}\G_3}\G_{2'}^{-1}\G_3^{-1}
\G_3^{-1}\\
&\overset{\text{by Claims 2 and II.4(a)}}\to  =\G_{3'}\G_3
(\G_3^{-1}\G_{3'}^{-1}\G_3\G_{3'}\G_3)^{-1}\G_{2'}(\G_3^{-1}\G_{3'}^{-1}
\G_3\G_3\G_{3'}\G_3)\G_3^{-1} \G_{3'}^{-1}\\
& =\G_3^{-1}\G_{3'}\G_3\G_{2'}\G_3^{-1}
\G_{3'}^{-1}\G_3    \\
&\overset{\text{by Claims 2
and II.4(a)}}\to=
  \G_{2'}^{-1}\G_3^{-1}\G_{3'}\G_3\G_{2'}.\endalign$$\enddemo

\demo{Proof of Claims 5,6,7}\ Arguments symmetric  to $2,$ $3,$ $4.$\enddemo

\demo{Proof of Claim 8}\ Easy to see from the geometric observation of
the action of $b^\bullet$.\enddemo

\demo{Proof of Claim 9}\ $b^\bullet$ does not affect those loops.\enddemo

\demo{Proof of Claim 10}\ By $H_4'$ we have $Z^2_{2', i}$ and
$Z^2_{2i}$ for $i=8,8', 7,7', 5, 5'.$
Thus, by the Van Kampen Theorem, $[\G_{2'},\G_{ii'}]=[\G_2,
\G_{ii'}]=1$\ \ $i=5,7,8.$
By Claim II.5, we have $[\G_{2'},\G_{\un{i}}]=[\G_2,\G_{\un{i}}]=1$\ \
$i=5,7,8.$
By Claim II.5, again $[\G_{\un{2}},\G_{\un{i}}]=1.$\enddemo

\demo{Proof of Claim 11}\ By $H_4'$ we have $Z^{(2)}_{44',8'}.$
Using the Van Kampen method $[\G_{4},\G_{8'}]=1$ and $[\G_{4'},
\G_{8'}]=1$
By Lemma II.5, $[\G_{\un{4}},\G_{8'}]=1.$
We have in $H_4',$\
$(Z_{33',8'})^{2\bullet}.$ By the Van Kampen method we have $[\G_3^\bullet,
\G_{8'}^\bullet]=[\G_{3'}^\bullet, \G^\bullet_{8'}]=1.$
  By Claims 4, 3, and 9, $\G_{3'}^\bullet=\G_{2'}\1\G_3\G_{2'},$\
$\G_3 =\G_{2'}\1\G_3\1\G_{3'}\G_3\G_{2'},$ \
$\G_{8'}^\bullet=\G_{8'}$ and  thus
$[\G_{2'}\1\G_3\G_{2'},\G_{8'}]=
[\G_{2'}\1\G_3\1\G_{3'}\G_3\G_{2'}, \G_{8'}]=1.$
 $[\G_{3'};\G_{8'}]=1.$
  Since $[\G_{2'},\G_8]=1,$
(Claim 10) we get by Lemma II.4(g),  $[\G_{3'},\G_{8'}]=1$ and
$\G_3\1\G_{3'}\G_3, \G_{8'}]=1.$ By Lemma II.5, $[\G_{\un{3}},
\G_{8'}]=1.$\enddemo

\demo{Proof of Claim 12}\ By Claim 11, $[\G_{\un 3}, \G_{8'}]=1.$
We apply on it Corollary I.5 with $ (Z_{88'}Z_{22'})^m$ to
get $[\G_{\un 3},\G_{\un 8}]=1.$
Similarly, we get $[\G_{\un 4},\G_{\un 8}]=1.$
\enddemo

\demo{Proof of Claim 13}\ We shall only prove $(\a,\be)=(1, 2).$
In $H_1'$ we have $\tilde{Z}_{22'(1)}.$
By the Van Kampen Theorem we get
$\G_2=\G_1\1\G_{1'}\1\G_2\1\G_{2'}\G_2\G_{1'}\G_1.$ We apply on it $Z_{22'}^2$
to get $$\G_{2'}\G_2\G_{2'}\1=\G_1\1\G_{1'}\1\G_{2'}\G_{1'}\G_1.$$
Thus, $E'_2 = E_1\1 E_{1'}\1 E_{2'}E_{1'} E_1.$
The other relations are induced from $H_2', H_3', H_5', H_6',H_7'.$
\enddemo

\demo{Proof of Claim 14}\ In $\hat{F}$ we have $\left(\ov{Z}^2_{3'7'}
\right)^\bullet.$
After complex conjugation this braid transforms to
$\left(\un{Z}^2_{3'7'}\right)^\bullet.$
By Corollary I.4, $\left( \un Z^2_{3'7'}\right)^\bullet$ implies, via the Van
Kampen
method, the following relation on $G$:\ $[\G_{3'}^\bullet,\G_{7'}^\bullet]=1.$
\enddemo

\demo{Proof of Claim 15}
$$\align
 1 \overset{\text{Claim 14}}\to = [\G_{3'}^\bullet, \G_{7'}^\bullet ]
&=[\G_{2'}\1 \G_3 \G_{2'}, \G_8\G_{7'}\G_7\G_{7'}\1 \G_8]
\overset{\text{Claims 10 and II.4(g)}}\to =
[\G_{3}, \G_8\G_{7'}\G_7\G_{7'}\1\G_8]\\
&\overset{\text{Claims 12 and II.4(g)}}\to =
[\G_{3},\G_{7'}\G_7\G_{7'}\1].\endalign$$  \enddemo

\demo{Proof of Claim 16}

By Claim 1:	$[\G_1,\G_{\un{7}}]=1.$ Thus,
$[\G_{1'}\G_1,\G_{7'}\G_7\G_{7'}\1]=1.$

By Claim 15: $[\G_3,\G_{7'}\G_7\G_{7'}\1]=1.$

Thus, $[\G_1\1\G_{1'}\1\G_3\G_{1'}\G_1, \G_{7'}\G_7\G_{7'}\1]=1.$

But by Claim 13: $\G_{3'}=\G_1\1\G\1_{1'}\G_3\G_{1'}\G_1.$

Thus, $[\G_{3'},\G_{7'}\G_7\G_{7'}\1]=1.$\enddemo

\smallskip

\demo{Proof of Claim 17}
\ Claim 16, Claim 15 and Claim II.4(b).\enddemo

\demo{Proof of Claim 18}\ We apply $(Z_{33'} Z_{77'})^m$ on $[\G_{3'},
\G_{7'}\G_7\G_{7'}\1]=1$ and on $[\G_3,\G_{7'}\G_7\G_{7'}\1]=1$ to get
$[\G_3, \G_{7'}]= [\G_3\G_{3'}\G_3\1, \G_{7'}]=1.$
By Claim II.4(b) $\G_{\un{3}},\G_{7'}]=1.$
We use Claim 17 and Claim II.4(b) to get $[\G_{\un{3}},\G_{\un{7}}]=1.$
\enddemo

\demo{Proof of Claim 19}\ $[\G_{2'},\G_8]=[\G_{2'}\G_7]=[\G_{2'}\G_{7'}]$
(Claim
10), thus $[\G_{2'},\G_{7'}^\bullet]=[\G_{2'},\G_{7}^\bullet]=1\Rightarrow
[\G_{2'},\G_7^\bullet]=1.$  Now,
$[\G_{\un{3}},\G_{\un{7}}]=[\G_{{\un{3}}},\G_8]=1$ (Claim 11, Claim 18). Thus,
$[\G_{\un 3}, \G_{\un{7}}^\bullet]=1.$ From $[\G_{2'},\G_{\un
7}^\bullet]=[\G_{\un 3},\G_7^\bullet]=1,$ we get
$[\G_3^\bullet,\G_{\un{7}}^\bullet]=[\G_{3'}^\bullet, \G_{\un{7}}^\bullet]=1.$
Thus, $[\G_{\un{3}}^\bullet,\G_{\un{7}}^\bullet]=1.$\enddemo

\demo{Proof of Claim 20}
In $\hat{F}$ we have $\left(Z_{33',44'}^{(3)}\right)^\bullet,$ which by the Van
Kampen Theorem implies $\la\G_{3}^\bullet, \G_{4}^\bullet\ra=1.$ We apply on it
$(Z_{33'} Z_{77'})^{m_3} (Z_{44'} Z_{55'})^{m_4}$ for all possible $m_4$ and
$m_3$ to get $\la\G_{\un 3}^\bullet, \G_{\un 4}^\bullet\ra=1$\quad
$(\{(\G_3^\bullet)\rho^m\}=\G_{\un 3}^\bullet).$\enddemo

\demo{Proof of Claim 21}\ Same proof as Claim 20.
\enddemo

\medskip

\demo{Proof of Claim 22}\ In $\hat{F}$ we have $\a_{2}^\bullet$
where $\a_2   $ is described in Fig.I.0(d).

 By the  Van Kampen Theorem we get
$$\G_{7}^\bullet\G_{5'}^\bullet\G_{7}^{\bullet{\1}}=\G_{4'}^\bullet\G_{4}^\bullet
\G_{3'}^\bullet\G_{4}^\bullet\G_{3'}^{\bullet{\1}}\G_{4}^{\bullet{\1}}
\G_{4'}^{\bullet{\1}}.$$

Since, $\la\G_{\un{4}}^\bullet,\G_{\un{3}}^\bullet\ra=1$
we can apply Claim II.4(a) to get
$$\align
&\quad =
\G_{4'}^\bullet\G_4^\bullet\G_4^{\bullet\1}\G_{3'}^\bullet\G_4^\bullet\G_4^{\bullet\1}
\G_{4'}^{\bullet\1}\\
&\quad = \G_{4'}^\bullet\G_{3'}^\bullet\G_{4'}^{\bullet\1}.\endalign$$

 We apply
Claim II.4(a) again to get
$$=\G_{3'}^{\bullet\1}\G_{4'}^\bullet\G_{3'}^\bullet.$$\enddemo

\demo{Proof of Claim 23}\ Directly from $\left(\overset{(44')}\to
{\un{Z}}_{3'7}^\bullet\right)^2$ in $\hat{F}.$ \enddemo

\demo{Proof of Claim 24}\

By Claim 20 $\la\G_{3'}^\bullet\G_3^\bullet\G_{3'}
^{\bullet \1},
\G_4^\bullet\ra=1.$

Thus$\la\G_3^\bullet,\G_{3'}^{\bullet
\1}\G_4^\bullet\G_{3'}^\bullet\ra =1.$

Since $ [\G_3^\bullet, \G_7^\bullet]=1,$\ $\la\G_3^\bullet,\G_{7}^{\bullet \1}
\G_{3'}\1\G_4^\bullet\G_{3'}^\bullet\G_7^\bullet\ra=1.$

By Claim 22\
$\la\G_{3'}^\bullet, \G_5^\bullet\ra=1.$\enddemo

\demo{Proof of Claim 25}\ Follows from Claim 21, as in Claim 24. \enddemo

\demo{Proof of Claim 26}\ Follows from Claim 24, using $[\G_{2'},
\G_{\un{5}}]=1.$
\enddemo

\demo{Proof of Claim 27}\ From Claim 25, using $[\G_8,\G_{\un{4}}]=1.$
\enddemo

\demo{Proof of Claim 28}\

By Claim 23 $[\G_7^\bullet, \G_{4'}^\bullet\G_4^\bullet
\G_{3'}^\bullet\G_{4}^{\bullet^{-1}}\G_{4'}^{\bullet^{-1}}]=1.$

Thus $[\G_{4'}^{\bullet^{-1}}\G_7^\bullet\G_{4'}^\bullet,
\G_4^\bullet\G_{3'}^\bullet \G_{4}^{\bullet^{-1}}]=1.$

{}From Claim 27, Claim 20 and Claim II.4(a)\quad
$[\G_7^\bullet\G_{4'}^\bullet\G_{7}^{\bullet^{-1}},
\G_{3'}^{\bullet^{-1}}\G_4^\bullet\G_{3'}^\bullet]=1.$

By Claim 22 $[\G_7^\bullet\G_{4'}^\bullet\G_{7}^{\bullet^{-1}},
\G_7^\bullet\G_5^\bullet\G_{7}^{\bullet^{-1}}]=1.$

Thus, $[\G_{4'}^\bullet, \G_5^\bullet]=1.$

Since $\G_{4'}^\bullet=\G_{4'}, \G_5^\bullet = \G_5$ we get the Claim.
\enddemo

\demo{Proof of Claim 29}\

By Claim 13 $\G_{6'}=\G_4\1 \G_{4'}\1\G_6\G_{4'}\G_4.$

Thus, $\G_4\G_{6'}\G_4\1=\G_{4'}\1\G_6\G_{4'}.$

By Claim 1 and Claim II, 4(a)
$$\G_{6'}\1\G_4\G_{6'}=\G_6\G_{4'}\G_6\1.$$

Thus,
$$\G_{4}=\G_{6'}\G_6\G_{4'}\G_6\1\G_{6'}\1.$$

We substitute the last equation in Claim 29 to get
$$[\G_{6'}\G_6\G_{4'}\G_6\1\G_{6'}\1,\G_5]=1.$$

By Claim 9 $[\G_5,\G_{\un{6}}]=1.$\ Thus,
$$[\G_{4'},\G_5]=1.$$ \enddemo

\demo{Proof of Claim 30}\ By Claim 28, Claim 29 and Claim II.5(a).\enddemo

\demo{Proof of Claim 31}\ We apply $Z_{44'} Z_{55'}$ on $[\G_{4'},\G_5]=1$ and
 on $[\G_4, \G_5]=1$ to get
$[\G_{4'}\G_4 \G_{4'}\1, \G_{5'}]=[\G_{4'},\G_{5'}]=1.$
We then use Lemma II.5(b) to get $[\G_{\un{4}}, \G_{5'}]=1.$
Together with Claim 30, we get the Claim. \enddemo

\demo{Proof of Claim 32}\  By Claim 7 and Claim 10.\enddemo

\demo{Proof of Claim 33}\ $$[\G_3,\G_7^\bullet]=[\G_3,\G_8\G_{7'}\G_8\1]
\overset{\text{By Claims 12 and II.4(g)}}\to = [\G_3,\G_{7'}]_{\G^\bullet}\
\overset{\text{By Claim 18}}\to = 1.$$\enddemo

\demo{Proof of Claim 34}

$  [\G_7^\bullet, \G_5^2\G_3\G_5^{-2}]
\overset{\text{Claims 10 and 32}}\to =
[\G_7^\bullet, \G_5^2\G_{2'}\1\G_3\G_{2'}\G_5^{-2}]
\overset{\text{Claim 3}}\to= [\G_7^\bullet, \G_5^2\G_{3'}^\bullet \G_5^{-2}].$

Thus, $[\G_7^\bullet, \G_5^2\G_{3'}^\bullet \G_5^{-2}]\Leftrightarrow
[\G_{5\1}\G_7^\bullet\G_5,\G_5\G_{3'}^\bullet\G_5\1]=1.$

Now: $$\align
&[\G_5\1\G_7^\bullet\G_5, \G_5\G_{3'}^\bullet\G_5\1]\\
&\quad\overset{\text{Claim 9}}\to = [\G_5^{\bullet{}\1}\G_7^\bullet
\G_5^{\bullet{}\1}, \G_5^\bullet\G_{3'}^\bullet\G_5{^\bullet{}\1}]\\
 &\quad\overset{\text{Claims 21,  \ 24,
and II.4(a)}}\to =
[\G_7^\bullet\G_5^\bullet\G_7^{\bullet{}^{-1}},\G_{3'}^{\bullet{}^{-1}}\G_5
^\bullet\G_{3'}^\bullet]=\\
&\quad\overset{\text{Claim 22}}\to =
[\G_{3'}^{\bullet{}\1}\G_4^\bullet\G_{3'}^\bullet,\G_{3'}^{\bullet{}\1}\G_5^\bullet
\G_{3'}^\bullet]=1\\
&\quad\overset{\text{Claim II.4(f)}}\to =
[\G_{4}^\bullet,\G_5^\bullet]_{\G_{3'}{}^\bullet}=1\\
&\quad\overset{\text{Claim 9}}\to = [\G_{4}, \G_5]_{\G_{3'}{}^\bullet}\\
&\quad\overset{\text{Claim 29}}\to =1.\endalign$$\enddemo

\demo{Proof of Claim 35}\ Let $f:B_5\ri G$ be as follows:

$B_5=\la X_1,\dots, X_4\bigm| [X_i,X_j]=1\ |i-j|>2,\ \la X_i,X_{i+1}\ra=1,\
i=1,2,3\ra.$ $f(X_1)=\G_{2'}, f(X_2)=\G_3, f(X_3)=\G_5, f(X_4)=\G_7^\bullet.$

By Claims 10, 32, 33, 2, 26, 21, $f$ is well-defined.

Let $d$ be the braid that satisfies $(X_i)_d=X_{5-i},\ i = 1, 2, 3.$

Then
$$\align
1&=[\G_7^\bullet, \G_5^2\G_3\G_5^{-2}]=f[X_4, X_3^2X_2X_3^{-2}]\\
&= f[(X_1)_d, (X_2)_d^2 (X_3)_d(X_2)_d^{-2}]\\
&=f ([X_1, X_2^2X_3X_2^{-2}]_d)\\
&=[f(X_1), f(X_2)^2 f(X_3) f(X_2)^{-2}]_{f(d)}.\endalign$$

Thus, $[f(X_1), f(X_2)^2 f(X_3) f(X_2^{-2}]=1.$

Thus, $[\G_{2'}, \G_3^2\G_5\G_3^{-2}]=1.$\enddemo
\smallskip

\demo{Proof of Claim 36}

$[\G_3, \G_4] \overset{\text{Claim 9}}\to =[\G_3,\G_4^\bullet]
\overset{\text{Claim 22}}\to =
[\G_3,\G_{3'}^\bullet\G_7^\bullet\G_5^\bullet\G_7^{\bullet^{-1}}\G_{3'}
^{\bullet^{-1}}].$

Thus, $$\align &[\G_3,\G_4]=1\\&
\Leftrightarrow
[\G_3,\G_{3'}^\bullet\G_7^\bullet\G_5^\bullet\G_7^{\bullet^{-1}}
\G_{3'}^{\bullet^{-1}}]=1\\ &\overset{\text{Claims 18
and 19}}\to\Leftrightarrow\
[\G_3,\G_{3'}^\bullet\G_5^{\bullet}\G_{3'}^{\bullet^{-1}}]=1\\
&\Leftrightarrow\ [\G_3,\G_{2'}\1 \G_3\G_{2'}\G_5\G_{2'}\1\G_3\1\G_{2'}]=1\\
&\Leftrightarrow\ [\G_{2'}\G_3\G_{2'}\1, \G_3\G_{2'}\G_5\G_{2'}\1\G_3\1]=1\\
&\overset{\text{Claim 10}}\to\Leftrightarrow\ [\G_{2'}\G_3\G_{2'}\1,
\G_3\G_5\G_3\1]=1\\ &\overset{\text{Claims 2 and II.4(a)}}\to\Leftrightarrow\
[\G_3\1\G_{2'}\G_3,\G_3\G_5\G_3\1]=1\\
&\Leftrightarrow\ [\G_{2'},\G_3^2\G_5\G_3^{-2}]=1\endalign$$

which is true by Claim 35.

We get $[\G_{\un{3}},\G_{\un{4}}]$ by Corollary I.5.\enddemo

\demo{Proof of Claim 37}\ Similar to Claim 34.\enddemo

\demo{Proof of Claim 38}\ Claim 18, Claim 10, Claim 12, Claim 3.\enddemo

\demo{Proof of Claim 39}\ Let $f:B_5\ri G$ be as follows:
$$f(X_1)=\G_8, f(X_2)=\G_{7'}, f(X_3)=\G_4, f(X_4)=\G_{3'}^\bullet.$$

By Claim 12, Claim 38, Claim 5, Claim 27, Claim 20, $f$ is well defined.
Let $d$ be a braid such that $(X_i)_d=X_{4-i},\ i=1,2,3,4.$

Now:
$$\align
1=[\G_{7'},\G_4^2\G_{3'}^\bullet\G_4^{-2}]
&=[f(X_2), f(X_3)^2 f(X_4) f(X_3)^{-2}]\\
&=f[X_2, X_3^2 X_4 X_3^{-2}]\\
&=f[(X_3)_d, (X_2)^2_d (X_1)_d (X_2)^{-2}_d]\\
&=[f(X_3), f(X_2)^2 f(X_1) f(X_2)^{-2}]_{f(d)}\\
&=[\G_4,\G_{7'}^2\G_8\G_{7'}^{-2}]_{f(d)}.\endalign$$

Thus, $[\G_4,\G_{7'}^2\G_8\G_{7'}^{-2}]=1.$\enddemo
\smallskip
\demo{Proof of Claim 40}\ Follows from Claim 30, like in Claim 36.\enddemo

\demo{Proof of Claim 41}\ Apply $(Z_{44'} Z_{55'})^{m_4}\ (Z_{77'}
Z_{33'})^{m_3}$
on $[\G_5, \G_{7'}]=1$ from Claim 40.\enddemo

\demo{Proof of Claim 42}\ In $H'_4$ there is $\ov{Z}^{(3)}_{2', 44'}.$
We take its complex conjugate $\un{Z}^{(3)}_{2', 44'}$ and apply Corollary I.4
on $\un{Z}^{(3)}_{2', 4}$ to get $\la\G_{2'}, \G_4\ra=1$ and then apply
Corollary I.5 to get $\la\G_{\un{2}}, \G_{\un{4}}\ra=1.$\enddemo

\demo{Proof of Claim 43}\ We apply Corollary I.5 on Claims 2 and 5.\enddemo

\demo{Proof of Claim 44}\ Similar to Claim 42.\edm

\demo{Proof of Claim 45}\ By Claim 22,
$\G_7^\bullet\G_5^\bullet\G_7^{\bullet^{-1}}=\G_{3'}^{\bullet^{-1}}
\G_4^\bullet \G_{3'}^\bullet.$
Since $\la\G_5^\bullet,\G_7^\bullet\ra=1$ we have
$$\G_5^{\bullet^{-1}}\G_7^\bullet\G_5^\bullet=\G_{3'}^{\bullet^{-1}}
\G_4^\bullet\G_{3'}^\bullet.$$

Thus, $\G_7^\bullet =
\G_5^\bullet\G_{3'}^{\bullet^{-1}}\G_4^\bullet\G_{3'}^\bullet
\G_5^{\bullet^{-1}}.$

We substitute the formulas for $\G_{3'}^\bullet$ and $\G_7^\bullet$ to
get
$$\G_{7'}\1\G_8\G_{7'}=\G_5\G_{2'}\1\G_3\1\G_{2'}\G_4\G_{2'}\1
\G_3\G_{2'}\G_5\1.$$

Since $\la\G_{2'},\G_4\ra=1$ (Claim  42)
$$\G_{7'}\1\G_8\G_{7'}=\G_5\G_{2'}\1\G_3\1\G_4\1\G_{2'} \G_4 \G_3
\G_{2'}\1\G_5\1.$$

Since $[\G_5,\G_{3'}]=1$
$$\G_{7'}\1\G_8\G_{7'}=\G_{2'}\1\G_5\G_3\1\G_4\1\G_{2'}\G_4\G_3\G_5\1\G_{2'}.$$

Since $[\G_{2'},\G_8]=[\G_{2'},\G_{7'}]=1$
$$\G_{7'}\1\G_8\G_{7'}=\G_5\G_3\1\G_4\1\G_{2'}\G_4\G_3\G_5\1$$

and
$$\G_8=\G_{7'}\G_5\G_3\1\G_4\1\G_{2'}\G_4\G_3\G_5\1\G_{7'}\1.$$\enddemo

\demo{Proof of Claim 46}\ By $\tilde{Z}_{22'}$ of $H'_4$ we have:
$$\G_2=\G_{4'}\G_4\G_{3'}\G_3\G_{2'}\G_3\1\G_{3'}\1\G_4\1\G_{4'}\1.$$

We apply on it $Z_{88'}Z_{22'}$ using Corollary I.5 to get (recall that
$(\G_2) Z_{22'}'Z_{88'}=\G_{2'}$
and $(\G'_2) Z_{22'}Z_{88'}=(\G'_2) Z_{22'}=\G_{2'}\G_2\G_{2'}\1)$
$$\G_{2'}=\G_{4'}\G_4\G_{3'}\G_3\G_{2'}\G_2\G_{2'}\1\G_3\1\G_{3'}\1
\G_4\1\G_{4'}\1.$$\enddemo

\demo{Proof of Claim 47}\ By $\tilde{Z}_{88'}$ of $H'_4$ we get, using the
Van Kampen Theorem:
$$\G_{8'}=\G_{5}\1
\G_{5'}\1\G_{7}\1\G_{7'}\1\G_{8}\G_{7'}\G_7\G_{5'}\G_5.$$\enddemo

\demo{Proof of Claim 48}\

In $C'_5$ we have $\tilde{Z}_{i9}^2 = H(\tilde z_{i9})^2$ for $ i=
1,2,3,4,6,7,$  where:

\bk

$\tilde z_{i9} =\underset(8,8')\to{\ov{z}}_{i9}$

{\it{\ub{Case 1}}}.\quad  $i\neq 7$

The complex conjugation of $\tilde{Z}_{i9}$ is $\ov{\tilde
Z}_{i9}=H(\ov{\tilde z}_{i9})$ and

 $\ov{\tilde z}_{i9}= \overset{(8,8')}\to{\underline {z}}_{i9}$
which implies on $G$ (as in Claim
0) $\quad[\G_i,\G_8\1\G_{8'}\1\G_9\G_{8'}\G_8]=1.$

By Claim 10, for $i=2,3,4$ we have $[\G_i,\G_{\un{8}}]=1.$

By Claim 0, for $i=1,6$ we have $[\G_i,\G_{\un{8}}]=1.$

Thus, $[\G_i,\G_9]=1.$
We use Corollary I.5 to get $[\G_{\un{i}},\G_{\un{9}}] = 1.$

{\it{\ub{Case 2}}}.\quad $i= 7$

In $C'_5$ we have $\un{Z}_{7,9}^2$ which implies $[\G_7,\G_9]=1$ and
thus, by Corollary I.7,\newline $[\G_{\un{7}}, \G_{\un{9}}]=1.$

\mk

{}From Case 1 and Case 2 we get Claim 48.\enddemo

Thus we have proved all the Claims.
\mk

We shall now prove the statements of the Proposition.
We use the above claims, the definition of $E_i, E_{i'}$ and the
facts $\{\G_{\un{i}}\}=\{E_{\un{i}}\}$ and\ \  $E_{i'} E_i=\G_{i'} \G_i\
\forall
i.$ \roster \item From Claim 1, Claim 26, Claim 27, Claim 42, Claim 43, Claim
44. \item From Claim 0, Claim 10, Claim 12, Claim 18, Claim 31, Claim 36, Claim
41, Claim 48.
\item From Claim 45, definifions of $E_i$ and Corollary I.5.
\item From Claim 46, definitions of $E_i$ and Corollary II.1.
\item From Claim 47, definition of $E_i$ and the above fact.
\item From Claim 13.
\item From Claim 13. \ \ \ \ \ \endroster

\flushpar\qed for Proposition II.6

We need the following corollary in order to obtain in Chapter IV, \S8 a smaller
set of generators for $G.$
\proclaim{Corollary II.7}
Let $E_i,E_{i'}$ be as in the beginning of the chapter.
\flushpar Let
$A_i=E_{i'}E^{-1}_i$. Then:
\roster\item $A_5=(A_4)_{E_2^{-1}E_3E_7^{-1}E_8}$
\item $A_7=(A_3)_{E_2\1 E_4E_5\1 E_7\1}$
 \item $A_8=(A_2)_{E_4E_3E_5\1 E_7^{-1}}$
\item $A_2=E_1^{-2}A_1\1 (A)_{E_2^{-1}}(E_1^{2})_{E_2^{-1}}$
\item $A_6=E_4^{-2}A_4^{-1}(A_4)_{E_6^{-1}}(E_4^2)_{E_6^{-1}}$
\item $A_3=E_1^{-2}A_1\1 (A)_{E_3^{-1}}(E_1^{2})_{E_3^{-1}}$
\item $A_9=E_5^{-2}A_5^{-1}(A_5)_{E_9^{-1}}(E_5^2)_{E_9^{-1}}$
\item
$(A_4)_{E_2^{-1}E_4^{-1}}=E_4^2A_3E_3^2A_2(E_3^{-2})_{E_2^{-1}}(A_3^{-1})_{E_2^{-1}}
(E_4^{-2})_{E_2^{-1}}$
\item
$(A_7)_{E_6^{-1}E_7^{-1}}=E_7^2A_6(E_7^{-2})_{E_6^{-1}}$
\item
$(A_9)_{E_8^{-1}E_9^{-1}}=E_9^2A_8(E_9^{-2})_{E_8^{-1}}$
 \endroster
 \ep

\sk

\demo{Proof} We use proposition II.6 $(1),\dots,(7).$
The claims are grouped according to the similarity of their proofs and not
according to the order that we use them in Proposition IV.8.1.

{\it Recall}:
$$\align
(E_i)\rho_i&=E_{i\1}\\
(E_{i'})\rho_i&=E_{i'}E_iE_{i'}\\
E_{i'}&=A_iE_i\\
E_{i'}E_i&=A_iE_i^2.\endalign$$

 For (1), (2), (3)  we use (5) of  Proposition II.6.

For Claims (4)--(7) we use (6) of Proposition II.6.
The 4 claims are symmetric and we shall only prove the first one.

For (8) we use (4) of Proposition II.6.

For (9)--(10) we use (7) of Proposition II.6.
The claims are symmetric and we shall only prove the first one.
\roster
\item\ By (5) of Proposition II.6
$$E_8=(E_4\1E_2E_4)_{E_3E_5\1E_7\1}.$$
By (1), (2) and Claim II.4(a)
$$E_8=(E_2E_4E_2\1)_{E_3E_7\1E_5\1}.$$
Thus,
$$E_5\1E_8E_5=(E_4)_{E_2\1E_3E_7\1}.$$
By (1) and Claim II.4(a)
$$E_8E_5E_8\1=(E_4)_{E_2\1E_3E_7\1}.$$
Thus,
$$E_5=(E_4)_{E_2\1E_3E_7\1E_8}.$$
We apply on it  $\rho_4\rho_5$  to get
$$E_{5'}=(E_{4'})_{E_2\1E_3E_7\1E_8}.$$
We multiply the 2 results to get
$$A_5=(A_4)_{E_2\1E_3E_7\1E_8}.$$
\item \ By (5)  of Proposition II.6
$$E_7\1E_8E_7=E_5E_4\1E_3\1E_2E_3E_4E_5\1.$$
By (1) and Claim II.4(a)
$$E_8E_7E_8\1=E_4E_4\1E_2E_3E_2\1E_4E_5\1.$$
Thus,
$$E_7=(E_9)_{E_2\1E_4E_5\1E_8}.$$
We apply on this $\rho_2\rho_3$ to get
$$E_{7'}=(E_{3'})_{E_2\1E-4E-5\1E_8}.$$
We multiply the results to get
$$A_7=(A_3)_{E_2\1E_4E_5\1E_8}.$$
\item\ By (5) of Proposition II.6
$$E_8=(E_2)_{E_4E_3E_5\1E_7\1}.$$
We apply on it $\rho_2\rho_8$ to get
$$E_{8'}=(E_{2'})_{E_4E_3E_5\1E_7\1}.$$
We multiply the 2 results to get
$$A_8=(A_2)_{E_4E_3E_5\1E_7}.$$
\endroster

\roster\item"(4)"
{}From (6) of Proposition II.6
$$\align &E_{2'}=E_1\1E_{1'}\1E-2E_{1'}E_1\\
&\Ri A_2E_2=E_1\2 A_1\1E_2A_1E_1^2\\
&\Ri A_2=E_1\2A_1\1(A_1)_{E_2\1}(E_1^2)_{E_2\1}.\endalign$$
\item"(8)"
{}From (4) of Proposition II.6
$$\align &E_2=E_{4'}E_4E_{3'}E_3E_{2'}E_3\1E_{3'}\1E_4\1E_{4'}\1\\
&\Ri
E_2=A_4E_4^2A_3E_3^2A_2E_2E_3\2A_3\1E_4\2A_4\1\\
&\Ri
A_4\1E_2A_4=E_4^2A_3E_3^2A_2E_2E_3\2A_3\1E_4\2\\
&\Ri
A_4\1(A_4)_{E_2\1}=E_4^2A_3E_3^2A_2(E_3\2A_3\1E_4\2)_{E_2\1}\\
&\Ri
A_4\1(A_4)_{E_2\1}=E_4^2A_3E_3^2A_2(E_3\2)_{E_2\1}(A_3\1)_{E_2\1}(E_4\2)_{E_2\1}.
\endalign$$
By Lemma II.4(h)

$ (A_4)_{E_2\1E_4\1}=A_4\1(A_4)_{E_2\1}.$
\newline Thus,

$(A_4)_{E_2\1E_4\1}=E_4^2A_3E_3^2A_2(E_3\2)_{E_2\1}(A_3\1)_{E_2\1}(E_4\2)_{E_2\1}.$
\endroster

\roster\item"(9)"
{}From (7)
$$\align &E_6=E_{7'}E_7E_{6'}E_7\1E_{7'}\1\\
&\Ri E_6=A_7E_7^2A_6E_6E_7\2A_7\1\\
&\Ri A_7\1E-6A_7=E_7^2A_6E_6E_7\2\\
&\Ri
A_7\1(A_7)_{E_6\1}=E_7^2A_6(E_7\2)_{E_6\1}.\endalign$$
By Lemma II.4(h), $A_7\1(A_7)_{E_6\1}=(A_7)_{E_6\1E_7\1}.$
\newline Thus,

$(A_7)_{E_6\1E_7\1}=E_7^2A_6(E_7\2)_{E_6\1}.$\quad $\square$\endroster\edm

\end\magnification=1200
\parindent 20 pt

\NoBlackBoxes
\define\pf{\demo{Proof}}
\define\bk{\bigskip}
 \define\mk{\medskip}
\define\sk{\smallskip}
\define \a{\alpha}
\define \be{\beta}

\define \g{\gamma}
\define \lm{\lambda}
\define\ri{\rightarrow}
\define\Ri{\Rightarrow}
\define \s{\sigma}
\define\si{\sigma}
\define \fa{\forall}

\define \iy{\infty}

\define \r{\rho}

\define \tY{\tilde{Y}}\define \tg{\tilde{g}}\define \tT{\tilde{T}}
\define \tZ{\tilde{Z}}
\define \tB{\tilde{B}}\define \tP{\tilde{P}}
\define \tX{\tilde{X}}\define \tb{\tilde{b}}
\define \la{\langle}
\define \ra{\rangle}

\define \G{\Gamma}

\define \Ce{\Cal E}

\define \CP{\Bbb C\Bbb P}

\define \CPt{\Bbb C\Bbb P^2}

\define\BR{\Bbb R}
\define \F{\Bbb F}

\define \Z{\Bbb Z}
\define\BZ{\Bbb Z}
\define \ve{\varepsilon}
\define \vp{\varphi}

\define \Dl{\Delta}
\define \dl{\delta}
\define \C{\Bbb C}

\define \1{^{-1}}
\define \2{^{-2}}

\define \Aff{\operatorname{Aff}}
\define \Aut{\operatorname{Aut}}
\define \Gal{\operatorname{Gal}}

\define \Ab{\operatorname{Ab}}
\define \Center{\operatorname{Center}}
\define \un{\underline}
\define \ov{\overline}
\define \ub{\underbar}
\define \df{demo{Proof}}

\define \edm{\enddemo}
\define \ep{\endproclaim}

\define \Ss{S^{(6)}}
\baselineskip 20pt

\heading{\bf CHAPTER III.\ Construction theory for $\tB_n$}\rm\endheading

Let $D$ be a disc, $K\subset D,$ \ $\# K=n.$
Let $B_n=B_n[D,K]$ (see definition in Chapter 0).

\subheading{\S1. Definition of $\bold{\tilde B_n}$}
\demo{Definition}\ Let $D,K$ be as before.
Let $H(\s_1)$ and $H(\s_2)$ be 2 half-twists in $B=B_n[D,K].$
We say that $H(\s_1)$ and $H(\s_2)$ are:
\roster
\item"(i)" \ub{weakly disjoint} \ if  $\s_1\cap\s_2\cap K=\emptyset.$
\item"(ii)"  \ub{transversal}  \ if $\s_1$ and $\s_2$ are weakly disjoint
and  intersect each other exactly
once (and not in any point of $K),$ i.e.,
$\s_1\cap\s_2=$ one point,\ $\s_1\cap\s_2\cap K=\emptyset.$
\item"(iii)"  \ub{disjoint} \ if $\s_1\cap\s_2=\emptyset.
$
\item"(iv)" \ub{adjacent} \  if  $\s_1\cap\s_2\cap K=$ one point.
\item"(v)" \ub{consecutive} \  if they are adjacent and $\s_1\cap \s_2$ do
not intersect outside of $K,$ i.e.,  $\s_1\cap\s_2=$ point $\in K.$
\item"(vi)" \ub{cyclic} \quad if $\s_1\cap \s_2=$ 2 points $\in K.$
\endroster

\edm

\midspace{3.00in}\caption{Fig. III.1.0}

\proclaim{Claim III.1.0}
Let $X,Y$ be 2 half-twists in $B_n.$ Then:
\roster
\item"(i)"If $X,Y$ are disjoint, then $[X,Y]=1,$ i.e., $XY=YX.$
\item"(ii)" If $X,Y$ are consecutive, then $\la X,Y\ra=
XYXY\1X\1Y\1=1$ and \newline $XYX=YXY$ or
$X_{Y\1}=Y_X.$\newline  We say then that $X$ and $Y$ satisfy the triple
relation. \item"(iii)" If $X=H(x)=\ov\be,$\ $Y=H(y),$ then  $Y_X =
X\1YX=H((y)\be).$ \endroster\ep
\demo{Proof} (i)\ (ii)\ Basic properties of a braid group.
See, for example, Chapter 3 of \cite{MoTe4}.

(iii) See Fig. III.1.1 for a geometrical presentation of $Y_X.$
The common point of $x$ and $y$ travels under $\be$ counterclockwise to the
other end of $x.$ Thus $(y)\be$ connects the 2 ends of $x$ and $y,$ which is
not a common end of either of them.
By IV.3.4 of \cite{MoTe12}, $Y_X=X\1YX=(Y)X^\vee $ (see definition or $f^\vee$
in \cite{MoTe4}, Chapter 4) which is equal by Claim IV.3.0 of \cite{MoTe12} to
$H((y)\be_1).$\quad $\square$\edm

\midspace{2.00in}\caption{Fig. III.1.1}

\smallskip

\demo{Definition}\  $\underline{\tilde{B}_n}$

Let $\tilde{B}_n$ be the quotient of $b_n,$ the Braid group of order $n,$
by the subgroup generated by the commutators $[H(\s_1), H(\s_2)]$ where
$H(\s_1)$ and $H(\s_2)$ are transversal half-twists.
\edm

\smallskip

\demo{Notations}\

 Let $Y\in B_n.$  We denote the image of $Y$ in $\tilde{B}_n$
by  $\tilde{Y}.$
When possible, we shall abuse notation and denote $\tilde{Y}$
by $Y.$
If $Y$ is a half-twist in $B_n$ we call $\tilde{Y}$  a half-twist in
$\tilde{B}_n.$
We call two half-twists $\tilde{Y}, \tilde{X}$ in $\tilde{B}_n$ disjoint (or
weakly disjoint, adjacent, consecutive, transversal) if $Y,X$ are disjoint
(or
weakly disjoint, adjacent, consecutive, transversal).
If $\{X_i\}$ is a frame of $B_n,$ then $\{\tilde{X}_i\}$ is a frame in
$\tilde{B}_n.$
We also refer to $\{\tX_i\}$ as a standard base of $\tB_n.$\edm

\demo{Definition}\  $\un{\text{Polarized half-twist, polarization}}$

We say that a half-twist $X\in B_n$ (or $\tX$ in $\tilde B_n$) is polarized if
we choose an order on the end points of $X.$ The order is called the
polarization of $X$ or $\tX.$\edm

\demo{Definition }\ \underbar{Orderly adjacent}

Let $X,Y$ be two adjacent polarized half-twists in $B_n$ (resp.
in $\tilde{B}_n$).  We say that $X,Y$ are  {\it orderly adjacent}
if their common point is the ``end'' of one of them and the
``origin'' of another.
\edm

\demo{Definition}\ \underbar{Good quadrangle}

Let $H(\s_i)\ \ i=1,2,3,4$ be \ 4\ half-twists such that $H(\s_i)$ and
$H(\s_{i+1})$ are consecutive,
$H(\s_4)$ and $H(\s_1)$ are consecutive, $H(\s_1)$ and $H(\s_3)$ are disjoint
and $H(\s_2)$ and $H(\s_4)$ are disjoint, and in the interior of
$\bigcup\limits^4_{i=1}\s_i$ there is no point of $K.$ We say that
$\{H(\s_i)\}$ is a good quadrangle in $B_n$ and $\{\widetilde{H(\s}_i)\}$ is a
good quadrangle in $\tilde{B}_n.$\edm

\midspace{1.50in}\caption{Fig. III.1.2}

\proclaim{Remark III.1.1}
\roster\item"(a)" transversal, disjoint $\Rightarrow$ weakly disjoint.
\newline consecutive $\Rightarrow$ adjacent.
\item"(b)" Any two pairs of disjoint (or transversal, consecutive, cyclic)
half-twists are conjugate to each other by an element $b\in B.$
\item"(c)" Any two
half-twists    in $B_n$ (or $\tilde{B}_n$) are conjugate to each other by an
element of $B_n$ (or $\tB_n).$
\item"(d)" Every 2 transversal or disjoint half-twists in $\tilde{B}_n$
commute.
\newline Every 2 consecutive half-twists in $B_n$ or in $\tB_n$ satisfy the
triple relation $(XYX=YXY).$ \item"(e)" Any two good quadrangles in $B_n$ or in
$\tilde{B}_n$ are conjugate.
\item"(f)" Every 2 pairs of orderly adjacent (non-orderly adjacent) consecutive
half-twists are conjugate to each other by an element $b\in B$ preserving
polarization.\endroster
\ep \demo{Proof} Geometric observation in $B_n$ and $\tilde B_n,$ and Claim
III.1.0. \qed \edm
\proclaim{Lemma III.1.2}\ If $\{\tilde{Y}_i\}\ \ i=1,\dots, 4$ is
 a good quadrangle in $\tB_n,$ then\linebreak \rom{(a)}\ \ $\tilde{Y}_1
\tilde{Y}_3 =\tilde{Y}_3 \tilde{Y}_1,$\ \ \ \rom{(b)}\ \ $\tY_1^2\
\tY_3^3=\tY_2^2\ \tY^2_4.$\ep

\smallskip

\demo{Proof}
\roster
\item"(a)" Since they are disjoint.
\item"(b)" Let $X_1, X_2, X_3$ be 3 half-twists such that $X_1$ and $X_2$ are
consecutive, $X_2$ and $X_3$ are consecutive and $X_1$ and $X_3$ are disjoint.
Denote $X_i = H(x_i),\ \ i=1,\dots, 4.$
Clearly, $X_1, (X_3)_{X_{2}\1}, X_3, (X_1)_{X_{2}\1}$ is a good quadrangle (see
Fig. III.1.3).\endroster

\midspace{1.50in}\caption{Fig. III.1.3}

In $B_n:[X_1,X_3]=1$ and $\la X_1,X_2\ra=\la X_3, X_2\ra=1.$
Thus, we can use Claim II.4(a).

It is clear that $(X_2)_{X_{1} X_{3}}$ is transversal to $X_2.$
Thus, in $\tilde B_n:$\newline
$\left[\tX_2,(\tX_2)_{\tX_1\tX_3}\right]=1.$
Thus,

$\qquad\qquad\quad 1\ =[\tX_2,\tX_3\1 \tX_1\1\tX_2\tX_1\tX_3]$
\smallskip
$\qquad\qquad\quad \ \ \ = \tX_2\tX_3\1
\tX_1\1\tX_2\tX_1{\underbrace{\tX_3\tX_2\1\tX_3\1}}\tX_1\1\tX_2\1\tX_1\tX_3$
\baselineskip .03pt
\flushpar $\qquad  \qquad\qquad \qquad \qquad\qquad \qquad\quad\qquad\quad\ \
\overbrace{\ \quad \  \qquad}$

$\text{(By Claim\ II.4)} = \tX_2\tX_3\1
\tX_1\1\tX_2\tX_1\tX_2\1\tX_3\1{\underbrace{\tX_2\tX_1\1\tX_2\1}}\tX_1\tX_3$

$\text{(By Claim\ II.4)} = \qquad \qquad
\qquad \qquad \qquad \qquad \ {\overbrace{\tX_1\1\tX_2\1\tX_1}}\tX_1\tX_3$
\baselineskip .03pt
\flushpar $\qquad \qquad \qquad \qquad\qquad \qquad \qquad\quad
\qquad\quad\qquad \  \underbrace{\ \quad \  \qquad}$
\flushpar $\qquad \qquad \qquad \qquad\qquad
\qquad \qquad\qquad \qquad \quad\qquad  \overbrace{\ \quad \  \qquad}$
\flushpar $\qquad\   ([X_1,X_3]=1)  =\  \tX_2\tX_3\1\tX_1\1\tX_2
\underbrace{\tX_1\tX_2\1\tX_1\1}\ \tX_3\1\tX_2\1\tX_1\tX_1\tX_3$
\baselineskip 10pt
\flushpar $\qquad\qquad \qquad\qquad   \quad \qquad \qquad
\qquad  \quad \ \ {\overbrace{\tX_2\1\tX_1\1\tX_2}}$
\baselineskip 20pt

$\qquad\qquad\qquad = \quad \
\tX_2\tX_3\1\tX_1^{-2}\ {\underbrace{\tX_2\ \tX_3\1\ \tX_2\1}}\tX_1^2\tX_3$

$ \text{(By Claim\ II.4)} =\qquad\qquad\qquad\quad\ \
{\overbrace{\tX_3\1\tX_2\1\tX_3}}$

(Since $[X_1,X_3]=1)$

$\qquad\qquad\qquad =\  \tX_2\tX_3^{-2}\tX_1^{-2}\tX_2\1
\tX_1^2\tX_3^2$

$\qquad\qquad\qquad =\ (\tX_3)^{-2}_{\tX_{2}\1}\ (\tX_1)^{-2}_{\tX_{2}\1}
\tX_1^2 \tX_3^2.$

Thus, $\tX^2_1\tX_3^2 = (\tX_1)^2_{\tX_2\1}\cdot (\tX_3)^2_{\tX_2\1}.$

By Remark III.1.1, $\{\tY_i\}$ is conjugate to $\{\tX_1,
(\tX_3)_{\tX_2\1},\tX_3, (\tX_1)_{\tX_2\1}\}$ and thus satisfies $\tilde Y_1^2
\tilde Y_3^2= \tilde Y^2_4\tilde Y_2^2.$ \qed \edm

\bigskip
\subheading{\S2. $\bold{\tB_n}$-groups and prime elements}

\demo{Definition}\ $\un{\tB_n-\text{group}.}$

A group $G$ is called a $\tB_n$-group if there exists a homomorphism
$\tB_n\ri\Aut (G).$
We denote $(g)_b$ by $g_b.$\edm

\demo{Definition}\ $\un{\text{Prime element, supporting half-twist (s.h.t.)
corresponding central element}}
$

Let $G$ be a $\tB_n$-group.

An element $g\in G$ is called a prime element of $G$ if there exists a
half-twist $X\in B_n$ and $\tau\in \Center(G)$ with $\tau^2=1$ and
$\tau_b=\tau\ \forall\ b\in \tilde B_n$ such that:
\roster
\item $g_{\tX\1}=g\1 \tau.$
\item For every half-twist $Y$ adjacent to $X$ we have:
\item"" (a)\ $g_{\tX\tY\1\tX\1}=g_{\tX}\1 g_{\tX\tY\1}$
\item"" (b)\ $g_{\tY\1
\tX\1}=g\1 g_{\tY\1}.$
\item For every half-twist $Z$ disjoint from $X,$\ $ g_{\tZ}=g.$\endroster

The half-twist $X$ (or $\tX$) is called the {\it{supporting half-twist}} of
$g,$
\  ($X$ is the s.h.t. of $g.$)

The element $\tau$ is called the {\it{corresponding central element}}.\edm

\proclaim{Lemma III.2.1}

Let $G$ be a $\tB_n$-group.

Let $g$ be a prime element in $G$ with supporting half-twist $X$ and
corresponding central element $\tau.$
Then:
\roster
\item $g_{\tX}=g_{\tX\1}=g\1\tau,$\quad $g_{\tX^2}=g.$
\item $g_{\tY^{-2}}=g\tau \quad \forall \ \ Y$ consecutive half-twist to $X.$
\item $[g, g_{\tY\1}]=\tau \quad \forall \ \ Y$ consecutive half-twist to $X.$
\endroster\ep

\smallskip

\demo{Proof}
\roster
\item $g_{\tX^{-2}} = (g_{\tX\1})_{\tX\1} = (g\1\tau)_{\tX\1} =
(g\1)_{\tX\1}\cdot = (g\1\tau)\1 \tau=g\Ri
g_{\tX}=g_{\tX\1}\overset\text{Axiom (1)}\to=g\1\tau.$
 \item $g_{\tX\1 \tY\1 \tX\1}=
(g_{\tX\1})_{\tY\1\tX\1}= (g\1\tau)_{\tY\1\tX\1}= (g_{\tY\1\tX\1})\1\cdot\tau
\overset\text{Axiom (2)}\to =
g\1_{\tY\1}\cdot g\cdot \tau.$ \endroster

On the other hand:
$$\align g_{\tX\1 \tY\1 \tX\1}&=g_{\tY\1 \tX\1
\tY\1}=(g_{\tY\1\tX\1})_{\tY\1}\overset\text{Axiom (2)}\to =
 (g\1 g_{\tY\1})_{\tY\1}\\
&=  g\1_{\tY}\cdot g_{\tY^{-2}}.\endalign$$
Thus, $g_{\tY^-2}=g\tau.$
\roster
\item"(3)" $g_{\tX \tY\1 X\1} \overset{\text{Axiom (2)}} \to =
g_{\tX}\1\cdot g_{\tX\tY\1} \overset{\text{by (1)}} \to = g\1_X
g_{\tX\tY\1}= g\cdot g\1_{\tY\1}\cdot \tau^2.$\endroster

On the other hand:
$$\align g_{\tX \tY\1 \tX\1}& \overset\text{by (1)}\to=
(g\1\tau)_{\tY\1\tX\1}\\ & \overset\text{Axiom (2)}\to = (g\1\cdot
g_{\tY\1})\1\cdot\tau = g\1_{\tY\1}\cdot g\cdot \tau.\endalign$$

Thus, $g\cdot g\1_{\tY\1}\cdot  = g\1_{\tY\1}\cdot g\tau.$

Thus, $g_{\tY\1}\cdot g\cdot g_{\tY\1}\cdot g\1 =\tau\1 = \tau.$

Thus, $[g_{\tY\1}, g] = \tau.$ \qed \edm

\proclaim{Lemma III.2.2}
Let $G$ be a $\tB_n$-group.
Let $g$ be a prime element in $G$ with s.h.t. $X$ and corresponding central
element $\tau.$
Let $b\in \tB_n.$
Then $g_b$ is a prime element with s.h.t. $X_b$ and central element $\tau.$
\ep
\pf We use the fact that $(a_b)_c=(a_c)_{b_c}$ and $(ab)_c=a_cb_c.$
We have to prove 3 properties:
\roster\item $g_{\tX\1}=g\1\tau \Ri (g_{\tX\1})_b=(g\1\tau)_b \Ri
g_{\tX\1b}=g_b\1\tau\Ri g_{bb\1\tX\1b}=(g_b)\1\tau\Ri
(g_b)_{\tX_b\1}=(g_b)\1\cdot \tau.$
\item Let $Y$ be a half-twist adjacent to $X_b.$
Then $Y_{b\1}$ is adjacent to $X$ and satisfies axiom (2) of prime elements
for $g,$ $X$ and $Y_b\1.$
Namely: $g_{\tY_{b\1}\1\tX}=g\1g_{\tY_{b\1}\1}$ and
$g_{\tX\tY_{b\1}\1\tX\1}=g\1_{\tX}g_{\tX\tY_{b\1}\1}.$ \item"" (a) \
$g_{\tY_{b\1}\1\tX}=g\1g_{\tY_{b\1}\1}\Ri$ \newline

$(g_{\tT_{b\1}\1\tX})_b=(g\1g_{\tY_{b\1}\1})_b\Ri$\newline

$(g_b)_{(\tY_{b\1}\tX)_b}=(g_b)\1(g_b)_{\tY\1}\Ri$\newline

$(g_b)_{\tY\tX_b}=g_b\1\cdot(g_b)_{\tY\1}.$
\item"" (b)\ $g_{\tX\tY_{b\1}\1\tX\1}=g\1_{\tX}g_{\tX\tY_{b\1}\1}\Ri$
\newline

$(g_b)_{\tX_{b}\tY\1\tX_b\1}=(g_b\1)_{\tX_b}(g_b)_{\tX_b\tY\1}.$ \item"(3)" Let
$Z$ be a half-twist disjoint from $X_b.$ Then $Z_{b\1}$ is disjoint from $X.$
Then $g_{\tZ_{b\1}}=g.$
We conjugate $g_{\tZ_b\1}=g$ by $b$ by get:\
$(g_b)_{\tZ}=g_b.$\quad $\square$\endroster\edm

 We need the following lemma on $B_n$ to prove later a criterion for prime
element in a $\tB_n$-group.
\proclaim{Claim III.2.3}
Let $(X_1, \ldots, X_{n-1})$ be a frame in
$B_n = B_n[D,{K}]$.  Let $C(X_1) = \{b \in B_n|[b,X_1] =
1\}$ (centralizer of $X_1$), $C_p(X_1) = \{b \in
B_n\bigm|(X_{1})_{b} = X_1,$ preserving polarization\}. Let
 $\sigma = X_2X^2_1X_2$.
Then $C(X_1)$ is generated by $\{X_1,\sigma,X_3, \ldots, X_n\}$,
$C_p(X_1)$ is generated by   $\{X^2_1, \sigma, X_3, \ldots,
X_n\}$.
\ep
\demo{Proof of the Claim} Let $K=\{a_1,\dots,a_n\}.$  Let $x_1, \ldots ,
x_{n-1}$ be a system of
consecutive simple paths in $D,$ s.t. $X_i=H(x_i)$ ($H(x_i)$ is
the half-twist corresponding to $x_i;\
x_i$ connects $a_i$ with
$a_{i+1}$).  Let $\Gamma_1, \ldots, \Gamma_n$ be a free geometric base
of $\pi_1(D-K,*)$ consistent with $(X_1, \ldots, X_{n-1})$ (that
is, $(\Gamma_{i+1})X_i = \Gamma_i$, $(\Gamma_i)X_i =
\Gamma_i\Gamma_{i+1}\Gamma^{-1}_i$, $(\Gamma_j)X_i = \Gamma_j$
for $j \neq i$, $i + 1$).
We can assume that the $x_i$ does not
intersect the ``tails'' of $\Gamma_1, \ldots, \Gamma_n$.

Let $K_1$ be a finite set of $D$ obtained from $K\cup\{x_1\}$ by contracting
$x_1$ to a point $\tilde a_2\in x_1.$
$K_1=\{\tilde a_2,a_3,\dots, a_n\}.$
Let $B_{n-1}=B_{n-1}[D,K_1].$
Let $Y_2,\dots ,Y_{n-1}$ be a frame of $B_{n-1}$ where $Y_i$ can be
identified with $X_i$ for $i\ge 3.$

Let $H=\{b\in B_{n-1}\bigm| (\tilde a_2)b=\tilde a_2\}.$
{}From the short exact sequence $1\rightarrow P_{n-1}\hookrightarrow
B_{n-1} \rightarrow S_{n-1} \rightarrow 1$ (see \cite{MoTe 4}) we can conclude
that $H$ is generated by $Y_3,\dots ,Y_{n-1}$ and by the generators of
$P_{n-1}.$
We remove the generators of $P_{n-1}$ that can be expressed in terms of
$Y_3,\dots,Y_{n-1}$ (see \cite{A},\cite{B}, and \cite{MoTe4}, Section IV) and
conclude that $H$ is generated by $Y_2^2,Y_3,\dots ,Y_n.$
The element $Y_2^2$ corresponds to the motion $\Cal M'$ of $\tilde
a_2,a_3,\dots,a_n $ described as follows:
\ $\tilde a_2,a_4,\dots, a_n$ stays in place and $a_3$ is moving around
$\tilde a_2$ in the positive direction (see Fig. III.2.(a)).

\midspace{2.100in}\caption{Fig. III.2}

We  define a homomorphism $\Phi: C_p(X_1)\ri H$ as
follows:

Let $U$ be a ``narrow'' neighborhood of $x_1$ such that $\lambda =
\partial U$ is a simple loop.  Take $b \in C_p(X_1).$ There exists a
representing diffeomorphism $\beta:  D \rightarrow D$ ($\beta({K}) = {K}$,
$\beta|_{\partial D} = \operatorname{Id}_{\partial D})$ s.t.
$\beta|_{\overline{U}} = \operatorname{Id}_{\overline{U}}$
$(\overline{U} = U \cup \lambda)$.

The diffeomorphism $\be$ also defines an element of $B_{n-1}[D,K_1].$
This element is in fact in $H$ since $\tilde a_2\in x_1$ and thus
$(\tilde a_2)\be=\tilde a_2.$
Denote this element by $\Phi(b).$
The map $\Phi$ constructed in this way is obviously a homomorphism,
$\Phi: C_p(x_1)\ri H.$
Clearly, $X_3,\dots, X_{n-1}\in C_p(X_1).$
Clearly, $\Phi(X_i)=Y_i$ for $i\ge 3.$
Let $\Cal M$ be the following motion in $(D,K)$: \
$a_1,a_2,a_4,\dots,a_n$ are stationary and $a_3$ goes around $a_1,a_2$ in the
positive direction (Fig. III.2(b)).
Let $u$ be the braid in $C_p(X_1)$ induced from the motion in $\Cal M.$
Clearly, $\Phi(u)=Y_2^2.$
Thus, $\Phi$ is onto and $\Phi(u),\Phi(X_3),\dots,\Phi(X_{n-1})$ generate
$h.$
One can check that $u=Z_{31}^2Z_{32}^2.$
But $Z_{31}=X_2X_1A_1\1$ and $Z_{32}=X_2.$
Thus $u=\si.$
Thus $C_p(X_1)$ is generated by $\si,X_3,\dots, X_{n-1}$ and a set of
generators for $\ker\Phi.$

 Consider $\pi_1(D-K\cup x_1,*).$ Let $\tilde
\G_2$ be the path obtained from connecting $\lambda$ with $*\in \partial
D$ by a simple path intersecting each of $\Gamma_3, \ldots, \Gamma_{n }$
only at $*$. We get a (free) geometric base $\tilde{\Gamma}_2,\Gamma_3,
\ldots, \Gamma_n$ of $\pi_1 (D - ({K} \cup x_1),*)$.
It is obvious that
$\tilde\G_2=\G_1\G_2.$ \ $\Phi(b)$ defines in a natural way an
automorphism of $\pi(D_k\cup\{X_1\},*)$ s.t.
 $\Phi(b)$ does not
change the product $\tilde{\Gamma}_2 \Gamma_3 \ldots \Gamma_n$, and
$(\tilde{\Gamma}_2)\Phi(b)$ is a conjugate of $\tilde{\Gamma}_2$.

Consider now any $Z \in \ker\Phi$.
We have
$(\tilde{\Gamma}_2) Z = \tilde{\Gamma}_2$ $(\tilde{\Gamma}_2 =
\Gamma_1\Gamma_2)$, $(\Gamma_j)Z = \Gamma_j \; \forall j = 3,
\ldots , n$.
This implies that $Z$ can be represented by a
diffeomorphism which is the identity outside of $U$, that is, $Z =
X^l_1$, $l \in {\Bbb Z}$.  Since $Z \in C_p(X_1)$,
we get $l \equiv 0 \pmod 2$.

Thus, $C_p(X_1)$ is generated by $X^2_1,\sigma,
X_3, \ldots, X_{n-1}$.   Clearly, $C(X_1)$ is generated by
$ {C}_p(X_1)$ and $X_1$.  \qed\quad for the Claim\edm

\proclaim{Lemma III.2.4}
Let $\{\tX_1, \ldots, X_{n-1}\}$ be a frame in
$B_n$, $(\tX_1, \ldots, \tX_{n-1})$ their images in $\tilde{B}_n$.
Let $u \in G$ ($G$ is a $\tilde{B}_n$-group) be such that
\roster\item\  $u_{\tX^{-1}_1} = u^{-1} \tau$ with $\tau^2 = 1$, $\tau \in$
$\Center(G)$, $\tau_b = \tau \; \forall b \in \tilde{B}_n$;
\item"(2$_a$)" \ $u_{\tX^{-1}_2\tX^{-1}_1} = u^{-1}u_{\tX^{-1}_2}$;
\item"(2$_b$)" \
$u_{\tX_1\tX^{-1}_2\tX^{-1}_1} = u^{-1}_{\tX_1} u_{\tX_1\tX^{-1}_2}$;
\item"(3)"   $u_{\tX_j} = u \; \forall j = 3, \ldots, n - 1$.\endroster

Then $u$ is a prime element in $G$, and $\tX_1$ is a supporting
half-twist for $u$.
\ep

\demo{Proof}
Let $Z \in B_n$ be any half-twist disjoint from $X_1$, $\tZ$ be the
image of $Z$ in $\tilde{B}_n$.  $\exists b \in B_n$ such that
$(X_{1})_{b} = X_1$, $(X_{3})_{b} = Z$.  By Claim  III.2.3, $b$
belongs to the subgroup of $B_n$ generated by $X_1,X_3, \ldots, X_{n-1}$ and
$\sigma = X_2X^2_1X_2$.  Let $\tilde{b}$ and $\tilde{\sigma}$ be the images of
$b$ and $\sigma$ in $\tilde{B}_n$.   We have $u_{\tilde\si\1}=
u_{\tX^{-1}_2\tX^{-2}_1\tX^{-1}_2} = (u_{\tX^{-1}_2\tX^{-1}_1})_{\tX^{-
1}_1\tX^{-1}_2} = (u^{-1}u_{\tX^{-1}_2})_{\tX^{-1}_1\tX^{-1}_2} = (\tau
u)_{\tX^{-1}_2} \cdot u_{\tX^{-1}_2\tX^{-1}_1\tX^{-1}_2} = \tau u_{\tX^{-
1}_2} \cdot u_{\tX^{-1}_1\tX^{-1}_2\tX^{-1}_1} = \tau u_{\tX^{-1}_2}
\cdot (\tau u^{-1})_{\tX^{-1}_2\tX^{-1}_1} = \tau^2 u_{\tX^{-1}_2} u^{-
1}_{\tX^{-1}_2} \tau(\tau u) = u$.  Then $u_{\tilde\si}=u.$  Now: $u_{\tX_j}=u$
for $j\ge 3$ (by assumption (3)) and $u_{\tX_1^2}=u$ (by assumption (1)).
Thus, if $X_1$ appears in $b$ an even number of times, then $u_b=u.$ Otherwise
we replace $b$ by $bX_1.$ The ``new'' $b$ satisfies the same requirement for
$b,$ as above and $u_{\tilde b}=u.$ Thus, we can assume $u_{\tilde{b}} = u$.
We have $$u_{\tZ} = u_{\tilde{b}^{-1}\tX_3\tilde{b}} = u_{\tX_3\tilde{b}} =
u_{\tilde{b}} = u.$$

Let $Y$ be a half-twist in $B_n$ adjacent to $X_1$.  $\exists b_1
\in B_n$ s.t. $(X_{1})_{b_1} = X_1$, $(X_{2})_{b_1} = Y$.  Let
$\tilde{b}_1$ and $\tY$ be the images of $b_1$ and $Y$ in
$\tilde{B}_n$.  As above, we can choose $b_1$ so that
$u_{\tilde{b}_1} = u$.  Applying $\tilde{b}_1$ on the assumptions (2$_a$) and
(2$_b)$  we get (since $u_{\tilde{b}_1} = u$,
$(\tX_{1})_{\tilde{b}_1} = \tX_1$,\ $(\tX_{2})_{\tb_1} = \tY$):
$$u_{\tY^{-1}\tX^{-1}_1} = u^{-1}u_{\tY^{-1}} \ \text{and} \
u_{\tX_1\tY^{-1}\tX^{-1}_1} = u^{-1}_{\tX_1} u_{\tX_1\tY^{-1}}.\quad
\square$$ \edm

\newpage

\subheading{\S3. Polarized pairs and uniqueness of coherent pairs}

\demo{Definition}\ $\un{\text{Polarized pair}}$

Let $G$ be a $\tilde B_n$-group, $h$ a prime element of $G,$\ $X$
its supporting half-twist.
If $X$ is polarized, we say that $(h, X)$ (or $(h,\tX)$ is a polarized pair
with central element $\tau,$\quad $\tau=hh_{\tX\1}.$\edm

\demo{Definition}\ $\un{\text{Coherent pairs, Anti-coherent pairs}}$

We say that two polarized pairs $(h_1, \tX_1)$ and $(h_2, \tX_2)$ are coherent
(anti-coherent) if $\exists \tb \in \tB_n$ such that
$(h_1)_{\tb}=h_2,$\ $(\tX_1)_b=\tX_2,$ and $b$ preserves (reverses) the
polarization.\edm
\proclaim{Corollary III.3.1} Coherent and anti-coherent polarized pairs have
the same central element.\ep
\demo{Proof} The prime elements of coherent and anti-coherent pairs are
conjugate to each other.
Thus by Lemma III.2.2 we get the corollary.\quad $\square$\edm

We need the following Lemma to prove later the unique existence of a prime
element  with given s.h.t. conjugate to a given prime
element.

\proclaim{Lemma III.3.2}

Let $h\in G,$ a prime element with s.h.t. $\tX.$
Let $b\in B_n.$
Then: $\tX_{\tilde b}=\tX\Rightarrow h_{\tilde b}=h.$\ep

\demo{Proof}
  We can choose a set of
standard generators for $B_n[D,{K}],$ $\{X_1, \ldots,
X_{n-1}\}$ with $X_1 = X$.  Let $\sigma=X_2X_1^2X_2.$ Consider $C_p(X_1),$ the
centralizer of $X_1$ preserving polarization. By Lemma III.2.3, $C_p(X_1)$
is the subgroup of $B_n$ generated by $X^2_1, \sigma, X_3, \ldots,
X_{n-1}$.   Since  $\tX_3,
\ldots, \tX_{n-1}$   are disjoint from $X_1$, they do  not change $h$ (by
axiom(3) of prime elements).  By Lemma III.2.1, $h_{\tX_1^2}=h$.  Consider
 $h_{\sigma^{-1}} = h_{\tX^{-1}_2\tX^{-2}_1\tX^{-
1}_2}$.  We have:
$$\align
h_{\tilde\sigma^{-1}} &=
h_{\tX^{-1}_2\tX^{-2}_1\tX^{-1}_2}= (h_{\tX^{- 1}_2\tX^{-1}_1})_{\tX^{-1}_1
\tX^{-1}_2}\overset\text{by Axiom(2) of prime
element}\to  =
(h^{-1}h_{\tX^{- 1}_2})_{\tX^{-1}_1\tX^{- 1}_2} \\
  & =
(h^{-1}_{\tX^{-1}_1}h_{\tX^{-1}_2\tX^{-1}_1})_{\tX^{-1}_2}\overset\text{by
Axiom(2) of prime
element}\to =
(\tau h \cdot h^{-1}h_{\tX^{-1}_2})_{\tX^{-1}_2} = \tau h_{\tX^{-2}_2}\\
&\overset\text{by Lemma III.2.1(2) }\to
  =  \tau h \tau = h.\endalign$$
Thus $h_{\tilde\sigma} = h.$
Thus, for every generator $g$ of $C_p(X),$\ $h_{\tg}=h.$ Since
$ b \in \tilde{C}(X),$ \ $h_{\tB} = h$.\quad
$\square$\edm
\proclaim{Proposition III.3.3}
Let $\{h,\tX\}$ be a polarized pair, $h \in G$, $\tX \in
\tilde{B}_n$. Let $\tT$ be a polarized half-twist in $\tilde{B}_n$.
Then there exists a {\it unique} prime element
$g \in G$ such that $\{g,\tT\}$ and
$\{h,\tX\}$ are coherent.
\ep

\demo{Proof} Let $X,T \in B_n$ be polarized half-twists
representing $\tX$ and $\tT$.  $\exists b \in B_n$ such that $T =
X_b$ preserving polarization.  Let $\tilde{b}$ be the image of
$b$ in $\tilde{B}_n$.  Taking $g = h_{\tilde{b}}$ we
obtain a polarized pair $\{g,\tT\}$ such that
$\{g,\tT\}$ and $\{h,\tX\}$ are coherent.  To prove the {\it
uniqueness} of $g$, assume that
$\{g_1,\tT\}$ is another polarized pair coherent with $\{h,\tX\}$.  Then
$\exists b_1 \in B_n$ with $g = h_{\tilde{b}_1}$  and $T = X_{b_1}$,
preserving polarization.  We have $T = X_{b_1} = X_b$ and
$X_{b_1b^{-1}} = X$.  Denote  $b_2 = b_1b^{-1}$, so $X_{b_2} =
X$ (preserving  polarization).
By the previous lemma, $h_{b_2}=h.$
Thus, $h_{b_1}=h_b$ or $g=g_1.$\quad $\square$\edm
\demo{Definition} $\un{L_{h,\tX}(\tT)}$

Let $(h,\tX)$ be a polarized pair.
$\tT\in\tB_n.$
$L_{(h,\tX)}(\tT)$ is the unique prime element s.t. $(L_{(h,\tX)}(\tT),\tT)$ is
coherent
with $(h,\tX).$\edm
\mk
{}From uniqueness we get a simultaneous conjugation:
\proclaim{Lemma III.3.4}
Assume $(h,\tX) $ and $(g,\tX)$ are polarized pairs.  Let $\tau$ be the
central
element of $(g,\tX).$ If $(h,\tX)$ is anti-coherent to $(g,\tX)$ then
$h=g\1\cdot\tau.$\ep

\demo{Proof}
By assumption, $\exists b\in B_n$ s.t. $g=h_{\tilde b}$ and
$X=X_b,$ reversing polarization. Thus $X_{bX\1}=X,$ preserving
polarization.
Thus $(h_{\tilde b\tX\1},\tX)$ is coherent with $(h_{\tilde
b\tX\1},\tX_{\tilde b\tX\1}).$
Clearly, $(h,\tX)$ is coherent with $(h_{\tilde
b\tX\1},\tX_{\tilde b\tX\1}).$
{}From uniqueness, $h=h_{\tilde b\tX\1}=g_{\tX\1}.$
Since $\tau$ is the central element of $(g,\tX),$\ $\tau=g\cdot
g_{\tX\1}.$
Thus, $g_{X\1}=g\1\tau.$
So $h=g\1\tau.$\quad\qed\edm \proclaim{Corollary III.3.5} If
$(a_i,\tX)$ is coherent with $(g_i,\tY)$\ $i=1,2$, then there exist $b\in\tB$
s.t. $(a)_b=g_i$\ $i=1,2.$\ep \demo{Proof} Let $b$ be the element of $\tB_n$
s.t. $(a_1)_{b}=g_1,$\ $(\tX)_{b}  =\tY.$ Now, $((a_2)_{b},(X)_{b})$ is
coherent
with $(a_2,\tX).$ Since $(\tX)_b=\tY,$\ $((a_2)_b,\tY)$ is coherent with
$(a_2,\tX).$ The pair $(g_2,\tY)$ is also coherent with $(a_2,\tX).$ From
uniqueness, $(a_2)_{b}=g_2.$ \qed\edm \bk

\subheading{\S4. $\bold{\tB_n}$-action of nondisjoint half-twists}

\proclaim{Proposition III.4.1}
Let $\tT,\tY$ be $2$ orderly adjacent polarized half-twists in
$\tilde{B}_n$, $\{h,\tX\}$ be a polarized pair, $h \in G$, $\tX \in
\tilde{B}_n$.  Denote by $\tY^\prime$ the polarized half-twist
obtained from $\tY$ by changing polarization (that is, $\tT,
\tY^\prime$ are {\it  not orderly adjacent}).  Denote by
$$\gather L(T) = L_{\{h,\tX\}}(\tT), \\
L(Y) = L_{\{h,\tX\}}(\tY),\\
L(Y^\prime) = L_{\{h,\tX\}}(\tY^\prime).\endgather$$
Then
\roster\item $L(T)_{\tT^{-1}} = L(T)^{-1}\tau;$
\item  $L(T)_{\tY^{-1}} = L(T)L(Y)$;
\item  $L(T)_{(\tY'){}^{-1}} = L(Y)^{-1}L(T)$.\endroster
\ep

\demo{Proof}

(1) By Lemma III.2.1(1).

(2)  Let $b\in B_n$ be s.t. $L(T)=h_b,$\ $T=X_b,$ preserving polarization.
Let $Y_1 = Y_{{b}^{-1}}.$  Then $\{X,Y_1\}$\  is a pair of
adjacent half-twists $(X=T_{b\1},\ Y=Y_{b\1}),$ and so $h_{\tY^{-1}_1\tX^{-1}}
=
h^{-1}h_{\tY^{-1}_1}$.  Applying  $\tilde{b}$ to that equation, we get
$(L(T))_{\tY^{-1}\tT^{-1}} = L(T)^{-1} L(T)_{\tY^{-1}}$, or $$
L(T)_{\tY^{-1}} = L(T) \cdot L(T)_{\tY^{-1}\tT^{-1}}
$$

Let $b_1 = bY^{-1}T^{-1}$.  Then $X_{b_1} = X_{bY^{-1}T^{-1}} =
T_{Y^{-1}T^{-1}} = Y$ (since $T_{Y\1}=Y_T$ by III.0).  Using that $T,Y$ are
orderly adjacent and $X_b = T,$ preserving polarization, one can easily check
that actually $X_{b_1} = Y,$ preserving polarization.  Because of the {\it
uniqueness} of $L(Y) = L_{\{h,\tX\}}(\tY)$, we get
$L(Y) = h_{\tilde{b}_1} = h_{\tilde{b}\tY^{-1}\tT^{-1}} =
L(T)_{\tY^{-1}\tT^{-1}}$.  Together with the previous equation we get:
$L(T)_{\tY^{-1}} = L(T)L(Y)$, which is (2).

(3) Using $Y_{Y^{-1}} = Y'$ (preserving polarization) and
uniqueness, we can write
$L(Y') = L(Y)_{\tY^{-1}}.$ \ By (1), $L(Y)_{\tY\1} = L(Y)^{-1}
\tau.$ Thus, $L(Y') = L(Y)^{-1} \tau$.  From (2) we get:
$  L(T)_{\tY'{}^{-1}} = L(T)L(Y') =
L(T)L(Y)^{-1}\tau = L(Y)^{-1}L(T)$, which is (3).

(We used $[L(T),L(T)_{\tY^{-1}}] = \tau,$ from Lemma III.2.1, which implies
\newline $\tau = [L(T),L(T)L(Y)] = [L(T),L(Y)]$. ) \qed  \edm

\proclaim{Lemma III.4.2}
Let $h$ be a prime element in $G$,\ $\tX \in \tilde{B}_n$ a
supporting half-twist of $h$, $\tZ$ a half-twist in $\tilde{B}_n$
transversal to $\tX$.  Then $h_{\tZ} = h$.\ep

\demo{Proof} Let $X,Z$ be transversal half-twists in $B_n$,
representing $\tX,\tZ.$  Let $x,z$ be 2 transversally intersecting
simple paths corresponding to $X,Z$ (see Fig. III.4).

\midspace{2.00in}
\caption{Fig. III.4}

There exists a simple path $y$ such that the corresponding half-twist $Y$ is
adjacent to $X$ and $Z$, and $Z_1 = Z_{Y^{-2}}$ is disjoint from $X.$ Let $z_1$
be the path corresponding to $Z_1$ (see Fig. III.4).  Denote by $\tY, \tZ_1$
the images of $Y,Z_1$ in $\tilde{B}_n$.  We have $h_{\tZ} =
h_{\tY^{-2}\tZ_1\tY^2} \overset\text{by Lemma III.2.1}\to= (h
\tau)_{\tZ_1\tY^2} = (h \tau)_{\tY^2} = h \tau \cdot \tau = h$.  \qed\quad\edm

\bigskip

\subheading{\S5. Commutativity properties}

\proclaim{Proposition III.5.1}
Let $\{g_1,\tY_1\}, \{g_2,\tY_2\}$ be $2$ polarized pairs in $G$.
Assume that they are coherent or anti-coherent.
Let $\tau$ be the corresponding central element of $(g_1,\tY_1)$
\ $(\tau = g_1(g_{1})_{\tY^{-1}_1})$.

Then

\roster\item  if $\tY_1,\tY_2$ are adjacent, then $[g_1,g_2] = \tau$;
\item  if $\tY_1,\tY_2$ are disjoint or transversal, then $[g_1,g_2] =
1$.\endroster
\ep

\demo{Proof}

(1)  Assume first that $\{g_1,\tY_1\},\{g_2,\tY_2\}$ are coherent.
Take $b \in \tilde{B}_n$ with $g_2 = (g_1)_{b}$, $\tY_2 = (\tY_1)_{b}$
(preserving polarization).  Let $b_1 = \tY^{-1}_2\tY^{-1}_1$.  Then
$(\tY_1)_{b_1} = \tY_2$.  Assume that $b_1$ preserves polarization of
$\tY_1,\tY_2$.  We have $\{(g_{1})_{b_1},\tY_2\}$ and $\{g_2,\tY_2\}$
coherent with $\{g_1,\tY_1\}$.  By Proposition III.3.3 (the {\it
uniqueness} part) we get $(g_1)_{b_1} = g_2$.  Thus we have $g_2 =
(g_1)_{b_1} = (g_1)_{\tY^{-1}_2\tY^{-1}_1} = g^{-1}_1
(g_{1})_{\tY^{-1}_2}$, and $[g_1,g_2] =
[g_1,g^{-1}_1(g_{1})_{\tY^{-1}_2}] = g_1g^{-1}_1(g_{1})_{\tY^{-1}_2}
g^{-1}_1(g^{-1}_{1})_{\tY^{-1}_2}g_1 = [g_1,(g_{1})_{\tY^{-1}_2}]_{g_1}
\overset\text{by Lemma III.2.1(3)}\to=  \tau_{g_1} = \tau$.

If $b_1$ does not preserve the polarization of $\tY_1,\tY_2,$ consider $b_2 =
b_1\tY_2$.  Then $\tY_2 = (\tY_1)_{b_2}$, preserving polarization.  As above,
we
get $g_2 = (g_1)_{b_2} = (g_1)_{\tY
^{-1}_2 \tY^{-1}_1\tY_2} = (g_1)_{\tY_1\tY^{-
1}_2\tY^{-1}_1} = (g^{-1}_1 \tau)_{\tY^{-1}_2\tY^{-1}_1} = \tau (g^{-
1}_{1})_{\tY^{-1}_2}g_1$, and then $[g_1,g_2] = [g_1,\tau (g^{-
1}_{1})_{\tY^{-1}_2}g_1] = [g_1,(g^{-
1}_{1})_{\tY^{-1}_2}] = [g_1,(g_1)_{\tY^{-
1}_2}]^{-1}_{(g_1)_{\tY^{-1}_2}} \overset\text{by Lemma III.2.1(3)}\to= \tau$.

If $\{g_1,\tY_1\}$, $\{g_2,\tY_2\}$ are anti-coherent, denote $\tY'_2$
the half-twist obtained from $\tY_2$ by changing polarization.  One
can then check that $\{g_1,\tY_1\}$, $\{g_{2;\tY^{-1}_2},\tY'_2\}$ are
coherent.  We have from the above that $\tau = [g_1,(g_{2})_{\tY^{-1}_2}].$
By Corollary III.3.1 $\tau$ is also the central element of $(g_2,\tY_2).$
Thus $\tau=g_2(g_2)_{\tY_2\1}$ which implies $(g_2)_{\tY_2\1}=g_2\1\tau.$
Thus $\tau=[g_1,(g_2)_{\tY_2\1}] =
[g_1,g^{-1}_2 \tau] = [g_1,g_2^{-1}] = [g_1,g_2]^{-1}_{g_2}$.  Thus, $[g_1,g_2]
= \tau^{-1}_{g_2} = \tau$.

(2)  We can assume that $\{g_1,\tY_1\},\{g_2,\tY_2\}$ are coherent.
(Otherwise, we replace $\tY_2$ by $\tY'_2$ and $g_2$ by $(g_{2})_{\tY^{-
1}_2}$ and use $[g_1,(g_{2}){\tY^{-1}_2}] = [g_1,g_2]^{-1}_{g_2}$ (see
above).)

Consider first the case where $\tY_1,\tY_2$ are disjoint.  We can
choose a standard base of $\tilde{B}_n$, say $(\tX_1,\tX_2, \ldots,
\tX_{n-1})$ such that $\tX_1 = \tY_1$, $\tX_3 = \tY_2$ and the given
polarizations of $\tY_1,\tY_2$ coincide with ``consecutive''
polarizations of $\tX_1,\tX_3$ (``end'' of $\tX_1$ = ``origin'' of
$\tX_2$, ``end'' of $\tX_2 =$ ``origin'' of $\tX_3$).  Let $b_1 = \tX^{-
1}_2\tX^{-1}_1\tX^{-1}_3\tX^{-1}_2$.  Then $\tY_2 = (\tY_1)_{b_1}$, preserving
polarization.  From Proposition III.3.3 ({\it uniqueness}) it follows that
$g_2 = (g_1)_{b_1} = (g_1)_{\tX^{-1}_2\tX^{-1}_1\tX^{-1}_3\tX^{-1}_2} =
(g^{-1}_1(g_1)_{\tX^{-1}_2})_{\tX^{-1}_3\tX^{-1}_2} = (g^{-1}_{1})_{\tX^{-
1}_2}(g_1)_{\tX^{-1}_2\tX^{-1}_3\tX^{-1}_2} = (g^{-1}_{1})_{\tX^{-
1}_2}(g_1)_{\tX^{-1}_3\tX^{-1}_2\tX^{-1}_3} = (g^{-1}_{1})_{\tX^{-
1}_2}(g_1)_{\tX^{-1}_2\tX^{-1}_3}$.  By III.2.1(3)
$[g_1,(g_{1})_{\tX^{-1}_2}] = \tau$, which implies
$[g_1,(g_{1})_{\tX^{-1}_2\tX^{-1}_3}] = [g_1,(g_{1})_{\tX^{-
1}_2}]_{\tX^{-1}_3} = \tau_{\tX^{-1}_3}  = \tau$.  We can write $[g_1,g_2]
= [g_1,(g^{-1}_{1})_{\tX^{-1}_2}(g_1)_{\tX^{-1}_2\tX^{-1}_3}] =
[g_1,(g^{-1}_{1})_{\tX^{-1}_2}] \cdot [g_1,(g_1)_{\tX^{-1}_2\tX^{-
1}_3}]_{(g_1)_{\tX^{-1}_2}} = \tau \cdot \tau = \tau^2 = 1$.

Assume now that $\tY_1,\tY_2$ are transversal.  As in the proof of
Lemma III.4.2, we can find a half-twist $\tT \in \tilde{B}_n$ such that
$\tT$ is adjacent to $\tY_1,\tY_2$ and $ \tY_2' = (\tY_{2})_{\tT^{-2}}$ is
disjoint from $\tY_1$.  Let $b \in \tilde{B}_n$ be such that $\tY_2 =
(\tY_1)_{b}$, $g_2 = (g_1)_{b}$.  Let ${b}' = b\tT^{-2}$,\
 ${g}_2' = (g_1)_{b'} = (g_{2})_{\tT_2^{-2}}$.  Then
$\{g_2',\tY_2'\}$ is coherent, or anti-coherent, with
$\{g_1,\tY_1\}$.  Since $\tilde{Y}_2',\tY_1$ are disjoint, we get from
 the above $[g_1,g_2'] = 1$.  By Lemma III.2.1 $ g_2' =
(g_2)_{\tT^{-2}} = g_2\tau$, or $g_2 = g_2'\tau$.  Therefore,
$[g_1,g_2] = [g_1, g_2' \tau] = [g_1,  g_2'] = 1$. \quad
\qed\edm

Recall that here exists a natural homomorphism $\psi_n: B_n\ri S_n.$
$\psi_n(X_i)$ is the transposition $(i\ i+1)$ for $X_i$ a half-twist
connecting the points $q_i$ and $q_{i+1}.$

\demo{Definition}\ $\un{P_n}$

$P_n=\ker\psi_n.$

Recall from \cite{MoTe4} that: $P_n$ is generated by $Z_{ij}^2,$ where:

$Z_{ij}=(X_i^2)_{X_{i+1}\dots X_{j-1}}.$\edm

\demo{Definition}\ $\un{\tilde{P}_n}$

$ \tilde{P}_n =
\ker(\tilde{B}_n\overset \tilde\psi_n\to \rightarrow S_n)$ where $\psi_n$ is
induced naturally from $\psi_n.$\edm

\proclaim{Proposition III.5.2}
Assume $n \geq 4$. Let $\tX_1,\tX_2$ be $2$ adjacent half-twists in $\tB_n.$
Let $c=[\tX_1^2,\tX_2^2].$  Then the commutant $\tilde{P}^\prime_n$ of
$\tilde{P}_n$ is generated by $c$ where $c_b = c \;
\forall b \in \tilde{B}_n$, and $c^2=1.$
Moreover, if $(\tY_1,\tY_2)$ and $(\tZ_1,\tZ_2)$ are two pairs of adjacent
half-twists, then
$[\tZ_1^2,\tZ_2^2]=[\tY_1^2,\tY_2^2]=[\tZ_1^2,\tZ_2\2]=[\tZ_1\2,\tZ_2\2]=c.$
\ep

\pf Let $B_n = B_n(D, {K})$). Complete $\tX_1$ and $\tX_2$ to
 $\tX_1, \ldots, \tX_{n-1},$ a standard base of $\tilde{B}_n$,
\ $X_i=(H(x_i)$ and $x_1,\dots,x_{n-1}$ are simple paths
  in $D.$  Let $c=[\tX_1^2,\tX_2^2].$
 Let $x = (x_{1})_{\tX_2\tX_3}$.
We have a quadrangle formed by $x_1,x_2,x_3,x$, (see
Fig. III.5(a)).

\midspace{2.00in}
\caption{Fig. III.5}

Denote by $X \in {B}_n$ the half-twist defined by $x$.
Evidently, $\tX_1,\tX_2,\tX_3,\tX$ form a good quadrangle in $\tB_n.$  Thus
by Lemma III.1.2
$$
\tX^2_1\tX^2_3 = \tX^2_2\tX^2.
\tag1.10$$
Denote  $y_1 = \tX^2_1$,\quad $y_2 = \tX^2_2$,\quad $y_3 = \tX^2_3$,\quad  $y_4
= \tX^2,$\ (the squares of the edges),\newline $d_1 = \tX_1\tX^2_2\tX^{-1}_1$,\
$d_2 = \tX_2\tX_3^2\tX^{-1}_2$,\  (the squares of
the diagonals), \newline $y'= \tX_2\tX_1^2\tX^{-1}_2$ (the square of the outer
diagonal)
\flushpar See Fig. III.5(b) where we denote the paths corresponding to the
half-twists whose squares we considered here.

Clearly:

$y_1y_3=y_3y_1$

$d_1=(y')_{y_1\1}$

$(y_3)_{x_2}=(d_2)_{x_3\2}=(d_2)_{y_3\1}$

$y_1y_2d_1=y_2y_1y'=\Delta_3^2$ (a central element of $P_3).$

We rewrite (1.10) to get
$$y_4=y_1y_3y_2\1
\tag1.11$$

Conjugating (1.11) by $\tX_1$, we get:
$$
d_2 = y_1y_3y'{}^{-1};
\tag1.12$$
conjugating (1.11) by $\tX_2$ we get
$y_4 = d_1(d_2)_{y_3\1}\cdot y_2\1.$
Since $d_1 =
y_1y'y^{-1}_1$,

$$y_4 = y_1 y'y^{-1}_1 y_3d_2y^{-
1}_3y^{-1}_2 \overset\text{by (1.12)}\to=  y_1y' y^{-1}_1 y_3 y_1
y_3 y'{}^{-1} y_3^{-1}y^{-1}_2 .$$
We compare the last expression with (1.11) to get:$$y'y^2_3 y'{}^{-1}y^{-1}_3 =
y_3, \; \quad\text{or}\quad  [y',y^2_3] = 1.
\tag1.13$$

Since $y',y_3$ are squares of two adjacent half-twists in
$\tilde{B}_n$, and any two pairs of adjacent half-twists are
conjugate, we conclude from (1.13) that:
$$\forall\ \text{pairs of
adjacent half-twists, say}\ \tZ_1,\tZ_2\ \text{in}\ \tilde{B}_n:
[\tZ^2_1,\tZ^4_2] = 1,
\tag1.14$$
which also implies that $[\tZ^2_1,\tZ^2_2] = [\tZ^{-2}_1,\tZ^2_2] =
[\tZ^2_1,\tZ^{-2}_2] = [\tZ^{-2}_1,\tZ^{-2}_2].$

Conjugating (1.11) by $\tX^{-1}_3$ we get:
$$d_1=y_1y_3(d_2)_{y_3\1} = y_1y_3\cdot y_3d^{-1}_2y^{-1}_3 = y_1y^2_3 \cdot
d^{- 1}_2y^{-1}_3\overset (1.14)\to = y_1d^{-1}_2y_3,$$

We have by (1.14) that $1 = [y^2_1,y_2] = [y_1,y_2]_{y^{-1}_1} \cdot
[y_1,y_2]$.  Denoting $c = [y_1,y_2]$, we can write:
$$
c_{y_1^{-1}} = c^{-1},\quad\text{or}\quad c_{y_1} = c^{-1}.
\tag 1.15$$
  Denote by $\tilde{P}_3$ the
subgroup of $\tilde{B}_n$ generated by $y_1,y_2,d_1$, and by
$\alpha = \Delta^2_3 = y_1y_2d_1 = y_2y_1y'$ (a central
element of $\tilde{P}_3$), so that
$$
y' = y_1^{-1}y^{-1}_2 \alpha.
\tag 1.16$$
So, $c_{\tX_1} = [y_1,y_2]_{\tX_1} = [y_1,y']
\overset\text{by (1.16)}\to=  [y_1,y^{-1}_1y^{-1}_2 \alpha] = y_1 \cdot
y^{-1}_1y^{-1}_2 \alpha \cdot y^{-1}_1 \cdot \alpha^{-1}y_2y_1 =
y^{-1}_2y^{-1}_1y_2y_1 = [y^{-1}_2,y^{-1}_1] = [y_2,y_1] = c^{-
1}.$
Thus we have $c_{\tX^2_1} = (c^{-1})_{\tX_1} = c.$

 By (1.15)
$c_{\tX^2_1} = c_{y_1} = c^{-1}$.

We compare the last two results to get\ $c = c^{-1}$ or
$$
c^2 = 1 \quad\text{and} \quad c_{\tX_1} = c.
\tag1.17$$

Using a conjugation which sends $(\tX_1,\tX_2)$ to $(\tX_2,\tX_1)$, we
obtain from (1.17)
$$
c^{-1}_{\tX_2} = c^{-1},\quad{or} \quad c_{\tX_2} = c.
\tag 1.18$$

(1.17) and (1.18) show that $\forall$  $z \in \tilde{B}_3$ (the
subgroup of $\tilde{B}_n$ generated by $\tX_1,\tX_2$) we have
$$
c_z = c.
\tag1.19$$

Consider now $c_{\tX_3} = [y_1,y_2]_{\tX_3} = [y_1,d_2]\overset\text{by
(1.12)}\to =
[y_1,y_1y_3y'{}^{-1}] = $\newline
$y_1y_1y_3y'{}^{-1}y^{-1}_1y' y^{-1}_3y^{-1}_1 = y_1y_3 \cdot
(y_1y'{}^{-1}y^{-1}_1 y')y^{-1}_3y_1^{-1}\overset\text{by (1.19)}\to=
y_1 y_3 c y^{-1}_3y^{-1}_1 = c_{y^{-1}_1y^{-1}_3} = c_{y^{-3}_1} = c_{y_3} =
c_{\tX^2_3}$, shortly $c_{\tX^2_3} = c_{\tX_3}$.  This implies $c_{\tX_3} = c$.

Since $c = [\tX^2_1,\tX^2_2]$, we have $\forall \tX_j$, $j \geq 4$,
$c_{\tX_j} = c$.  Thus $\forall b \in \tilde{B}_n$ $c_b = c$ and
$c^2 = 1$.

Let $(\tY_1,\tY_2)$ be a pair of adjacent half-twists.
Since every 2 pairs of adjacent half-twists are conjugate in $\tB_n,$\quad
$\exists b\in \tB_n$ s.t. $[\tY_1^2,\tY_2^2]=[\tX_1,\tX_2^2]_b=c_b.$
Since\linebreak $c_b=c\quad\fa b\in\tB_n,$\ $[\tY_1^2,\tY_2^2]=c.$
Since $c^2=1,$\ $c\in\Center(\tB_n)$ we also have
$[\tY_1^2,\tY_2\2]=[\tY_1\2,\tY_2\2]=c.$  In particular, if $(\tZ_1,\tZ_2)$ is
another pair of adjacent half-twists, $[\tZ_1 ^2,\tZ_2^2]=[\tY_1^2,\tY_2^2]=c.$
Because any two disjoint and transversal half-twists of $\tilde{B}_n$ commute,
and $\tilde{P}_n$ is generated by $\tZ^2_{ij} = (\tX^2_{i})_{\tX_{i+1} \cdots
\tX_{j-1}},$ $1 \leq i < j \leq n$ (see \cite{MoTe4}), we conclude that
$\tilde{P}^\prime_n$ is generated by $c$.  \quad \qed\enddemo

\bk
\subheading{\S6. $\bold{\tilde P_n}$ as a $\bold{\tilde B_n}$-group}

Recall: \ $\Ab(B_n)\simeq \Bbb Z$\ $(B_n$ is generated by the half-twists  and
every 2 half-twists are conjugate).
 \demo{Definition}\ $\un{P_{n,0}}$

 $P_{n,0} = \ker(P_n \rightarrow \Ab B_n)$ (``degree zero'' pure braids).
\edm

\demo{Definition}\ $\un{\tilde{P}_{n,0}}$

$\tilde{P}_{n,0}$ is
the image of $P_{n,0}$ in $\tilde{P}_n$.  \edm

\proclaim{Lemma III.6.1} Let $X_1,X_2$ be $2$ consecutive half-twists in $B_n.$
Let\newline $u = (\tX^2_{1})_{\tX^{- 1}_2} \tX^{-2}_2$. Then $u\in\tP_{n,0}$ \
$u$ is a prime element in $\tilde{P}_n$ (considered as a
$\tilde{B}_n$-group), and $\tX_1$ is the supporting half-twist of
$u$.
\ep

\demo{Proof} Clearly, $u \in \tilde{P}_{n,0}$.
Since $X_1X_2X_1=X_2X_1X_2,$ \ $(X_1)_{X_2\1}=(X_2)_{X_1}$ and thus,
$u=(\tX_1)_{\tX_2\1}^2\tX_2\2=(\tX_2)_{\tX_1}^2\tX_2\2.$
We often use here the fact that $(X_1)_{X_2\1}=(X_2)_{X_1}$ as well as the fact
that $[\tX_1^{\pm 2},X_2^{\pm 2}]=c,$ i.e., $\tX_2\2\tX_2^2=c\tX_2^2\tX_1\2$
and $(\tX_1^2)_{X_2\1}=(\tX_1^2)_{X_2}c$\ for $c\in\Center(\tB_n),$\ $c^2=1.$
Complete $X_1,X_2$ to a frame of $B_n:$\ $X_1,\dots,X_{n-1}.$ $(\la
X_i,X_{i+1}\ra=1$ and $[X_i,X_j]=1$\ $|i-j|>2.).$ We shall use Lemma III.2.4,
that is, we must check  conditions (1), $(2_a)$, $(2_b)$, (3) of Lemma
III.2.4.

(1)  We have $u_{\tX^{-1}_1} = (\tX^2_{2})_{\tX_1\tX^{-1}_1} \cdot
(\tX^{-
2}_{2})_{\tX^{-1}_1} =\tX^{2}_{2}\cdot (\tX_1\2)_{\tX_2}.$
 Since $c = [\tX^2_1,\tX^2_2]$
(see Proposition III.5.2), $(\tX_1\2)_{\tX_2}=(\tX_1\2)_{\tX_2\1} c.$ Thus:
$u_{\tX_1\1}=  \tX^2_2 \cdot (\tX^{-2}_{1})_{\tX^{-1}_2} c = u^{-1}c$.

$(2_a)$ Since $[\tX_1^2,\tX_2^2]=c,$\ $u_{\tX^{-1}_2} =
(\tX^2_{1})_{\tX^{-2}_2}
\cdot \tX^{-2}_2 = c\tX^2_1\tX^{-2}_2 = \tX^{-2}_2 \cdot \tX^2_1$, and
$u_{\tX^{-1}_2\tX^{-1}_1} = (\tX^{-2}_{2})_{\tX^{-1}_1}\cdot \tX^2_1 =
(\tX^{-2}_{1})_{\tX_2} \cdot \tX^2_1 = c(\tX^{-
2}_{1})_{\tX^{-1}_2} \cdot \tX^2_1.$

On the other hand, $u^{-1}u_{\tX^{-1}_2} = \tX^2_2
\cdot (\tX^{-2}_{1})_{\tX^{-1}_2} \cdot \tX^{-2}_2 \cdot \tX^2_1 = c(\tX^{-
2}_{1})_{\tX^{-1}_2} \cdot \tX^2_1.$ We get
$$u_{\tX^{-1}_2\tX^{-1}_1} = u^{-1}u_{\tX^{-1}_2}.\qquad\quad$$

$(2_b)$  Using (1) and $(2_a)$, we can write
$$\align&u_{\tX_1\tX^{-1}_2} = (u^{-1}c)_{\tX^{-1}_2} = u^{-1}_{\tX^{-1}_2}
\cdot c,\\&u_{\tX_1\tX^{-1}_2\tX^{-1}_1} = (u^{-1}c)_{\tX^{-1}_2\tX^{-1}_1} =
u^{- 1}_{\tX^{-1}_2} u \cdot c,\\&u^{-1}_{\tX_1} u_{\tX_1\tX^{-1}_2} = cu \cdot
u^{-1}_{\tX^{-1}_2} c = u^{-1}_{\tX^{-1}_2} u \cdot c\endalign$$
(we use $[u,u_{\tX^{-1}_2}] = [(\tX^2_{1})_{\tX^{-1}_2} \cdot \tX^{-2}_2,
\tX^{-2}_2\tX^2_1] = c \cdot c \cdot c = c)$.   Thus,
$$u_{\tX_1\tX^{-1}_2\tX^{-1}_1} = u^{-1}_{\tX_1} u_{\tX_1\tX^{-1}_2}.\qquad$$

(3)  Clearly, $\forall j \geq 4$, $u_{\tX_j} = u$.  Consider
$u_{\tX^{-1}_3} = (\tX^2_{1})_{\tX^{-1}_2\tX^{-1}_3}
(\tX^{-2}_{2})_{\tX^{-1}_3}$.

Since $u$ can also be written as $u = (\tX^2_{1})_{\tX^{-1}_2} \cdot \tX^{-2}_2
= \tX^{-2}_2 \cdot (\tX^2_{1})_{\tX^{-1}_2} \cdot c = \tX^{-2}_2 \cdot
(\tX^2_{1})_{\tX_2}$, we have:  $$\align u_{\tX^{-1}_3} &= u \Leftrightarrow
(\tX^2_{1})_{\tX^{-1}_2\tX^{-1}_3} \cdot( \tX^{-2}_{2})_{\tX^{-1}_3} =
\tX^{-2}_2 \cdot( \tX^2_{1})_{\tX_2} \Leftrightarrow \tX^2_2 \cdot
(\tX^2_{1})_{\tX^{-1}_2\tX^{-1}_3}\\& = (\tX^2_{1})_{\tX_2} \cdot
(\tX^2_{2})_{\tX^{-1}_3}\endalign$$ which is true, because
$\{(\tX_{1})_{\tX_2},
\tX_2, (\tX_{2})_{\tX^{-1}_3}, (\tX_{1})_{\tX^{-1}_2\tX^{-1}_3}\}$ form a good
quadrangle (see Fig. III.6).   \quad\qed\edm

\midspace{2.00in}
\caption{
Fig. III.6}

\subheading{Construction of $\bold{\underline G(n)}$}

For $n \geq 3$ we define the group $\underline{G}(n)$ as follows:

Generators:  $s_1,u_1,u_2, \ldots, u_{n-1}$.

Relations:  $$\align&[s_1,u_i] = 1\quad \forall i = 1,3, \ldots n-1;
[u_i,u_j] = 1\ \text{when}\ |i-j| \geq 2;\\
&[s_1,u_2] = [u_i,u_{i+1}] = [u_1,u_2] \; \forall i = 2,3,
\ldots n - 2;\\
&[u_1,u_2] = [u_1,u_2]_{s_1} = [u_1,u_2]_{u_i} \; \forall i =
1,2, \ldots, n - 1;\\
&[u_1,u_2]^2 = 1.\endalign$$

\subheading{Equivalent construction of $\bold{\underline G(n)}$}

Consider a free abelian group $A(n)$ with generators
$S_1,V_1, \ldots, V_{n-1}$ and a skew-symmetric ${\Bbb Z}/2$-
valued bilinear form $Q(x,y)$ on $A(n)$ defined by: $Q(S_1,V_i) =
0 \; \forall i = 1,3, \ldots, n - 1;$ $Q(V_i,V_j) = 0$ when $|i-
j| \geq 2$, $Q(S_1,V_2) = Q(V_i,V_{i+1}) = 1$ $\forall i =
1,2,\ldots, n -2$.
One can check that there exists a unique central extension $G$ of
$A(n)$ by ${\Bbb Z}/2$ with $\operatorname{Ab}(G) \simeq A(n)$, $G^\prime
\simeq {\Bbb Z}/2$ and such that $\forall x,y \in G [x,y] =
Q(\overline{x},\overline{y})$ where $(\overline{x}$ and $\overline{y}$ are
the images of $x,y$ in $A(n)$).
\proclaim{Claim III.6.2}
\roster\item
The above  central extension is isomorphic to  $\underline{G}(n)$.
\item $\Ab(\underline{G}(n))$ is a free abelian group with $n$
generators (i.e., $A(n)$) and $\underline G(n)^\prime \simeq {\Bbb Z}/2$
(generated by $[u_1,u_2])$.
\item
The following formulas define a
$\tilde{B}_n$-action on $\underline{G}(n)$ for  $(\tX_1, \ldots,
\tX_{n-1}),$   a standard set of generators in $\tilde{B}_n$, and
$\nu=[u_1,u_2].$  $$\alignat3
&\underline{\tX_1\text{-action}}\quad &&
\underline{\tX_2\text{-action}}\quad && \underline{\tX_k\text{-action}, k
\geq 3} \\
&s_1 \rightarrow s_1 && s_1 \rightarrow u_2s_1&&
s_1 \rightarrow s_1 \\
&u_1 \rightarrow u^{-1}_1 \nu && u_1 \rightarrow
u_2u_1&& u_{k-1} \rightarrow u_ku_{k-1}
\\& u_2
\rightarrow u_1u_2 && u_2 \rightarrow u^{-1}_2 \nu
&& u_k \rightarrow u^{-1}_k \nu\\
&u_j \rightarrow u_j \; \forall j \geq 3\quad && u_3 \rightarrow
u_2u_3; && u_{k+1} \rightarrow u_ku_{k+1}\\
& &&u_j \rightarrow u_j \ \forall j \geq 4\quad && u_j
\rightarrow u_j \ \forall j \neq k-1,k,k+1
\endalignat$$
 \item
 Let $b \in \tilde{B}_n$, $y = (\tX_{1})_{b}$.  Then the $y^2$-action on
$\underline{G}(n)$ coincides with the conjugation by $(s_{1})_{b}$.
\item Let
$$ s_{ij}=\cases
(s_1)\
\qquad\quad\quad\qquad\qquad\qquad\quad\qquad\qquad\qquad\qquad\qquad\qquad
\text{if}\ (i,j)=(1,2); \\  (s_{1})_{\tX_2\dots \tX_{j-1}}\quad
\left(\overset\text{Claim III.6.2(3)}\to=u_{j-1}\dots u_2s_1\right)\qquad\quad\
\ \ \ \quad \text{if}\ i=1,\ j\geq 2;\\ (s_1)_{\tX_2\dots \tX_{j-1}\tX_1\dots
\tX_{i-1}}=\cases \nu\cdot u_{j-1}\dots u_1\cdot u_{i-1}\dots u_2s_1&\
\text{if}\ i\geq 3,\ j>i\\ \nu\cdot u_{j-1}\dots u_1s_1&\ \text{if}\ i=2,\
j>i.\endcases\endcases\tag1.20$$ Then:   $$[s_{ij},s_{kl}] =
\cases
\nu, \quad& \text{if} \ (\{i,j\} \cap \{k,l\}) = 1);\\
1,\quad& \text{otherwise}.
\endcases$$
\item Let $\tilde F_{n-1}$ be the subgroup
  of $\underline{G}(n)$ generated by $(s_{n-
1,n},s_{n-2,n}, \ldots, s_{1n})$.  $\tB_{n-1}$ acts on $\tilde F_{n-1}$ as
follows: $$\align   (s_{jn})\tX_k &=
s_{jn} \qquad   j \neq k, k + 1; \quad k=1 \ldots n - 1\tag1.21\\
(s_{k,n})\tX_k &= s_{k+1,n} \\
(s_{k+1,n})\tX_k &= s_{kn} \nu = s_{k+1,n} s_{kn}s^{-1}_{k+1,n}.\endalign$$
This action is a $\tB_{n-1}$-action, where the action of the
generators $\tX_{n-2} ,\ldots \tX_1$ of $\tB_{n-1}$ correspond to
standard {\it Hurwitz moves} on ($s_{n-1,n},s_{n-1,n}, \ldots,
s_{1n})$. (See the definition of Hurwitz moves in Chapter $0$.)
\item There is a natural chain of embeddings $\underline{G}(3)
\subset \underline{G}(4) \subset \cdots \subset \underline{G}(n-
1) \subset \underline{G}(n)$ corresponding to the chain:
$(s_1,u_1,u_2) \subset (s_1,u_1,u_2,u_3) \subset \ldots \subset
(s_1,u_1, \ldots u_{n-1})$.\endroster
\endproclaim

\pf

(1), (2), and (3) are easy to verify.

(4)
Consider first the case $b =$ Id.  From (3) we get for the
$\tX^2_1$-action:
$$s_1 \rightarrow s_1, \; u_1 \rightarrow u_1, \; u_2 \rightarrow
u_2\nu, \; u_j \rightarrow u_j \; \forall j \geq 3.$$
At the same time by the first construction:
$$(s_{1})_{s_1} = s_1,\ (u_1)_{s_1} = u_1,\; (u_2)_{s_1} = s^{-
1}_1u_2s_1 = u_2 \nu,\; (u_{j})_{s_1} = u_j \; \forall j \geq 3.$$
Thus $\tX^2_1$-action and $s_1$-conjugation coincide.  Consider now
any $b \in \tilde{B}_n$ and any $g \in \underline{G}(n)$.  Let $h
= g_{b^{-1}}$.  We have:
$$g_{(\tX^2_{1})_{b}} = g_{b^{-1}\tX^2_1b} = ((h)_{\tX^2_1})_b =
(h_{s_1})_b = (h_{b})_{(s_{1})_{b}}=g_{(s_{1})_{b}}.$$

(5), (6), (7) are easy to verify. \quad $\square$

\edm

\proclaim{Lemma III.6.3}
Let $n \geq 3$.

Let $\{X_1,\dots,X_{n-1}\}$ be a frame of $B_n.$

Let
$$Z_{ij}=\cases X_1 & \text{if}\ (i,j)=(1,2);\\
(X_{1})_{X_2\dots X_{j-1}}& \text{if}\ i=1,\ j\geq 3;\\
(X_{1})_{X_2\dots X_{j-1}X_1\dots X_{i-1}}\quad&\text{if}\ i\geq
2,\ j>i.\endcases$$

Let $\tZ_{ij}$ be the image of $Z_{ij}$ in $\tilde B_n.$

Consider $\underline G(n)$ as a $B_n$-group as in Claim \rom{III.6.2}.

Then there exists a unique $\tilde{B}_n$-surjection $\Lambda_n: \tilde{P}_n
\rightarrow \underline{G}(n)$ with $\Lambda_n(\tX^2_1) = s_1$
and $\Lambda_n(\tZ^2_{ij}) = s_{ij}$ for $1 \leq i < j \leq
n$.\ep

\demo{Proof}
Use induction on $n$.

For $n = 3,  $ $\Lambda_3\: \tilde{P}_3 \rightarrow
\underline{G}(3)$ must be defined by $\lambda_3(\tZ^2_{ij}) =
s_{ij}$, $1 \leq i < j \leq 3$.  One can check directly that
$\Lambda_3$ is well
defined, and that it is a $\tilde{B}_3$-surjection.  Uniqueness of
such $\Lambda_3$ is evident.

Assume now that $n \geq 4$ and that the desired $\Lambda_{n-1}:
\tilde{P}_{n-1} \rightarrow \underline{G}(n-1)$ exists.

Considering $(X_1,X_2) \subset (X_1, X_2, X_3) \subset \ldots
\subset (X_1, \ldots, X_{n-1})$, we get a chain of embeddings
$B_3 \subset B_4 \subset \ldots \subset B_n$ and the
corresponding chain $P_3 \subset P_4 \subset \ldots \subset P_n$.
To the latter corresponds a chain of homomorphisms:  $\tilde{P}_3
\overset{i_3}\to{\longrightarrow} \tilde{P}_4
\overset{i_4}\to{\longrightarrow} \ldots \rightarrow \tilde{P}_{n-1}
\overset{i_{n-1}}\to{\longrightarrow} \tilde{P}_n$, where $\tilde{P}_3$
 is obtained (by definition) from $P_3$ by adding the relations:
$[\tZ^2_{12},\tZ^2_{23}] = [\tZ^2_{12},\tZ_{13}^2] = [\tZ^2_{23},
\tZ^2_{13}]$
and $[\tZ^2_{12},\tZ^2_{23}]^2 = 1$.

It is known that the set $\{Z^2_{ij}, 1 \leq i < j \leq n\}$
generates $P_n$, and $P_n \simeq P_{n-1} \ltimes
F_{n-1},$ where $P_{n-1}$ is the subgroup of $P_n$ generated by
$\{Z^2_{ij}, 1 \leq i < j \leq n - 1\}$, $F_{n-1}$ is a free
subgroup of $P_n$ generated by $\{Z^2_{in}, 1 \leq i \leq n - 1
\}$, and the semi-direct product $P_{n-1} \ltimes
F_{n-1}$ is defined according to the $P_{n-1}$-action on
$F_{n-1}$ which comes from the $B_{n-1}$-action by conjugation
(using $B_{n-1} \subset B_n \supset P_n$).  The latter coincides
with the standard $B_{n-1}$-action on $F_{n-1}$ (the generators
$X_{n-2}, \ldots, X_1$ of $B_{n-1}$
correspond to standard {\it Hurwitz moves} on $(Z^2_{n-1,n},
Z^2_{n-2,n}, \ldots , Z^2_{1n})$ (see [MoTe4], Chapter 4).

  Using
canonical $P_{n-1} \rightarrow \tilde{P}_{n-1}$, we obtain from
$\Lambda_{n-1}$ a $B_{n-1}$-surjection $\hat{\Lambda}_{n-1}:
P_{n-1} \rightarrow \underline{G}(n)$.  For the free subgroup
$F_{n-1}$ of $P_n$ generated by $\{Z^2_{in}, \; n - 1 \leq i \leq
1\}$ define $\mu_{n-1}: F_{n-1} \rightarrow \underline{G}(n)$ by
$\mu_{n-1}(Z^2_{in}) = s_{in}$.  Considering $P_n$ as $P_{n-1}
\ltimes F_{n-1}$, we  define
$\hat{\Lambda}_n: P_n
\rightarrow \underline{G}(n)$ which on $P_{n-1}$
coincides with
$\hat{\Lambda}_{n-1}: P_{n-1} \rightarrow \underline{G}(n-1)
\subset \underline{G}(n)$ (see Claim III.6.2(7)) and on $F_{n-1}$ coincides
with
$\mu_{n-1}: F_{n-1} \rightarrow \underline{G}(n)$.  To show that such
$\hat{\Lambda}_n$ exists one has to check the following:

1)  The conjugation of $\mu_{n-1}(F_{n-1})$ by elements of
$\hat{\Lambda}_{n-1}(P_{n-1})(\subset \underline{G}(n))$
coincides with the $P_{n-1}$-action defined by $P_{n-1} \subset
B_{n-1} \subset   B_n \rightarrow \tilde{B}_n$ and the given
$\tilde{B}_n$-action on $\underline{G}(n)$.  That
is, $\forall f \in \mu_{n-1}(F_{n-1})$ and $\forall h$ of the form $
\hat{\Lambda}_{n-1}(\tY)$ ($Y \in P_{n-1}$) we must have $h^{-1}fh
= f_{\tY}$.

2)  The $P_{n-1}$-action on $\mu_{n-1}(F_{n-1})$ (defined by
$P_{n-1} \subset B_{n-1} \subset B_n \rightarrow \tilde{B}_n$ and
the given $\tilde{B}_n$-action on $\underline{G}(n)$) comes from
$B_{n-1}$-action on $\mu_{n-1}(F_{n-1})$ in which $X_{n-2},
\ldots X_1$ correspond to the standard {\it Hurwitz moves} on
$(s_{n-1,n},s_{n-2,n}, \ldots s_{1n})$.\edm

\demo{Proof of 1)}
Since $\forall b \in B_{n-1},$\   $\hat{\Lambda}_{n-1}((X^2_{1})_{b})
= (s_{1})_{b}$ we see from Claim III.6.2 that $\forall f \in \mu_{n-1}(F_{n-
1}),$\quad   $f_{\hat{\Lambda}_{n-1}((X^2_{1})_{b})} = (f_{s_1})_{b} =
f_{(\tX^2_{1})_{b}}$.  Since $P_{n-1}$ is generated by $\{Z^2_{ij}, 1
\leq i < j \leq n - 1\}$, i.e., by $\{(X^2_{1})_{b}, b \in B_n\}$
we get 1).\edm

\demo{Proof of 2)} It follows immediately from Claim III.6.2(6).

Thus, 1) and 2) are true and we can extend $\hat{\Lambda}_{n-1}$,
$\mu_{n-1}$ to a homomorphism $\hat{\Lambda}_n: P_n(= P_{n-1}
\ltimes F_{n-1}) \rightarrow \underline{G}(n)$ such
that for $1 \leq i < j \leq n - 1$ $\hat{\Lambda}_n(Z^2_{ij}) =
\hat{\Lambda}_{n-1}(Z^2_{ij}) = s_{ij}$, and for $1 \leq i \leq n
- 1$ $\hat{\Lambda}_n(Z^2_{in}) = \mu_{n-1}(Z^2_{in}) = s_{in}$,
in short $\hat{\Lambda}_n(Z^2_{ij}) = s_{ij}$ for $1 \leq i < j
\leq n$.

Using induction, one can check directly that $\hat{\Lambda}_n$ is
a $B_n$-homomorphism (recall that by  Claim III.6.2 we have explicit formulas
for $s_{ij}$'s).

Because $s_{1n} = u_{n-1} \cdot u_2s_1$ (by 1.20), that is,
$u_{n-1} = s_{1n}(u_{n-2} \ldots u_2s)^{-1}$, we see that
$\underline{G}(n)$ is generated by $\underline{G}(n-1) =
\hat{\Lambda}_n(P_{n-1})$ ($= \hat{\Lambda}_{n-1}(P_{n-1}))$ and
$s_{1n} = \hat{\Lambda}_n(Z^2_{1n})$.  Therefore
$\hat{\Lambda}_n$ is a $B_n$-surjection.

Let ${N} = \ker (B_n \rightarrow \tilde{B}_n)\ (=\ker(P_n
\rightarrow \tilde{P}_n)$). Let $T = X^2_1X^2_3X^{- 2}_2Z^{-2}_{14}$.  Clearly,
$ {N}$ is generated by $\{T_b,b \in B_n\}$.  We have
$\hat{\Lambda}_n(T) = \hat{\Lambda}_4(T) = s_1 \cdot s_{34}
\cdot s^{-1}_{23} s^{-1}_{14} \overset\text{Claim III.6.2}\to= s_1 \cdot \eta
u_3u_2 u_1 \cdot u_2s_1 \cdot s^{-1}_1u^{-1}_1 u^{-1}_2 \eta
\cdot s^{-1}_1u^{-1}_2 u^{-1}_3 = s_1 \eta u_3u_2 \cdot s^{-1}_1
u^{-1}_2 u^{-1}_3 =$ Id.  Since $\hat{\Lambda}_n$ is a $B_n$-homomorphism, we
get $\hat{\Lambda}_n(T_b) =$ Id $\forall b \in B_n$, and thus $\hat{\Lambda}_n(
{N}) = $ Id.  Hence $\hat{\Lambda}_n$ defines canonically a
$\tilde{B}_n$-surjection $\Lambda_n: \tilde{P}_n \rightarrow \underline{G}(n)$
with $\Lambda_n(\tX^2_1) = s_1$.

Uniqueness of such $\Lambda_n$ follows from the fact that
$\tilde{P}_n$ is generated by the $B_n$-orbit of $\tX^2_1$. \qed\edm
 \proclaim{Theorem III.6.4}
There exists a unique $\tilde{B}_n$-isomorphism $\Lambda_n:
\tilde{P}_n \rightarrow \underline{G}(n)$  with
$\Lambda_n(\tX^2_1) = s_1$.  In
particular:
\roster\item $\Ab\tilde{P}_n$ is a free abelian group with $n$ generators,
\newline $\tilde{P}^\prime_n \simeq {\Bbb Z}/2$ (generated by
$c= [\tX^2_1,\tX^2_2]$);
\item $\tilde{P}_{n,0}$ is $\tilde{B}_n$-isomorphic to the subgroup
$G_0(n)$ of $\underline{G}(n)$, generated by\newline $u_1, \ldots, u_{n-
1}$, $\Ab\tilde{P}_{n,0}$ is a free abelian group with $n-1$
generators\newline  $\{\xi_1,  \ldots, \xi_{n-1}\}$,
$\tilde{P}^\prime_{n,0} \simeq {\Bbb Z}/2,$ generated by
$c=[\xi_1,\xi_2];\ (\nu\in G_0(n)$ corresponding
 to $c\in
\tilde P_{n,0}).$
\item  $\tilde{P}_{n,0}$ is a primitive $\tilde{B}_n$-group
generated by the $\tilde{B}_n$-orbit of a prime element
${u} = \tX^2\tY^{-2}$, where $\tX,\tY$ are adjacent half-twists
in $\tilde{B}_n$, $\tilde T = \tX\tY \tX^{-1}$ is a supporting half-twist for
${u}$.\endroster
\ep

\demo{Proof}
Clearly, $\tilde{P}_{n,0}$ is generated by $\{\tX^2_1\tZ^{-2}_{ij}, l
\leq i < j \leq n\}$.  Because $\tX^2_1 \cdot \tZ^{-2}_{ij} = \tX^2_1
\tZ^{-2}_{1i} \cdot \tZ^2_{1i} \cdot \tZ^{-2}_{ij}$ and both $\tX^2_1\tZ^{-
2}_{1i}$, $\tZ^2_{1i} \cdot \tZ^{-2}_{ij}$ are conjugates of $u$, we
see that $\tilde{P}_{n,0}$ is generated by the $\tilde{B}_n$-
orbit of $u$, and since $u$ is a prime element of $P_{n,0},$ this means that
$\tilde{P}_{n,0}$ is a primitive $\tilde{B}_n$-group.
By Lemma III.6.1, the s.h.t. of $u$ is $Y\1XY.$
Thus, we proved (3).

Polarize each $X_i$ (and $\tX_i$) according to the sequence $(X_1,
\ldots X_{n-1}$) (the ``end'' of $X_i$ = the ``origin'' of
$X_{i+1}$).  By Proposition III.3.3 $\forall i = 1, \ldots, n - 1$ $\exists$
unique prime element $\xi_i = L_{\{u,\tX_1\}}(\tX_i) \in
\tilde{P}_{n,0}$ such that $\{\xi_i,\tX_i\}$ is coherent
with
$\{u,\tX_1\}$.  Clearly $\xi_1 = u$.

By Proposition III.4.1 we have $\forall i = 1, \ldots, n - 1$:
$$
(\xi_{i})_{\tX^{-1}_i} = \xi^{-1}_ic; \;
(\xi_{i})_{\tX^{-1}_{i-1}} =
\xi_i \xi_{i-1}; \;
=(\xi_{i})_{\tX^{-1}_{i+1}} = \xi_i\xi_{i+1}.
\tag1.22$$
It is clear also that $\forall j \neq i$, $i - 1, i + 1$
$$
(\xi_{i})_{\tX_j} = \xi_i.
\tag1.23$$
We see from (1.22), (1.23) that the subgroup of $\tilde{P}_{n,0}$
generated by $(\xi_1, \ldots, \xi_{n-1})$ is
closed under the
$\tilde{B}_n$-action.  Since $\tilde{P}_{n,0}$ is generated by
the $\tilde{B}_n$-orbit of $u = \xi_1$, we conclude
that
$\tilde{P}_{n,0}$ is generated by $(\xi_1, \ldots,
\xi_{n-1})$.  This implies that $\tilde{P}_n$ is
generated by $(\tX^2_1,\xi_1,\xi_2, \ldots,
\xi_{n-1}$).

We have $\xi_2 = L_{\{u,\tX_1\}}(\tX_2) =
(\xi_{1})_{\tX^{-1}_2\tX^{-1}_1} = (\tX^{-2}_{1})_{\tX_2} \cdot \tX^2_1$,
which implies that
$$
[\tX^2_1,\xi_2] = c.
\tag1.24$$

By Lemma III.5.1 we have
$$
[\xi_i,\xi_j] =
\cases
c\quad & \text{if}\ |i-j| = 1 \\
1\quad &\text{if} \ |i-j| \geq 2.
\endcases\tag1.25
$$

Observe also that
$$
(\tX^2_{1})_{\tX_2} = \xi_2 \cdot \tX^2_1.
\tag1.26$$

Formulas (1.22)--(1.26) show that we can define a $\tilde{B}_n$-homomorphism
$M_n: \underline{G}(n) \rightarrow \tilde{P}_n$ with $M_n(s_1) = \tX^2_1$,
$M_n(u_i) = \xi_i$, $i = 1,\ldots, n - 1$.  (See Claim III.6.2.)

Since  $\tilde{P}_n$ is generated by the
$\tilde{B}_n$-orbit of $\tX^2_1$ and $\underline{G}(n)$ is generated by
the $\tilde{B}_n$-orbit of $s_1$, we conclude that $\Lambda_n$
and $M_n$ are inverses of each other.\quad $\square$\edm

\newpage

\subheading{\S7. Criterion for prime element}

\proclaim{Proposition III.7.1}
Assume $n \geq 5$.  Let $G$ be a $\tilde{B}_n$-group,
$(\tX_1,\tX_2, \ldots, \tX_{n-1})$ be a standard base of $\tilde{B}_n$.
Let $S$ be an element of $G$ with the following properties:
\roster\item"(0)"  $G$ is generated by $\{S_b,b \in \tilde{B}_n\}$;
\item"$(1_a)$"    $S_{\tX^{-1}_2\tX^{-1}_1} = S^{-1}S_{\tX^{-1}_2}$;
\item"$(1_b)$" $S_{\tX_1\tX^{-1}_2\tX^{-1}_1}$ = $S^{-
1}_{\tX_1}S_{\tX_1\tX^{-1}_2};$
\item"(2)"  For $\tau = SS_{\tX^{-1}_1},$ $T = S_{\tX^{-1}_2}$ we have:
\item"$(2_a)$"  $\tau_{\tX^2_1} = \tau$;
\item"$(2_b)$" $\tau_T = \tau^{-1}_{\tX_1}$;
\item"(3)" $S_{\tX_j} = S$  $\forall j \geq 3$;
\item"(4)"  $S_c = S$, where $c = [\tX^2_1,\tX^2_2]$.\endroster

Then $S$ is a prime element of $G$, $\tX_1$ is a supporting half-twist of
$S$ and $\tau$ is the corresponding central element.  In particular, $\tau^2 =
1$, $\tau \in\Center(G)$, $\tau_b = \tau$ $\forall$ $b \in \tilde{B}_n$. \ep

\pf
The proof includes  several lemmas. From Theorem III.5.2,
$c\in\Center(\tB_n),$\  $c^2 = 1.$  From Theorem III.6.4 it follows that
$\tilde{P}^\prime_n$ is generated by $c$.  $\tilde{P}^\prime_n$ is a normal
subgroup of $\tilde{B}_n$.  Denote by $\tilde{\tilde{B}}_n =
\tilde{B}_n/\tilde{P}^\prime_n$, $\tilde{\tilde{P}}_n =
\tilde{P}_n/\tilde{P}^\prime_n = \Ab \tilde{P}_n$. $\Tilde{\Tilde P}_n$ is a
commutative group.  We have $\Tilde{\Tilde\psi}_n:\Tilde{\Tilde
B}_n\twoheadrightarrow S_n.$
By abuse of notaion we use $\psi_n$ for $\Tilde{\Tilde\psi}_n.$ Let $Y\in B_n.$
By abuse of notation we denote the image of $Y$ in $\tB_n$ or in $\Tilde{\Tilde
B}_n$ or in $\Tilde{\Tilde P}_n$  by the same symbol $\tY.$  It is clear that
$\tilde{B}_n$ acts on $\Tilde{\Tilde{P}}_n$ (through conjugations) as the
symmetric group $S_n = \tilde{B}_n/\tilde{P}_n$.

Since $S_c = S$ and $c\in\Center\tilde B_n,$ we have $\forall b \in
\tilde{B}_n$
$(S_b)_c = S_{bc} = (S_{c})_{b} = S_b$. Since $G$ is generated by $\{S_b,b
\in \tilde{B}_n\}$ we have $\forall g
\in G$ \ $g_c = g.$   In particular, we conclude that $\tilde{B}_n$ acts on $G$
as its quotient $\Tilde{\Tilde{B}}_n$; in other words, $G$ is a
$\Tilde{\Tilde{B}}_n$-group.

  Let $(D,K)$ be a model for $B_n,$\
${K} = \{a_1, \ldots, a_n\}$,\ $B_n=B_n[D,K].$ Take any $a_{i_1},a_{i_2} \in
{K}$.  Let $\gamma_1,\gamma_2$ be two different simple paths in
$D - ({K} - a_{i_1}-a_{i_2})$ connecting $a_{i_1}$ with $a_{i_2}$, let
${H}(\gamma_1),{H}(\gamma_2)$ be the half-twists
corresponding to
$\gamma_1,\gamma_2$, and let ${\tilde H}(\gamma_1),
\tilde{H}(\gamma_2)$ be the images of ${
H}(\gamma_1), {H}(\gamma_2)$ in $\tilde{\tilde{B}}_n$. \qed\enddemo

\proclaim{Lemma 1} Let $\g_1,\g_2$ be $2$ simple paths in
$D-\{K-a_{i_{1}}-a_{i_{2}}\}$ connecting $a_{i_{1}}$ with $a_{i_{2}}.$
Then:
$\tilde{H}(\gamma_1)^2 = \tilde H(\gamma_2)^2$. \ep
\demo{Proof of Lemma 1}
Choose a frame of ${B}_n$ $(Y_1, \ldots,
Y_{n-1})$ s.t. $Y_1 =  {H}(\gamma_1).$ Let $b \in
B_n$ s.t. $\gamma_2 = (\gamma_1)b$, that is ${H}(\gamma_1)_b
= {H}(\gamma_2)$.
Let $\tY_i$ be the image of $Y_i$ in $\Tilde{\Tilde B}_n.$

Let $\sigma_1$ be the image of $b$ in $S_n$.  Since $(a_i)b =
a_i$, $(a_j)b = a_j$, $\sigma_1 \in
\operatorname{Stab}(i)\cap \operatorname{Stab}(j)$ in $S_n$. The subgroup of
$\Tilde{\Tilde{B}}_n$ generated by $\tY_3, \ldots, \tY_{n-1}$ is mapped by
$\Tilde{\Tilde \psi}_n:\Tilde{\Tilde B}_n\ri S_n$ onto
$\operatorname{Stab}(i)\cap \operatorname{Stab}(j).$  Choose $\tilde b_1$ in
this subgroup with its image in $S_n$ equal to $\sigma_1$.  Clearly,
$(\tY_{1})_{\tilde b_1} = \tY_1$.   Since the image of
$\tilde b^{-1}_1\tilde{b}$ in $S_n$ is equal to $\sigma^{-1}_1\sigma_1 = $ Id,
we have $\tilde b^{-1}_1 \tilde{b} \in \Tilde{\Tilde{P}}_n.$ Since
$\Tilde{\Tilde P}_n$ is commutative when considering $\tY_1^2$ as an element
of $\Tilde{\Tilde P}_n,$\  $(\tY^2_{1})_{b^{-1}_1\tilde{b}} = \tY^2_1$.  Thus,
we have
$$\tilde{H}(\gamma_2)^2 = \tilde H(\gamma_1)^2_{\tilde{b}} =
(\tY^2_{1})_{\tilde b_1\tilde b^{-1}_1\tilde{b}} =
(\tY^2_{1})_{\tilde b^{-1}_1\tilde{b}} = \tY^2_1 = \tilde H(\gamma_1)^2
\quad\square \quad\text{ for Lemma 1}$$\edm

\demo{Definition}\ $\underline{f_{ij}}$

$\forall i,j \in (1, \ldots, n)$, $i \neq j$, we define $f_{ij}
\in \Tilde{\Tilde{P}}_n$ as follows:  Take any simple path
$\gamma$ in $D - ({K}-a_i-a_j)$ connecting $a_i$ with
$a_j$.  Let $f_{ij} = \tilde{H}(\gamma)^2$.  Lemma 1 shows that this
definition does not depend on the choice of $\gamma$. We choose for $i<j:$
$$f_{ij}=\cases (\tX_i^2)_{X_{i+1}}\cdot\dots\cdot X_{j-1}\quad & 2<j-i\\
\tX_i^2 & i+1=j.\endcases$$

It
is clear that for  $\sigma_1$  the image   in $S_n$ of  $ b \in
\tilde{B}_n$ we have:
$$ \quad(f_{ij})_b = f_{(i)\sigma_1,(j)\sigma_1}.$$

It is clear from our choice of $\g$ for $f_{ij}= \tilde
H(\gamma)^2$  that: $$
\psi_n(\tilde{H}(\gamma)) = (i,j).$$
\medskip

It will be convenient to use the following notation for  $ g
\in G$ and $b \in \Tilde{\Tilde{B}}_n:$
\demo{Notation} \ $\un{[g,b]}:$

For $g\in G$, $ b\in\Tilde{\Tilde B}_n$ and the action of $\Tilde{\Tilde B}_n$
on $G$, we denote $[g,b] = g \cdot g^{-1}_{b^{-1}}.$
\medskip
One can check that:$$\align&g_b=g\Leftrightarrow [g,b]=1\\
&[g,b]_z=g_z(g_z)_{b_z\1}\1\\
&[g^{-1},b] = [g,b]^{-1}_g\\  &[g,b^{-1}] =
[g,b]^{-1}_b\\&[g_1g_2,b] = [g_2,b]_{g^{-1}_1} \cdot [g_1,b]\\
&[g,b_1b_2] = [g,b_1] \cdot [g,b_2]_{b^{-1}_1}.\endalign$$

\demo{Notation}\ $\underline{Q_{b,l,m}}$

$\forall b \in \Tilde{\Tilde{B}}_n,\quad \forall l,m \in (1, \ldots, n),$ $l
\neq m$, we denote $Q_{b,l,m} = [S_b,f^{-1}_{lm}].$
\edm

\proclaim{Lemma 2}

\roster\item"(i)"Let $b \in \Tilde{\Tilde{B}}_n$ be such that
$(\{1,2\})\psi_n(b) \cap \{l,m\} = \emptyset$.  Then $Q_{b,l,m} = 1$.
\item"(ii)" Let $Q = Q_{\operatorname{Id},1,3} = [S,f^{-1}_{13}]$.
Then
$Q_{\tX^{-1}_2} = Q$.\endroster
\ep

\demo{Proof of Lemma 2}

(i)\
Let $\{l_1,m_1\} = (\{l,m\})\psi_n(b)^{-1}$.  So
$(f_{lm})_{b\1}=f_{l_{1}m_{1}}.$ We have $\{1,2\} \cap \{l_1,m_1\} =
\emptyset$,
that is, $3 \leq l_1   , 3\le m_1$.  By our choice, $f_{lm}$ is a product o
$X_j$ for $j\ge 3.$
Thus, using property (3) of $S$  we get
$S_{f_{l_1,m_1}} = S.$  In other words,  $[S,f_{l_1,m_1}\1] = 1$.  We get
$$(Q_{b,l,m})_{b^{-1}} = [S_b,f^{-1}_{lm}]_{b^{-1}} =
[S,(f_{lm}\1)_{b^{-1}}] = [S,f_{l_1,m_1}\1] = 1,$$
and so $Q_{b,l,m} = 1$.   \quad\qed \quad for Lemma 2(i)

(ii)\
{}From $S_{\tX^{-1}_2\tX^{-1}_1} = S^{-1}S_{\tX^{-1}_2}$ (assumption $(1_a$)
 of the
Proposition) it follows that $S_{\tX^{-1}_2} = S
S_{\tX^{-1}_2\tX^{-1}_1}$.  Applying $\tX^{-1}_3$, we get
$$
S_{\tX^{-1}_2\tX^{-1}_3} = SS_{\tX^{-1}_2\tX^{-1}_1\tX^{-1}_3}
\tag1.27$$
which, after applying $\tX^{-1}_2$, gives:
$$S_{\tX^{-1}_2\tX^{-1}_3\tX^{-1}_2} = S_{\tX^{-1}_2} S_{\tX^{-1}_2\tX^{-
1}_1\tX^{-1}_3\tX^{-1}_2}.$$
Since $S_{\tX^{-1}_2\tX^{-1}_3\tX^{-1}_2} = S_{\tX^{-1}_3\tX^{-1}_2\tX^{-
1}_3} = S_{\tX^{-1}_2\tX^{-1}_3}$, we obtain $S_{\tX^{-1}_2\tX^{-1}_3}
=$\newline $ S_{\tX^{-1}_2}S_{\tX^{-1}_2\tX^{-1}_1\tX^{-1}_3\tX^{-1}_2}$, or
$$
S_{\tX^{-1}_2} = S_{\tX^{-1}_2\tX^{-1}_3} S^{-1}_{\tX^{-1}_2\tX^{-1}_1\tX^{-
1}_3\tX^{-1}_2}
\tag 1.28$$

Let  $b_1=\tX^{-
1}_2\tX^{-1}_1\tX^{-1}_3\tX^{-1}_2$. Observing that $(f_{13})_{\tX^{-1}_2} =
f_{12}$, we get\linebreak from (1.28): $Q_{\tX^{-1}_2} =
[S_{\tX^{-1}_2},(f^{-1}_{13})_{\tX^{-1}_2}] = [S_{\tX^{-
1}_2\tX^{-1}_3}S^{-1}_{b_1},f^{-1}_{12}].$  Thus:
$$Q_{\tX^{-1}_2}
 = [S^{-1}_{b_1},f^{-
1}_{12}]_{S^{-1}_{\tX^{-1}_2\tX^{-1}_3}} \cdot [S_{\tX^{-1}_2\tX^{-1}_3},
f^{-1}_{12}].
\tag1.29$$

Since $\psi_n(b_1) = (2\quad 3)\ (1\quad 2)\ (3\quad 4)\ (2\quad 3)$ (products
of transpositions),\linebreak $(\{1,2\})\psi_n(b_n) = \{3,4\}.$ Since $\{3,4\}
\cap \{1,2\} = \emptyset$, we get from   (i)
 that $Q_{b_1,1,2} = 1.$ Thus, $[S^{-1}_{b_1},f^{-1}_{12}] =
[S_{b_1},f^{-1}_{12}]^{- 1}_{S_{b_1}} = (Q^{-1}_{b_1,1,2})_{S_{b_1}}$.
(1.29) now gives:
$$
Q_{\tX^{-1}_2} = [S_{\tX^{-1}_2\tX^{-1}_3},f^{-1}_{12}].
\tag1.30$$

Consider a quadrangle formed by $\{a_1,a_2,a_3,a_5\}$, as in
Fig. III.7.1.  By Lemma III.1.2, we can write in $\Tilde{\Tilde{P}}_n$:
$f_{35}f_{12} = f_{25}f_{13}$, or $f_{12} = f_{35}^{-
1}f_{25}f_{13}$, $f^{-1}_{12} = f^{-1}_{13} f^{-
1}_{25}f_{35}$.

\midspace{1.40in}
\caption{Fig. III.7.1}

{}From (1.30) we get:
$$
Q_{\tX^{-1}_2} =
[S_{\tX^{-1}_2\tX^{-1}_3},f^{-1}_{13}f^{-1}_{25}f_{35}] =
[S_{\tX^{-1}_2\tX^{-1}_3},f^{-1}_{13}] \cdot [S_{\tX^{-1}_2\tX^{-
1}_3},f^{-1}_{25}f_{35}]_{f_{13}}
\tag1.31$$

Consider $[S_{\tX^{-1}_2\tX^{-1}_3},f^{-1}_{25}f_{35}] = [S_{\tX^{-
1}_2\tX^{-1}_3}, f^{-1}_{25}] \cdot [S_{\tX^{-1}_2\tX^{-
1}_3},f_{35}]_{f_{25}} = Q_{\tX^{-1}_2\tX^{-1}_3,2,5} \cdot [S_{\tX^{-
1}_2\tX^{-1}_3},f^{-1}_{35}]^{-1}_{f^{-1}_{35}f_{25}} = Q_{\tX^{-
1}_2\tX^{-1}_3,2,5} \cdot (Q^{-1}_{\tX^{-1}_2\tX^{-1}_3,3,5})_{f^{-
1}_{35}f_{25}}.$  Since $\psi(\tX^{-1}_2\tX^{-1}_3) =
(2\quad 3)(3\quad 4),$ the images of $\{1,2\}$ under it are  $
\{1,4\}.$ But $\{1,4\} \cap \{2,5\} = \emptyset$ and $\{1,4\} \cap \{3,5\} =
\emptyset.$ Thus,  we get by (i) that
$Q_{\tX^{-1}_2\tX^{-1}_3,2,5} = Q_{\tX^{-1}_2\tX^{-1}_3,3,5} = 1$, and
so $[S_{\tX^{-1}_2\tX^{-1}_3},f^{-1}_{25}f_{35}] = 1$.  (1.31)
now implies $Q_{\tX^{-1}_2} = [S_{\tX^{-1}_2\tX^{-1}_3}, f^{-1}_{13}]$.
By (1.27) $S_{\tX_2^{-1}\tX^{-1}_3} = S \cdot
S_{\tX^{-1}_2\tX^{-1}_1\tX^{-1}_3}$ which gives
$$\align Q_{\tX^{-1}_2} &= [S \cdot S_{\tX^{-1}_2\tX^{-1}_1\tX^{-1}_3}, f^{-
1}_{13}] = [S_{\tX^{-1}_2\tX^{-1}_1\tX^{-1}_3},f^{-1}_{13}]_{S^{-1}}
\cdot [S,f^{-1}_{13}] \\&=
(Q_{\tX^{-1}_2\tX^{-2}_1\tX^{-1}_3,1,3})_{S^{-1}} \cdot Q.\endalign$$

The value of $ \psi(\tX^{-1}_2\tX^{-1}_1\tX^{-1}_3)\ ( =
 (2 \quad 3)(1\quad 2)(3\quad 4))$ on $\{1,2\}$ is $\{2,4\}$.  Since $\{2,4\}
\cap \{1,3\} = \emptyset$ we get from part (i) of the Lemma that
$Q_{\tX^{-1}_2\tX^{-1}_1\tX^{-1}_3,1,3} = 1$, therefore,
$$Q_{\tX^{-1}_2} = Q \qed\quad \text{for   Lemma 2(ii).}$$
\edm
\proclaim{Lemma 3} $\tau=Q\1.$\ep
\demo{Proof}
By the assumption on $\tau,$\ $S_{\tX^{-1}_1} = S^{-1}
\tau.$  By definition of $T,$\  $T_{\tX^{-2}_1} = S_{\tX^{-1}_2\tX^{-2}_1}.$
We apply assumption $(1_a)$ twice to get, using $S_{\tX_{1}\1}=S\1\tau,$ that
\newline  $(S^{- 1}S_{\tX^{-1}_2})_{\tX^{-1}_1} =
S^{-1}_{\tX^{-1}_1}S_{\tX^{-1}_2\tX^{- 1}_1} = \tau^{-1}S \cdot
S^{-1}S_{\tX^{-1}_2} = \tau^{-1} S_{\tX^{- 1}_2} = \tau^{-1}T$.  Thus
$$
T_{\tX^{-2}_1} = \tau^{-1}T,
$$
or
$$\tau^{-1} = T_{\tX^{-2}_1} T^{-1}.\tag1.32$$
Applying $\tX^2_1$ on (1.32) and using $\tau_{\tX^2_1} = \tau$
(assumption $(2_a)$), we get $$\align\tau^{-1}& = T \cdot
T^{-1}_{\tX^2_1} = [T,\tX^{- 2}_1] = [S_{\tX^{-1}_2},f^{-1}_{12}] =
[S,(f^{-1}_{12})_{{\tX_2}}]_{\tX^{- 1}_2} = [S,f^{-1}_{13}]_{\tX^{-1}_2}\\& =
Q_{\tX^{-1}_2} \overset\text{by
 Lemma 2}\to =
Q,\endalign$$ that is, $$
\tau^{-1}=Q, \quad \text{or} \quad \tau = Q^{-1}.\quad\square\quad\text{for
Lemma 3} \tag 1.33$$

\proclaim{Lemma 4}
$\forall j \geq 3$, $\tau_{\tX_j} = \tau.$
\ep

\demo{Proof of Lemma 4}
{}From $\tau = SS_{\tX^{-1}_1}$ and $S_{\tX_j} = S$ $\forall j \geq 3$
it follows that $\tau_{\tX_j} = \tau$ $\forall j \geq 3$.    \qed\quad for
Lemma 4.\edm

\proclaim{Lemma 5}
$\tau_{\tX_1} = \tau.$
\ep

\demo{Proof of Lemma 5}  Let us use
now $\tau^{-1}_{\tX_1} = \tau_T$ ($(2_b)$ of the Proposition).

By Lemmas 2 and 3  $\tau_{\tX_2} = \tau$.

 Thus,
$\tau_T = \tau_{S_{\tX^{-1}_2}} = S^{-1}_{\tX^{-1}_2} \tau
S_{\tX^{-1}_2} = (S^{-1}\tau S)_{\tX^{-1}_2} = (\tau_S)_{\tX^{-
1}_2}$.

 So $\tau_{\tX^{-1}_1} = \tau_{\tX_1} =\tau_T\1= (\tau^{-
1}_S)_{\tX^{-1}_2}$, or
$$
\tau_S = \tau^{-1}_{\tX^{-1}_1\tX_2}.
\tag1.34$$

Since $\tau_{\tX_3} = \tau$ and $S_{\tX_3} = S$, we get
$(\tau_S)_{\tX_3} = \tau_S$ and
$$
\tau_{\tX^{-1}_1\tX_2\tX_3} = (\tau^{-1}_S)_{\tX_3} = \tau^{-1}_S =
\tau_{\tX^{-1}_1\tX_2}.
\tag1.35$$

Applying $\tX^{-1}_2$ on (1.35), we get
$
\tau_{\tX^{-1}_1\tX_2\tX_3\tX_2^{-1}} = \tau_{\tX^{-1}_1}.
 $
Since $\tau_{X_3}=\tau$ and $\la X_2,X_3\ra=1,$\
$\tau_{\tX^{-1}_1\tX_2\tX_3\tX_2^{-1}} = \tau_{\tX^{-1}_1\tX^{- 1}_3\tX_2\tX_3}
= \tau_{\tX^{-1}_3\tX^{-1}_1\tX_2\tX_3} = \tau_{\tX^{-1}_1\tX_2\tX_3}$.
Thus:
$$\tau_{\tX_{1}\1\tX_{2}\tX_{3}}=\tau_{\tX_{1}\1}\tag1.36.$$
Combining
formulas (1.35)-(1.36) we get $\tau_{\tX^{-1}_1} =
\tau_{\tX^{-1}_1\tX_2}$.  Applying it to $\tX_1$ we get $\tau=
\tau_{\tX^{-1}_1\tX_2\tX_1 } = \tau_{\tX_2\tX_1\tX^{-1}_2}  =
\tau_{\tX_1\tX^{-1}_2}$. Thus $\tau = \tau_{\tX_1\tX^{-1}_2}$,
or $\tau_{\tX_1} = \tau_{\tX_2} = \tau$. \linebreak
\qed\quad for  Lemma 5.\edm

\proclaim{Lemma 6}
$\tau_{\tX_j} = \tau$ $\forall j = 1,2,\ldots, n - 1$.
\ep
\demo{Proof of Lemma 6} By Lemmas 2, 3, 4, 5.\quad \qed\edm
\proclaim{Lemma 7} $\tau_S = \tau^{-1}$.
\ep
\demo{Proof of Lemma 7}
{}From $\tau_{\tX_1} = \tau_{\tX_2} = \tau$ and (1.34).\quad \qed
\edm

\proclaim{Lemma 8}
$\tau_S = \tau$.
\ep

\demo{Proof of Lemma 8}
Consider assumption $(1_b)$ of the Proposition:
$$S_{\tX_1\tX^{-1}_2\tX^{-1}_1} = S^{-1}_{\tX_1}S_{\tX_1\tX^{-1}_2}.$$

Using $S_{\tX^{-1}_1} = S^{-1} \tau$ and $\tau_{\tX_1} = \tau$, we
get $S = S^{-1}_{\tX_1}\tau$, or $S^{-1}_{\tX_1} = S\tau^{-1}$,
$S_{\tX_1} = \tau S^{-1}$.  Assumption $(1_b)$ now gives (using $\tau_{\tX_i}
= \tau \; \forall i = 1,\ldots, n - 1)$:
$$\tau S^{-1}_{\tX^{-1}_2\tX^{-1}_1} = S\tau^{-1} \cdot \tau S^{-
1}_{\tX^{-1}_2} = SS^{-1}_{\tX^{-1}_2} = ST^{-1}.$$

On the other hand, by $(1_a)$ and (2), $S_{\tX^{-1}_2\tX^{-1}_1} = S^{-1}T.$
Thus $$\tau S_{X_{2}\1}\1X_1\1=\tau T\1S.$$

We compare the last 2 expressions to get
$\tau T^{-1}S = ST^{-1}$, or
$$
\tau = ST^{-1}S^{-1}T, \quad\text{or} \quad T^{-1}_{S^{-1}} =
\tau T^{-1}.
\tag1.37$$

By Lemmas 3 and 6 $  Q  =
Q_{\tX^{-1}_1\tX^{-1}_2\tX^{-1}_3} = [S_{\tX^{-1}_1\tX^{-1}_2\tX^{-1}_3},
(f^{-1}_{13})_{\tX^{-1}_1\tX^{-1}_2\tX^{-1}_3}].$
Thus:
$$Q= [S_{\tX^{-1}_1\tX^{-1}_2\tX^{-1}_3}, f^{-1}_{24}]
\tag 1.38 $$
(we use $(\{1,3\})\psi(\tX^{-1}_1\tX^{-1}_2\tX^{-1}_3) =
(\{1,3\}) (1\quad 2)(2\quad 3)(3\quad 4) = \{4,2\})$.

Considering a quadrangle formed by $a_2,a_3,a_4,a_5$ (see Fig. III.7.2)

\midspace{2.00in}
\caption{Fig. III.7.2}

\flushpar we can write in $\Tilde{\Tilde{P}}_n$ (Lemma III.2.2)  $f_{35}f_{24}
=
f_{25}f_{34}$, or $f_{24} = f^{-1}_{35}f_{25}f_{34}$, $f^{-
1}_{24} = f^{-1}_{34}f^{-1}_{25}f_{35}$.  From (1.38) we get,
denoting by $b = \tX^{-1}_1\tX^{-1}_2\tX^{-1}_3$,
$$\align Q &= [S_{\tX^{-1}_1\tX^{-1}_2\tX^{-1}_3},f^{-1}_{24}] = [S_b, f^{-
1}_{34}f^{-1}_{25}f_{35}]\tag1.39\\
&= [S_b,f^{-1}_{34}][S_b,f^{-1}_{25}f_{35}]_{f_{34}} = Q_{b,3,4}
\cdot [S_b,f^{-1}_{25}]_{f_{34}} \cdot
[S_b,f_{35}]_{f_{25}f_{34}}\\
&
= Q_{b,3,4} \cdot (Q_{b,2,5})_{f_{34}} \cdot [S_b,f^{-1}_{35}]^{-
1}_{f^{-1}_{35}f_{25}f_{34}}= Q_{b,3,4} \cdot (Q_{b,2,5})_{f_{34}} \cdot
(Q_{b,3,5})_X\1. \endalign$$
Now,  $(\{1,2\})\psi(b) = \{1,2\}(1\quad 2)(2\quad 3)(3\quad 4) =
\{4,1\}$.  Since $\{4,1\} \cap \{2,5\} = \emptyset$ and
$\{4,1\} \cap \{3,5\} = \emptyset$, we get by   Lemma 2 that
$Q_{b,2,5} = Q_{b,3,5} = 1$, and by (1.39)
$$Q = Q_{b,3,4}.$$

We can write: $S_{\tX^{-1}_1\tX^{-1}_2\tX^{-1}_3} = (S^{-1}\tau)_{\tX^{-
1}_2\tX^{-1}_3} = S^{-1}_{\tX^{-1}_2\tX^{-1}_3} \tau = T^{-1}_{\tX^{-
1}_3}\tau$ (using $\tau_{\tX_i} = \tau$ $\forall i$).  So $Q =
Q_{b,3,4} = [S_{\tX^{-1}_1\tX^{-1}_2\tX^{-1}_3},f^{-1}_{34}]
\overset\text{by}\ f_{34} = \tX^2_3\to = [T^{-1}_{\tX^{-1}_3} \tau, \tX^{-2}_3]
=
[\tau, \tX^{-2}_3]_{T_{\tX^{-1}_3}} \cdot [T^{-1}_{\tX^{-1}_3}, \tX^{-
2}_3] = [T^{-1}_{\tX^{-1}_3}, \tX^{-2}_3] =  [T^{-1},\tX^{-2}_3]_{\tX^{-
1}_3}$.  Since $Q_{\tX_3} = Q$, we get
$$
Q = [T^{-1},\tX^{-2}_3].
\tag1.40$$

This implies $Q_{S^{-1}} = (T^{-1}T_{\tX^2_3})_{S^{-1}}\overset\text{
assumption 3}\to = T^{- 1}_{S^{-1}} \cdot (T_{S^{-1}})_{\tX^2_3} =
[T^{-1}_{S^{-1}},\tX^{- 2}_3]$\linebreak $ \overset\text{by (1.37)}\to= [\tau
T^{-1}, \tX^{-2}_3] = [T^{-1},\tX^{- 2}_3]_{\tau^{-1}} [\tau,\tX^{-2}_3]
\overset\text{by Lemma 4}\to= [T^{-1},\tX^{-2}_3]_{\tau^{- 1}} \overset\text{by
(1.40)}\to= Q_{\tau^{-1}}.$

Using $Q=\tau\1$ we get $\tau_{S\1}\1=\tau_{\tau\1}\1.$
Thus, $\tau_{S^{-1}} = \tau$ and $\tau_S = \tau$.   \newline\qed\quad for
Lemma
8.

We can now finish the proof of Proposition III.7.1.

By Lemma III.2.4, we only have to prove that $\tau^2 = 1$,
$\tau_b = \tau$ $\forall b \in \Tilde{\Tilde{B}}_n$ and $\tau
\in\operatorname{Center}(G)$.

By the previous Lemma, $\tau_S = \tau$, and by Lemma 7, $\tau_S = \tau^{-1}$.
Thus, $\tau = \tau^{-1}$ and    $\tau^2 = 1$.

By Lemma 6, $\tau_{\tX_i} = \tau$
$\forall i \in (1, \ldots, n - 1)$.  Thus $\tau_b = \tau$ $\forall b \in
\Tilde{\Tilde{B}}_n$.

 By Lemma 8
$\tau_S = \tau$, i.e. $[\tau,S]=1.$
Let $b\in\tB_n:$ \ $[\tau,S_b] = [\tau_{b^{-1}},S]_b =
[\tau,S]_b  = 1$.

Thus $\tau$ commutes with $S_b$\ $\forall b\in \Tilde{\Tilde
B}_n.$ Since $
\ker(\tB_n\ri\Tilde{\Tilde B}_n)$ acts trivially on $G$,\linebreak
$\Tilde{\Tilde B}_n$   acts on $G$ via $\tilde B_n,$ and thus $\tau$ commutes
with $S_b$\ $\forall b\in \tB_n.$

 By assumption (0) of the proposition, $G$ is
generated by $\{S_b\}_{b\in\tB_n}.$ Thus   $\tau \in\Center(G).$

\qed \quad for Proposition III.7.1\edm

\end

\magnification=1200
\parindent 20 pt
\NoBlackBoxes
\define\Center{\operatorname{Center}}
\define \a{\alpha}
\define \be{\beta}
\define \Dl{\Delta}
\define \dl{\delta}
\define \g{\gamma}

\define \G{\Gamma}
\define \lm{\lambda}

\define \r{\rho}
\define \s{\sigma}\define \si{\sigma}

\define \ve{\varepsilon}
\define \vp{\varphi}

\define \cd{\cdot}
\define \df{\dsize\frac}

\define \fa{\forall}

\define \iy{\infty}
\define \la{\langle}
\define \ra{\rangle}

\define \ri{\rightarrow}
\define \Ri{\Rightarrow}

\define \ub{\underbar}
\define \un{\underline}
\define \ov{\overline}

\define \edm{\enddemo}
\define \ep{\endproclaim}

\define \sk{\smallskip}
\define \mk{\medskip}
\define \bk{\bigskip}

\define \1{^{-1}}
\define \2{^{-2}}

\define \hL{\hat{L}}

\define \CP{\Bbb C\Bbb P}

\define \CPt{\Bbb C\Bbb P^2}

\define \BR{\Bbb R}
\define \BZ{\Bbb Z}

\define \tB{\tilde{B}}
\define \tE{\tilde{E}}
\define \tT{\tilde{T}}
\define \tX{\tilde{X}}
\define \tY{\tilde{Y}}
\define \tP{\tilde{\Pi}}
\define \tp{\tilde{P}}
\define \tv{\tilde{v}}
\define \tZ{\tilde{Z}}
\define \tz{\tilde{z}}

\define \CG{\Cal G}

\define \Aff{\operatorname{Aff}}
\define \Gal{\operatorname{Gal}}

\define \Ss{S^{(6)}}
\define \Xab{X_{ab}}

\heading{\bf CHAPTER  IV. \ \ New Set of Generators for $G$}\rm\endheading
\baselineskip 20pt

Recall from Chapter II that $G=G(\ve_{18}) = \pi_1 (\CPt
 - S_3, *)$ is generated by $E_i, E_{i}'$ and satisfies the relations listed in
Theorem II.6.

In this chapter we shall introduce new generators for $G,$ using the braid
group $B_9$ and the quotient $\tB_9$ from Chapter III.
\bk
\subheading{\S1. New presentation of $\bold{B_9}$}
\demo{Definition}\ $\un{T_i \quad i=1,\dots, 9\ \ i\neq 4}$\edm

Consider a geometric model $(D,k)$ for $\#K=9$ as in Fig. IV.1.1.

Let $\left\{ t_i\right\}^0_{i=1\ i\neq 4}$ be paths connecting different
parts of $K$ as in Fig. IV.1.2.

Let $T_i$ be the half-twist corresponding to $t_i\quad i=1,\dots,9\ i\neq 4.$\
$(T_i=H(t_i))$

\midspace{1.00in}\caption{Fig. IV.1.1
\qquad\qquad\qquad\qquad\qquad\qquad\qquad Fig. IV.1.2}

\proclaim{Lemma IV.1.0}

$T_i$ and $T_j$ are adjacent for $(i,j)$ as follows:

$\qquad \quad i,j \in \{1, 2, 3\}$

$\qquad \quad i=5\qquad \qquad j=3,8,9$

$\qquad\quad i=6,7,8 \qquad j=i+1.$

$T_i$ and $T_j$ are disjoint for $(i,j)$ as follows:

$\qquad\quad i\in \{1,2,3\}\ \ \ j\in \{6,7,8,9\}$

$\qquad\quad i=5\qquad \quad\ \ j= 1,2,6,7$

$\qquad\quad i=6\qquad \quad\ \ j= 8,9$

$\qquad\quad i=7\qquad \quad\ \  j=9.$\ep
\sk
\demo{Proof}\ From Fig. IV.1.2. \qed\edm
\sk
\demo{Remark}\ The choice of the model comes from the configuration of planes
in the degeneration of $V_3$ to the union of planes as in Fig.II.1.
(We constructed this degeneration in BGT III \cite{MoTe7}.)
In each of the triangles we choose a point (Fig. ~IV.1.3).

\midspace{1.50in}\caption{Fig. IV.1.3}

We choose a path connecting 2 points in neighboring triangles as in Fig.
IV.1.4.

\midspace{1.25in}\caption{Fig. IV.1.4}

We then get a configuration which is basically equal to the one in Fig. IV.1.2.

Since we do not need all possible connections to get a set of generators for
$B_9,$ we skip the connection between the points in $P_3$ and $P_5.$\edm

\proclaim{Lemma IV.1.1}

There exists a presentation of the braid group $B_9$, as follows:
$$B_9 = \la T_{i.} \bigm| i=1,\dots, 9,\ \ i\ne 4\ra$$
and the following is a complete set of relations:
$$\alignat 2
&\la T_i, T_j\ra = 1 \quad &&\text{if}\ T_i\ \text{and}\  T_j\ \text{are
consecutive}\\
&[T_i, T_j] = 1 &&\text{if}\ T_i\ \text{and}\  T_j\ \text{are
disjoint}\\
&[T_1, T_2\1 T_3 T_2]=1 &&\\
&[T_5, T_8\1 T_9 T_8]=1.\endalignat$$\ep
\sk
\demo{Proof}\ Consider the geometric model $(D,K),$\ $ \# K=9$ as in Fig.
IV.1.1.
We choose a frame in $B_9 [D,K]$ where each half-twist in the frame corresponds
to a path, as in Fig. IV.1.5.\edm

\midspace{1.00in}\caption{Fig. IV.1.5}

In terms of $T_i,$ this frame is
$$T_1, T_2, T_2\1 T_3T_2, T_5, T_9, T_9 T_8 T_9\1, T_7, T_6.$$

By E. Artin's presentation of the braid group (see Chapter 0), we know that
$B_9$ is generated by the above frame and the only relations are triple
relations for non-neighboring elements and commutation  relations for
neighboring elements. Thus, a full set of relations is:
$$\align
\la T_1,T_2\ra&=\la T_2,T_2\1T_3T_2\ra = \la T_2\1T_3T_2,T_5\ra=\la T_5,
T_9\ra\\ &= \la T_9,T_9T_8T_9\1\ra = \la T_9T_8T_9\1, T_7\ra = \la T_7,T_6\ra =
1\endalign$$ and all possible commutation relations between other elements of
the frame.

Since $T_3=T_2\ (T_2\1 T_3T_2)\ T_2\1$ and $T_9=T_9\1\ (T_9T_8T_9\1)\ T_9,$
then
\ \ $T_1,\ T_2,\ T_3,\ T_5\ T_9,$ $T_8,\ T_7,\ T_6$ generates $B_9.$
Translating the above relations to these generators and using simple facts
about commutators, we obtain the Lemma.\qquad \qed

\sk
\demo{Definition}\ $\underline{T_4}$.

$T_4=T_2\1\ T_3\ T_7\1\ T_8\ T_5\ T_8\1\ T_7\ T_3\1\ T_2.$
It is possible to notice that $T_4$ is the half-twist that corresponds to the
path $t_4$ as in Fig. IV.1.6, and thus $T_4$ is adjacent to $T_2$
and $T_6,$ transversal to $T_7$ and $T_3$ and disjoint from the others.\edm

\midspace{1.00in}\caption{Fig. IV.1.6}

\bk
\subheading{\S 2. Presentation of $\bold{\tilde B_9}$}

Let $\tB_9$ be as in Chapter III, $\tB_9= B_9/T,$ where $T=\la[X, Y]\ra,$ and
$X, Y$ are transversal.

Let $T_i$ be as in \S 1.
Let $\tT_i$ be the images of $T_i$ in $\tB_9.$
\sk
\proclaim{Lemma IV.2.1}

$\tB_9$ is generated by $\la \tT_i \bigm| i=1\dots 9\ra$ and the only
relations are:
\roster
\item\ $\la \tT_i, \tT_j\ra = 1 \qquad\quad T_i, T_j$\ are consecutive $i, j\ne
4$ \item\ $[\tT_i, \tT_j] = 1 \qquad\quad\ T_i, T_j$\ are disjoint $i, j\ne 4$
\item\ $[\tT_1, \tT_2\1 \tT_3 \tT_2] = 1$
\item\ $[\tT_5, \tT_8\1 \tT_9 \tT_8] = 1$
\item\ $\tT_4= \tT_2\1 \tT_3 \tT_7\1 \tT_8 \tT_5 \tT_8\1 \tT_7 \tT_3\1
\tT_2.$\endroster\ep
\sk
\demo{Proof} \ By the previous Lemma, $\tB_9$ is generated by $\tT_i,\ i\ne 4$
with a full set of relations $(1)\dots (4).$
We add one generator and express it
in terms of the other relations to get $\tT_4$ and relation (5).\quad $\square$
\edm
 \sk
\proclaim{Lemma IV.2.2}
\roster
\item\ $\la\tT_4, \tT_2\ra = 1.$
\item\ $\la\tT_4, \tT_6\ra = 1.$
\item\ $\la\tT_7, \tT_6\ra = 1.$
\item\ $[\tT_4, \tT_i] = 1\qquad i=1,\ 3,\ 5,\ 9,\ 8.$
\endroster\ep
\sk
\demo{Proof}\ We use Theorem III.3.1 from BGT I.

Since $T_4$ is consecutive  to $T_i,$\quad $i=2,6$ and $T_6$ is consecutive to
$T_7$ then $\la T_4,T_2\ra = \la T_4,T_6\ra= \la T_6,T_7\ra
= 1.$
Thus, $\la\tT_4,\tT_2\ra = \la T_4,T_6\ra = \la \tT_6,\tT_7\ra = 1.$

$T_4$ and $T_i$ for $i=1,9$ are disjoint, therefore, $[T_4, T_i] = 1$ for
$i=1,9.$ Thus, $[\tT_4, \tT_i] = 1$ for $i=1,9.$ The half-twists
$T_4$ and $T_i\quad i=3, 5, 8$ are transversal and, thus, $[\tT_4, \tT_i] = 1$
for $i=3,5,8$ (Remark III.1.1).\quad $\square$\edm
\sk
We need the following relations of $\tT_i$ in order to get a smaller set of
generators for $G.$

\proclaim{Lemma IV.2.3}
$$\alignat2 &(\tT_4)_{\tT_2\1\tT_3\tT_7\1\tT_8}=\tT_5\quad &&\text{preserving
polarization}\\
  &(\tT_2)_{\tT_4\tT_3\tT_5\1\tT_7\1}=\tT_8\quad &&\text{preserving
polarization}\\
 &(\tT_3)_{\tT_2\1\tT_4\tT_5\1\tT_8 }=\tT_7\quad &&\text{preserving
polarization}\endalignat$$\ep
\demo{Proof}
It is actually true  for $T_i$ instead of $\tT_i.$
It can be verified geometrically using Fig. III.1.1 for a geometric
presentation of a half-twist conjugated by another half-twist. \qed \edm
\proclaim{Lemma IV.2.4}\ Another presentation of $\tB_9.$

Let $\tT_{1'} =\tT_2^{+2}\ \tT_1\ \tT_2^{-2}\quad \tT'_i = \tT_i\ \ i\ne 1.$
Then $\tilde B_9$ is generated by $\tT_i'$ and the following is a complete set
of relations: \roster
\item\ $\la\tT_i', \tT_j'\ra=1$\ \ if $\tT_i, \tT_j$\ \ are adjacent \quad
$i,j\neq 4.$
\item\ $[\tT_i', \tT_j']=1$\ \ if $\tT_i, \tT_j$\ \ are disjoint \quad
$i,j\neq 4.$
\item\ $[\tT_1', \tT_2' \tT_3'\ \tT_2'{}^{-1}]=1$
\item\ $[\tT_5', \tT'{}^{-1}_8\ \tT_9' \tT_8]=1$
\item\ $\tT_4'= \tT_{2}'{}\1 \tT_3' \tT_7'{}{\1} \tT_8' \tT_5'
\tT_8'{}{\1} \tT_7' \tT_3'{}{\1} \tT_2'.$\endroster\ep
\sk

\demo{Proof}\ Clearly, $T_i=T_i'\ i\neq 1$\ \ $T_1=\tT_{2}'{}^{-2}
\tT_1' \tT_2'{}^{2}.$ We substitute these expressions in the relations of Lemma
IV.2.1 to prove the Lemma. \qed\edm
\bk
\subheading{\S 3. $\bold{\a:\tilde B_9\rightarrow G}$}

We want to prove that there exists $\a: \tB_9 \ri G$ s.t. $\a (\tT_i) = E_i,$
for $E_i$ that were introduced in Chapter II. For that we prove certain
relations that $E_i$ satisfy, based on Proposition II.6. \sk
\proclaim{Lemma IV.3.1}
\roster
\item\ $\la E_i, E_j\ra=1$\ \ if $T_i, T_j$\ \ are adjacent \quad
$i,j\neq 4.$
\item\ $[E_i, E_j]=1$\ \ if $T_i, T_j$\ \ are disjoint \quad\
$i,j\neq 4.$
\item\ $[E_1, E_2 E_3 E_2\1]=1$
\item\ $[E_9, E\1_5 E_8 E_5]=1$
\item\ $E_4= E_{2}\1 E_3 E_7\1 E_8 E_5 E_8\1 E_7 E_3\1 E_2.$
\endroster\ep
\sk
\demo{Proof}\  We use Proposition II.6, which states a list of relations
satisfied by the  $E_i.$
\roster
\item"(1), (2)" By Fig. IV.1.3, $T_i, T_j\ \ (i,j\ne 4)$ are adjacent,
$\Leftrightarrow \hL_i$ and $\hL_j$ are edges of some triangle,
$ \Leftrightarrow$ (by Corollary II.3) $\psi (E_i)$ and $\psi(E_j)$ have one
common index. Moreover,  $T_i$
and $T_j$ are disjoint $\Leftrightarrow \psi (E_i)$ and $\psi (E_j)$ are
disjoint.
By (1) and (2) of Proposition II.6  we get (1) and (2) of this Proposition.
\item"(3)"   By Proposition II.6 (3),\
$E_8 = E_7 E_5 E_3\1 E_4\1 E_2 E_4 E_3 E_5\1 E_7\1.$\newline
By Proposition II.6 (2), $[E_3, E_4] = 1.$  \newline
Thus, $E_8 = E_7 E_5 E_4\1 E_3\1 E_2 E_3 E_4 E_5\1 E_7\1.$\newline
By Proposition II.6 (1), $E_3\1 E_2 E_3 = E_2 E_3 E_2\1.$\newline
Thus, $E_2\ E_3\ E_2\1\ = E_4\ E_5\1\ E_7\1\ E_8\ E_7\ E_5\ E_4\1.$\newline
By Proposition II.6 (2) $E_i$, for $i = 4,5,7,8,$ commutes with $E_1$.\newline
Thus, $[E_1, E_2\ E_3\ E_2\1] = 1.$
\item"(4)" By (3), $E_8=E_7 E_5 E_3\1 E_4\1 E_2 E_1 E_3 E_5\1
E_7\1.$\newline
By (2), $[E_5, E_7] = 1.$ newline
Thus,\ $E_5\1 E_8 E_5 = E_7 E_3\1 E_4\1 E_2 E_4 E_3 E_7\1.$
\newline
For $i=7,  3, 4,2,$\
$E_i$ commutes with $E_9$.
Thus, $[E_9, E_5\1 E_8 E_5]=1.$
\item"(5)" By Proposition II.6 (3),\
$ E_4\1 E_2 E_4 = E_3 E_5\1 E_7\1 E_8 E_7 E_5 E_3\1.$\newline
 By Proposition II.6 (2)(1),\
 $E_4\1 E_2 E_4 = E_2 E_4 E_2\1$\ \ and\ \ $[E_5, E_7] = 1.$\newline
 By Proposition II.6 (1)\ $E_5\1 E_8 E_5 = E_8 E_5 E_8\1.$\newline
Thus, $E_2\ E_4 E_2\1 = E_3 E_7\1 E_5\1 E_8 E_5 E_7
E_3\1.$\newline
$E_4 = E_2\1 E_3E_7\1 E_8 E_5 E_8\1 E_7 E_3\1
E_2.$ \qquad \quad \qed\endroster\edm
\medskip
We first prove that there exist $\a':\tB_9\ri G$ such that $\a'(\tT'_i) = E_i.$
\sk
\proclaim{Lemma IV.3.2}

There exists a homomorphism $\a',$\ \ $ \a': \tB_9\ri G$ such that $\a'(\tT'_i)
= E_i.$\ep

\demo{Proof}\ By Lemma IV.2.4, $\tB_9$ is generated by $\tT'_i.$
Thus, we define
 $\a' (\tT'_i) = E_i.$
To prove that $\a'$ induces a homomorphism, we have to show that $E_i$
satisfies
 any relation
that $\tT'_i$ satisfies.
In IV.2.4 we presented a full list of relations for $\tT'_i.$
In Lemma IV.3.1 we proved that these relations are satisfied when $\tT'_i$ is
replaced by $E_i.$ \quad \qed\edm
\sk
\proclaim{Lemma IV.3.3}
$$[E_1,\ E_2\1\ E_3\ E_2] = 1.$$\ep
\sk
\demo{Proof}\ $T'_3$ is transversal to $T'_2\ T'_1\ T_2\1.$  (See Fig.
IV.3.1)\edm

\midspace{1.00in}\caption{Fig. IV.3.1}

Thus, $[\tT'_2 \tT'_1 \tT_2\1, \tT_3'] = 1.$

Thus, $\a' [\tT'_2 \tT'_1 \tT_2\1, \tT_3'] = 1.$

Thus, $[E_2 E_1 E_2\1, E_3] = 1.$

Thus, $[E_1,\ E_2\1\ E_3 E_2] = 1.$ \qquad \qquad \qed
\sk
\proclaim{Lemma IV.3.4}
There exists $\a:\tB_9\ri G$ such that $\a (\tT_i) = E_i.$\ep
\sk
\demo{Proof}\ We use the presentation of $\tB_9$ from Lemma IV.2.1 where
$\tB_9$
is generated by $\tT_i.$
We define $\a(\tT_i) = E_i.$
The relations listed in IV.2.1 are satisfied when $\tT_i$ is replaced by $E_i$
by Lemma IV.3.1 and Lemma IV.3.3.
Thus we can extend the definition of $\a$ to the whole of $\tB_9$ in a natural
way.\edm
\bk

\subheading{\S 4. Prime elements in $\bold{B_9}$}

We now recall a few results from Chapters II and III concerning the
braid group $B_n$ and the quotient  group $\tB_n.$
We refer the reader to Chapter III for the definition of prime element with
s.h.t. (supporting half-twist) $T,$ and central element $c.$
\bk
We quote here a few results from  \cite{MoTe4}, Chapter II and Chapter III.
\proclaim{Lemma IV.4.0}\
If $X$ and $Y$ are 2 consecutive half-twists in $B_n$, then
\roster
\item"(a)" $XYX = YXY.$
\item"(b)" $XYX\1 = Y\1XY.$
\item"(c)" $(Y)_{X\1} = (X)_Y$.
\item"(d)" $(Y)_{X\1}$ is consecutive to $X$ and to $Y.$  It is the half-twist
corresponding to a path  connecting those ends of $X$ and $Y$ which are not  a
common index of $X$ and $Y$. \item"(e)"
$u=\tX_{\tY\1}^2\tY\2=\tilde Y_{\tilde X}^2\cdot \tilde Y^{-2},$ is an
element of $\tilde P_n;$\ $u$ is a prime element with s.h.t. $\tilde X$ and
central element $c=[\tY^2 , \tX^2]$\ \ (i.e., $c^2=1,\ c\in \Center(\tilde
P_n)$).\item"(f)" $[\tX'{}^2,  \tY'{}^2]
=[\tX'{}^2,\tY\2]=[\tX'{}\2,\tY'{}\2]=c$\ $\fa X',Y'$ a pair of consecutive
half-twists where $c^2=1,\ c\in\Center (G).$ \item"(g)"  $\tilde Y\2 (\tilde
Y^2) _{\tX\1}= c$ (inverse of a prime element of $\tilde P_n$ with s.h.t.
$X).$
\item"(h)" If $Z$ is transversal to $X,$\ $g$ a prime element
with s.h.t. $X,$ then $(g)_Z=g.$\endroster\ep

\demo{Proof}

(a) From Lemma III.3.1 \cite{MoTe4}.

(b) From Lemma II.4 and (a).

(c) From (b).

(d) Let $Y=H(y),$\ $y$ connects $a$ and $b$ and $X=H(x)$\ $x$ connects $b$ and
$c.$ Assume $X$ corresponds to a diffeomorphism $\be: D\ri D$ s.t. $\be(b)=c.$
Then $(Y)_{X\1}=H((\si)\be\1).$
Clearly, $(\si)\be\1$ is a path connecting $a$ and $c.$ (See Claim III.1.0 and
Fig. III.1.1.)

(e) Lemma III.6.1.

(f) Proposition III.5.2.

(g) We shall prove that $c(Y\2(Y^2)_{X\1})\1$ is a prime element:
$c(Y\2(Y^2)_{X\1})\1=c(Y\2)_{X\1}Y^2=Y^2(Y\2)_{X\1}$ (by (b)).
Let $T=(Y)_{X\1}.$
By (e) and (c), $(T)_X^2\cdot T\2$ is a prime element.
But $(T)_X=((Y)_{X\1})_X=Y.$
Thus $Y^2(Y)_{X\1}\2$ is a prime element.

(h) Lemma III.4.3. \qed\edm\medskip

\subheading\nofrills{$\bold \xi$-situation:}

Consider $\tB_9, \tilde P_9$ as $\tB_9$-groups by conjugation.

Let $\tT_1,\dots, \tT_9$ be as in \S 1.

We choose a polarization on $\tT_i$ from smaller end index to bigger end index
as shown in Fig. IV.4.1.

\midspace{1.50in}\caption{Fig. IV.4.1}

Let $\xi_1 = (\tT_2)^2_{\tT_1} (\tT_2)^{-2}.$

By Lemma IV.4.0 (e), $\xi_1$ is a prime element in $\tp_9$ with
s.h.t. $\tT_1$ and central element $[T_2^2,T_1^2].$

 Let $c$ be the corresponding central
element.
Thus, $c\in\Center (\tp_9), \quad c^2 = 1,$\ $c=[T_2^2,T_1^2].$

Let $\xi_i$ be the unique prime element in $\tp_9$ s.t. $(\xi_i,\tT_i)$ is
coherent with $(\xi_1, \tT_1)$\ $i = 2,\dots, 9.$ (See Proposition III.3.3).
\bk
\proclaim{Claim IV.4.1}
\roster
\item $c$ is the corresponding
central element of $(\xi_i, T_i)\ \ \fa i = 1\dots 9,$\ $c^2=1.$
\item $c\in\Center (\tB_9).$
\item $\xi_i =(Y)^2_{\tT_i} Y^{-2}$ for some $Y$ half-twist adjacent to
$T_i.$
\item $c=[\tT_k^2, \tT_\ell^2]\ \ \fa k,\ell$ s.t. $T_k$ and $T_\ell$
are consecutive.
\item $\xi_i$ is a prime element of $\tB_9.$
\item Let $X, Y, $ be 2 half-twists, \
$X=H(x),\quad Y=H(y), $\quad $T_i=H(t_i)$ s.t. $x,y,t_i$ make a triangle.
Assume that $x$ and $y$ meet in $v,$ and a counter clockwise
rotation around $v$ inside the triangle meets $x$ before it meets
$y$ (Fig. \rom{IV.4.2(a)} and \rom{(b)}).  Then $\tT_i,$ the s.h.t. of $\xi_i$
satisfies $T_i=XYX\1.$  And:
\newline \text{\rm(i)} If the polarization of $T_i$ goes from $x$ to ${y},$
then
$\xi_i=\tilde X^2\ \tY^{-2}.$
\newline \text{\rm(ii)} If the polarization of $T_i$
goes from $y$ to $x,$ then $\xi_i=X\2Y^2.$\endroster\ep

\demo{Proof}
\roster
\item The pair $(\xi_i, \tT_i)$ is coherent with $(\xi_1, \tT_1).$
Thus, $\xi_i$ is conjugate to $T_1$ and $\tT_i$ is conjugate to $f_1$ by some
$B_i.$
Then by Lemma III.3.1, their corresponding central element is equal.
\item A priori, $c\in \Center (\tp_9).$
We have to prove that $c\in\Center (\tB_9).$
Consider $\tp_9$ as $ \tilde B_9$-group.
$\xi_i$ is a prime element in $\tp_9$ as a a $\tilde B_9$-group where $c$ is
the
central element of $(\xi_i, T_i).$
Thus, we have $(c)b=c\ \ \fa
b\in\tB_9.$
$((\quad)b =$ action of
 $\tB_9$ on $\tp_9$), but $(c) b = c_b$ by definition.
Thus, $c_b=c\ \ \fa b\in \tB_9 \Rightarrow c\in\Center (\tB_9).$
\item The pair $(\xi_i, T_i)$ is coherent with the pair $(\xi_1, \tT_1).$
Thus, $\exists b_i\in B$ s.t. $\xi_i=(\xi_1)_{b_i}$ and $  \tT_i =
(\tT_1)_{b_i}.$ Denote $(\tT_2)_{b_i} = \tY$ and apply conjugation by $b_i$ on
$\xi_1=(\tT_2^2)_{T_1}\tT_2^{-2}.$ \item By Proposition III.5.2.
 \item By
construction, $\xi_i$ is a prime element of $\tp_9$ with central element $c.$
In (2) we proved that $c\in\Center(\tB_9).$
Thus, $\xi_i$ is also a prime element of $\tB_9.$
\item Consider the two triangles (Fig. IV.4.2(a) and (b)):

\midspace{.75in}\caption{Fig. IV.4.2}

If the polarization of $T_i$ goes from $X$ to $Y$ as in Fig. (a) (from $Y$ to
$X$ as in Fig. (b)) then we take $b\in\tB_n$ s.t. $(T_1)b=\tT_i$\ $(T_2)b=Y$
preserving polarization of $T_1$ (reversing polarization),
respectively.
Clearly, $((T_2)_{T_1})b=X.$
Consider the polarized pair $(\xi_1,\tT_1)=((\tT_2)_{T_1}^2\tT_2\2,\tT_1).$
We apply $b$ on it to get a coherent (anti-coherent) polarized pair
$(\tX^2\tY\2,\tT_i).$
Recall that $(\xi_i,\tT_i)$ is coherent with $(\xi_1,\tT_1).$
Thus by Proposition III.3.3 (or by Lemma III.4.1) $\tX^2\tY\2=\xi_i$ (or
$\tX^2\tY\2=\xi_i\1c),$ respectively. To get $\xi_i=\tX\2\tY^2$ from
$\tX^2\tY\2=\xi_i\1c$ we use (4) above. \quad$\square$\endroster\edm
\sk The next Corollary is
technical in nature, to be used later in order to obtain a smaller set of
generators for $G.$\sk \proclaim{Corollary IV.4.2}  \roster \item $\xi_1 =
(\tT_2)^2_{\tT_1\1}\tT_2\2.$ \item
$\xi_2 = \tT_1\2 (\tT_1)_{\tT_2\1}^2$\quad $\xi_2=T_4^2(T_4)_{T_2\1}\2\quad
\xi_2=\tT_3\2(\tT_3)_{\tT_2\1}^2$
 \item $\xi_3 =  \tT_1^2\ (\tT_1)\2_{\tT_3\1}$
\item $\xi_4=(\tT_4)_{\tT_6\1}^2\tT_6\2$
\item $\xi_5=\tT_8\2(\tT_8^2)_{\tT_7\1}$ \item $\xi_6 =
(\tT_4)\2\ (\tT_4^2)_{\tT_6\1},$ \quad $\xi\1_6 = c\tT_7^2 (\tT_7)_{T_6\1}\2$
 \item $\xi_7=\tT_8\2(\tT_8)_{\tT_7\1}^2$  \item  $\xi_8 = \tT_9\2
(\tT_9)_{\tT_8\1}^2$
 \item $\xi_9 =
\tT_5^2(\tT_5)_{\tT_9\1}\2 $ \item $(\xi_4)_{\tT_2\1\tT_3 \tT_7\1 \tT_8} =
\xi_5$ \item $(\xi_2)_{\tT_4 \tT_3\tT_5\1 \xi_7\1} = \xi_8$
\item$(\xi_3)_{\tT_2\1 \tT_4 \tT_5\1 \tT_8} = \xi_7$\endroster\ep
\demo{Proof}
In the entire proof we use Lemma IV.4.0(f) to interchange squares of
consecutive half-twists, and multiplying the product by $c.$
We also use the facts that $c^2=1,$\
$c\in\Center(G).$
The main tool is Lemma IV.4.1(6).  We also use Lemma IV.4.0(d) to present an
edge of a triangle as a conjugation of the other 2 edges.

(1) By definition.

(2) Consider the triangle from Fig. IV.4.3(a):

\midspace{1.00in}\caption{Fig. IV.4.3}

\flushpar Since the polarization of $T_2$ goes from $Y$ to $X,$ we apply Lemma
IV.4.1(6) to get:
$$\xi_2=X\2Y^2=\tT_1\2(\tT_1)_{\tT_2\1}^2.$$
Consider the triangle from Fig. IV.4.3(b).

\flushpar Since the polarization of $T_2$ goes from $X$ to $Y,$ we apply Lemma
IV.4.1(6) to get
$$\xi_2=X^2Y\2=\tT_4^2 (\tT_4)_{\tT_2\1}\2.$$

\flushpar Consider the ``triangle'' from Fig. IV.4.3(c)

\flushpar Since the polarization of $T_2$ goes from $Y$ to $X,$ then:
$$\xi_2=\tT_3\2(\tT_3)_{T_2\1}^2.$$

(3) Consider the triangle from Fig. IV.4.4:

\midspace{1.00in}\caption{Fig. IV.4.4}

\flushpar Since the polarization of $T_3$ goes from $Y$ to $X,$ then:
$$\xi_3=X\2Y^2=(T_1)_{T_3}\2T_1^2.$$
We apply $\tT_3\2$ on the above equation to get:
$$(\xi_3)_{\tT_3\2}=(\tT_1)\2_{\tT_3\1}(\tT_1)^2_{\tT_3\2}.$$

By Lemma III.2.1, $(\xi_3)_{\tT_3\2}=\xi_3.$

By Lemma IV.4.1(4), $(\tT_1^2)_{\tT_3\2}=cT_1^2.$

Thus, $$\align \xi_3&=c(T_1)_{T_3\1}\2T_1^2\\
&=T_1^2(T_1)_{T_3\1}\2\quad \text{(Lemma IV.4.0(5))}.\endalign$$

(4) Consider the triangle from Fig. IV.4.5:

\midspace{1.00in}\caption{Fig. IV.4.5}

\flushpar Since the polarization of $T_4$ goes from $X$ to $Y$
$$\xi_4=(\tT_4)_{\tT_6\1}^2\tT_6\2.$$

(5) Consider the triangle from Fig. IV.4.6:

\midspace{1.00in}\caption{Fig. IV.4.6}

\flushpar Since the polarization of $T_5$ goes from $X$ to $Y$
$$\xi_5=(\tT_8)_{\tT_5}^2\tT_8\2.$$

(6) Consider the triangle from Fig. IV.4.7(a):

\midspace{1.00in}\caption{Fig. IV.4.7}

\flushpar Since the polarization of $T_6$ goes from $Y$ to $X$
$$\xi_6=\tT_4\2(\tT_4^2)_{T_6\1}.$$

Consider the ``triangle'' from Fig. IV.4.7(b).

\flushpar Since the polarization of $T_6$ goes from $Y$ to $X$
$$\xi_6=\tT_7\2(\tT_7^2)_{T_6\1}\Ri\xi_6\1=c\tT_7^2(\tT_7\2)_{\tT_6\1}.$$

(7) Consider the ``triangle'' from Fig. IV.4.8:

\midspace{1.00in}\caption{Fig. IV.4.8}

\flushpar Since the polarization of $T_7$ goes from $Y$ to $X$
$$\xi_7=\tT_8\2(\tT_8^2)_{T_7\1}.$$

(8) Similar to (7).

(9) Similar to (3).

For (10)--(12) we shall only prove the first assertion.
The others are similar.
By Lemma III.2.2,
$((\xi_2)_{\tT_4\tT_3\tT_5\1\tT_7\1},(T_2)_{\tT_4\tT_3\tT_5\1\tT_7\1})$ is a
polarized pair coherent with $(\xi_2,\tT_2).$
By Lemma IV.2.3 $(\tT_2)_{\tT_4\tT_3\tT_5\1\tT_7\1}=\tT_8,$ preserving
polarization.  Thus $((\xi_2)_{\tT_4\tT_3\tT_5\1\tT_7\1},\tT_8)$ is coherent
with $(\xi_2,\tT_2).$
But also $(\xi_8,\tT_8)$ is coherent with $(\xi_2,\tT_2).$
{}From uniqueness, $(\xi_2)_{\tT_4\tT_3\tT_5\1\tT_7\1}=\xi_8.$ \qed
\edm
\bk
\subheading{\S 5. $\bold{\eta}$-situation}

In \S 3 we constructed a homomorphism of groups $\a:\tB_9\ri G$ s.t. $\a
(\tT_i) = E_i.$

We introduce $G$ as a $\tB_9$-group by $(g){\tY} = g_{\a (\tY)} ( = \a(\tY)\1
g\
\a(\tY)).$

The homomorphism $\a$ here then becomes a homomorphism of $\tB_9$-groups.
\bk
\proclaim{Claim IV.5.1} Let $\mu=\a(c),\quad c= [T_2^2,T_1^2],$ then:
\roster
\item"(a)" $\mu = [E_k^2, E_\ell ^2]$ $\fa k, \ell$ s.t. $T_k$ and $T_\ell$ are
adjacent.
\item"(b)" Let $\rho$ be any automorphism of the type $\r_{m_1, m_2, m_3, m_4,
m_6, m_9}$  s.t. $\exists j_0$ with  $m_{j_0}\ne 0$ and all other $m_i = 0.$
Then $(\mu) \rho = \mu.$
\item"(c)" $\mu \in\Center(G).$
\item"(d)" $\mu = \left[ (E_k^2) \rho_k^{m}. (E_\ell)^2 \rho_\ell^n\right] \ \
\fa k,\ell$ s.t. $T_k$ and $T_\ell$ are adjacent, and $(k,\ell) \ne (3, 7),
(2,8), (4,5).$
\item"(e)" $\mu^2 = 1.$
\item"(f)" If $\tX$ and $\tY$ are 2 consecutive half-twists, then
$\a(\tX)^{2} \a(\tY)^{ 2}=\mu\a(\tY)^2\a(\tX)^2.$
\endroster\ep \sk \demo{Proof}
\roster
\item"(a)" $c = [T_k^2, T_\ell^2]\ \fa k,\ell$ s.t. $
\tT_k$ and $\tT_\ell$ are adjacent (Claim IV.4.1  (4)).
Then $\mu=\a (c) = \a [\tT_k^2, \tT_\ell^2] = [E_k^2, E_\ell^2]\ \ \fa k, \ell$
such that  $\tT_k$ and $\tT_\ell$ are adjacent.
\item"(b)"
 $\rho=\rho_1^{m_1} (\rho_2\rho_8)^{m_2} (\rho_3 \rho_2)^{m_3} (\rho_4
\rho_5)^{m_4}\ \rho_6^{m_6} \rho_9^{m_9},$ and $(E_i) \rho = E_i$ if $\rho_i$
appears in  $\rho$ to  the power 0.
If $j_0\neq 1,2,$ then $(E_1)\rho=E_1$ and $(E_2)\rho = E_2.$
We take $\mu = [E_1^2, E_2^2]$ and apply $\rho$ on it to get $(\mu)\rho =
[E_1^2, E_2^2]\rho = [(E_1^2)\rho, (E_2^2)\rho] = [E_1^2, E_2^2] = \mu.$
If $j_0 = 1$ or 2, we apply $\rho$ on $\mu = [E_3^2, E_5^2]$ and continue
similarly to get the result.
\item"(c)" Since $c\in \Center (\tB_9),$ then $\mu \in\Center(\a(\tB_9)).$
Since $\a (\tT_i) = E_i,$ we get $[\mu, E_i] =1\ \ \fa i.$
We want to prove $[\mu, E_{i'}] =1\ \ \fa i$ and then we get
$\mu\in\Center(G).$
($G$ is generated by $\{ E_j,E_{j'}\}^9_{j=1}.$)
Take $\rho=\rho_1^{m_1} (\rho_2\rho_8)^{m_2} (\rho_3 \rho_7)^{m_3}$ $(\rho_4
\rho_5)^{m_4} \rho_6^{m_6}\rho_9^{m_9}$ such that $\rho_i$
appears in $\rho$ to the power 1 and all other $m_i= 0,$
Thus, $(E_i) \rho = E_{i'}$ and there is exactly one $m_{j_0}\neq 0.$
Thus, by $(b)$ $(\mu) \rho =
\mu.$
Thus,
$$1 = (1) \rho = [E_i, \mu] \rho = [(E_i)\rho, (\mu) \rho] = [E_{i'}, \mu].$$
\item"(d)" Since $(k,\ell)\neq (2,8), (3,7), (4,5),$ then $L_k$ and $L_\ell$
are not on the same line.
Let $\rho' = \rho_1^{m_1} (\rho_2 \rho_8)^{m_2} (\rho_3 \rho_7)^{m_2}
(\rho_4 \rho_5)^{m_5} \rho_6^{m_6} \rho_9^{m_6}$ where $\rho_k$ appears
to the power $m$ and all other $m_i=0.$
In particular, $\rho_\ell$ appears to the power $0\cd ((k,\ell) \neq (2,8),
(3,7), (4,5)).$
Then $(E_k) \rho = (E_k) \rho_k^m, (E_\ell) \rho = E_\ell.$
Let $\rho'$ be as above s.t. $\rho_\ell$ appears there to the power $n$ and
all other $m_i=0.$
Then, $(E_\ell) \rho' = (E_\ell) \rho_\ell^n, (E_k) \rho' = E_k.$
Thus, $\mu = (\mu) \rho\rho' = [(E_k)
\rho_k^m, (E_\ell) \rho_\ell^n].$
\item"(e)" $\mu^2 = \a (c)^2 = \a (c^2) =\a(1) = 1.$
\item"(f)"
{}From Lemma IV.4.1 (4). \ \ \ \ \qed\endroster \edm
\sk
\proclaim{Corollary IV.5.2}

$\eta_1=\a (\xi_1)$ is a prime element in $G$ with s.h.t. $\tT_1$ and
central  element  $\mu.$\ep
\sk
\demo{Proof}\ $\xi_1$ is a prime element of $\tB_9$ with s.h.t. $\tT_1$\ (see
Claim IV.4.1). Thus, $\eta_1 = \a(\xi_1)$ is a prime element of $\a (\tB_9)$
with s.h.t. $\tT_1$  and central element $\mu = \a (c),$ from the previous
lemma,  $\mu \in\Center(G).$ Thus, $\eta_i$ is a prime element of $G.$ \qed\edm
\bk
\subheading{$\bold\eta$-situation}

Consider $G$ as a $\tB_9$-group.

Let $\eta_1=\a(\xi_1)$ be a prime element of $G$ with s.h.t. $\tT_1$ and
central
element $\mu=\a(c).$

 Let $\eta_i = L_{\eta_{1}, \tT_{1}} (\tT_i) $  be the unique
prime element of $G$ such that $(\eta_i, \tT_i)$ is coherent with
$(\eta_1, \tT_1).$
\bk
\proclaim{Claim IV.5.3} Consider the $\eta$-situation. Then:
\roster
\item"(a)" $(\eta, \tT_i)$ is a polarized pair of $G$ with s.h.t. $\tT_i$ and
corresponding element $\mu.$ In particular, $\eta_i$ is a prime element with
s.h.t. $\tT_i.$
 \item"(b)" $\eta_i = \a(\xi_i).$
\item"(c)" Every $\eta_i$ is of the form $\a (X^2Y\2)$ where $X, Y$ are
adjacent half-twists and the s.h.t. $T_i$\ is
$X Y X\1$\ $(=Y\1XY).$ \item"(d)"  $\eta_1 = (E_2)^2_{E_1^{-1}} E_2\2.$
\item" " $\eta_2 = E_1\2 (E_1)^2_{E_2\1}.$
\item" "$\eta_3 = E_1^2  (E_1)^{-2}_{E_3\1}.$
\item" "$\eta_4 = (E_4 )^2_{E_6\1}E_6\2.$
\item" "$\eta_5 =(E_8)^2_{E_5}E_8\2.$
\item" "$\eta_6 = (E_4)\2 ({E_4})^{2}_{E_6\1}$\quad  $\eta_6\1 = \mu E_7^2
(E_7)\2_{E_6^{-1}}.$
\item" "$\eta_7 =E_8\2(E_8)^2_{E_7\1}.$
\item" "$
\eta_8=E_9\2(E_9)^2_{E_8\1}.$
\item" "$\eta_9 =  E_5^2(E_5)\2 _{E_9\1}.$
\item" "$\eta_5 = (\eta_4)_{\tT_2\1\tT_3\tT_7\1\tT_8}$
\item" "$\eta_7 = (\eta_3)_{\tT_2\1\tT_4\tT_5\1\tT_8}.$
\item" "$\eta_8 = (\eta_2)_{\tT_4\tT_3\tT_5\1\tT_7\1}.$
\endroster\ep
\sk
\demo{Proof}
\roster
\item"(a)"By the construction of $\eta_i,$
the pairs $(\eta_i,\tT_i)$ are coherent polarized pairs.
$(\eta_1,\tT_1)$ has $\mu$ as central element and
 coherent pairs have the same central element (Corollary III.3.1).
\item"(b)" $\xi_1$ is a prime element with s.h.t. $\tT_i$ s.t. $(\xi_i, \tT_i)$
coherent with $(\xi_1, \tT_1).$
Thus, $\a (\xi_i)$ is a prime element with s.h.t. $\tT_i$ s.t. $(\a(\xi_i),
\tT_i)$
is coherent with $(\a (\xi_1), \tT_1) = (\eta_1, \tT_1).$
{}From uniqueness $\a(\xi_i) =
\eta_i.$
\item"(c)"
{}From Lemma IV.4.1.(6). \item"(d)" Apply $\a$ on  the formulas of Corollary
IV.4.2. \endroster \qed\edm

\bk

\subheading{\S 6. $\bold N$-situation}

\subheading\nofrills{$\bold N$-situation:}

 Let $\CG$ be a $\tB_9$-group.

Assume there exists a $\tB_9$-homomorphism $\lm : \tB_9\ri \CG$ s.t.
\flushpar $(g)b = g_{\lm (b)}\quad \fa b\in \tB_9,\ \forall \ g\in\CG.$

Let $f$  be a prime element in $\CG$ with s.h.t. $\tT_1$ and central
element $\nu.$

Let $f_i=L_{\{(f, \tT_1)\}}(\tT_i)\}$ be the unique prime element s.t.
$(f_i, \tT_1)$ is coherent with $(f_i, \tT_1)$\ $i=1,\dots, 9.$

By Corollary III.3.1, $\nu$
 is the central element of $f_i$\  $\fa i.$

Let $\eta_1=\lm(\xi_1),\ \eta_1$ is a prime element with s.h.t. $\tT_1$ and
central element $\mu = \lambda (c).$

Let $\eta_i=L_{\{\eta_1,\tT_1\}}(\tT_i).$

It is easy to see that, similar to the situation in Claim IV.5.3,
$\eta_i=\lm(\xi_i).$

Let $\mu=\lm(c)$ be the central element of $\eta_i, \ \forall i.$

Let $N_1 = \{\eta_i \bigm| (\eta_i, \tT_i)$ be a polarized pair
coherent with $ (\eta_1, \tT_1),$\ $ i=1\dots 9\}.$

Let ${N_2} = \{f_i \bigm| (f_i, \tT_i)$ be  a polarized pair
coherent with $(f, \tT_1),$ \ $i=1\dots 9\}.$

\newpage

\demo{Definitions}

Let $a, b \in N_1 \cup N_2.$
We say that $a$ and $b$ are weakly disjoint (transversal, disjoint, adjacent,
consecutive or cyclic) if their s.h.t. are weakly disjoint (transversal,
disjoint, adjacent, consecutive or cyclic respectively).
(See definitions in the beginning of Chapter III.)

Let $a\in N_1\cup N_2$ with s.h.t. $\tX.$
Let $\tZ \in \tB_9.$
We say that $a$ and $\tZ$ are weakly disjoint (transversal, disjoint, adjacent,
consecutive, cyclic) if $X$ and $Z$ are weakly disjoint (transversal, disjoint,
adjacent, consecutive, cyclic respectively).\enddemo
\bk
\proclaim{Lemma IV.6.1}\
$a\in N_1\cup N_2,\ \tZ\in \tB_9\ \mu=\lm (c),$ then:

\text{\rm(i)} $a, \tZ\ \text{are adjacent} \Rightarrow a_{\tZ^2} = \cases
a\mu \qquad a\in N_1 \Rightarrow\\
a\nu \qquad a\in N_2.\endcases$

\text{\rm(ii)} $a,\tZ\ \text{are weakly disjoint or commonly supported}
\Rightarrow a_{\tZ^2} = a.$

In other words,

$\left[\lm(\tZ^2),\ a^{\pm 1}\right] =
\left[\lm(\tZ^{-2}),\ a^{\pm 1}\right] = \cases
\mu \quad& a\in N_1,\ \  a,Z\ \text{adjacent}\\
\nu \quad& a\in N_2,\ \  a,Z\ \text{adjacent}\\
1 \quad& a,Z\qquad\quad \text{weakly disjoint or
cyclic}.\endcases$\ep
\sk
\demo{Proof}\ Let $\tX$ be the s.h.t. of $a.$

(i) There exists $v\in P_9$ s.t. $Y=v\1 Zv$\ is
consecutive to $X.$
By the definition of prime element:
$$a_{\tY^{-2}}= \cases
a\mu \qquad a\in N_1\\
a\nu \qquad a\in N_2.\endcases$$
On the other hand, $Y\2 = v\1 Z\2 v Z^2
Z\2 = \left[v\1, Z\2\right] Z\2.$
Since $v\in P_9$ and $Z^2 \in P_9,$ then
$\left[\tv\1, \tZ\2\right] \in \tp_9.$
But $\tp'_9$ is generated by $c$ (Lemma III.5.2),
$c^2=1,$ so $\tY\2 = c^\ve \tZ\2\ \ve = 0, 1.$
Thus, $a_{\tY\2} = a_{\lm(\tY\2)} = a_{\lm
(c^\ve)\lm(\tZ\2)} = a_{\mu^\ve\lm (\tZ\2)} =
a_{\lm(\tZ\2)} = a_{\tZ\2}.$
So, $a_{\tZ\2} = a_{\tY\2}.$
Thus,
$$a_{\tZ\2} = \cases
a\mu \qquad a\in N_1\\
a\nu \qquad a\in N_2.\endcases$$
Since $\nu,\mu\in\Center (G)$
$$a_{\tZ^2}=\cases
a\mu \qquad a\in N_1\\
a\nu \qquad a\in N_2.\endcases$$

(ii) If $a$ and $Z$ are weakly disjoint (cyclic), then:\
There exists $v\in P_9$ s.t. $Y = v\1\ Z\ v$ and $X$ are
disjoint (commonly supported).
As above, $a_{\tY^2} = a_{\tZ^2}.$
By the definition of prime element $a_{\tY^2} = 1.$
Thus, $a_{\tZ^2} = 1.$\qed \edm

\proclaim{Proposition IV.6.2}

\text{\rm(i)} If $a, b\in N_1\cup N_2$ are adjacent, then
$$\left[a, b^{\pm 1}\right] = \cases
\mu \qquad a, b\in N_1\\
\nu \qquad \text{otherwise}.\endcases$$

\text{\rm(ii)} If $a, b\in N_1\cup N_2$ are commonly supported or weakly
disjoint, then
$$\left[ a, b^{\pm 1}\right] = 1.$$\ep
\demo{Proof}\ For $a, b\in N_1$ and $a, b\in N_2,$ the proof follows from Lemma
III.4.2.
Consider the case $a\in N_1,\ b\in N_2.$
Assume $a=\eta_i$ with s.h.t. $\tT_i,$ and $b\in N_2$ with s.h.t. $\tT_j.$
By Claim IV.4.1(6), $a$ is of the form $\lm (\tX^2 \tY\2)$ where
$X$ and $Y$ are consecutive half-twists s.t. $X,\ Y, T_i$ create a triangle as
in Fig. IV.4.2.

Now, $\left[ a, b^{\pm 1}\right] = \left[\lm(\tX^2\ \tY\1), b^{\pm 1}\right] =
\left[\lm(\tY\2), b\1\right]_{\lm(\tX\2)}\cd [\lm (\tX^2), b^{\pm
1}].$\linebreak ( by Claim II.4).
By the previous Lemmas and the fact that $\nu\in\Center(G),$ we get:

\flushpar$\left[\lm(\tY_2\2), b^{\pm 1}\right] _{\lm(\tX^2_2)} = \cases
\nu \quad Y\ \text{and}\ T_j\ \text{are adjacent}\\
1 \quad Y\ \text{and}\ T_j\ \text{are weakly disjoint or commonly supported.}
\endcases$
\mk

\flushpar $\left[\lm(X_2\2), b^{\pm 1}\right] = \cases
\nu\ &X\ \text{and}\ T_j\ \text{are adjacent}\\
1 \ &X\ \text{and}\ T_j\ \text{are weakly disjoint or commonly supported.}
\endcases$
\flushpar
These are the only possible values since $Y$ (or $X$) are never cyclic
with $T_j.$
Since $\nu^2 = 1,$ in order to   get $[a, b^{\pm 1}] \ne 1,$ we need one
factor to be $\nu$ and the other one to be 1.  Thus we need $T_j$  adjacent to
$Y,$ and $T_j$ is weakly disjoint or commonly supported to $X,$ (or
vice-versa).
This can happen in four different cases (Fig. IV.6).  But in all of the
four cases  $T_j$ is adjacent to $T_i$ (see Fig. IV.6). \edm

\midspace{2.50in}\caption{Fig. IV.6}

\proclaim{Lemma IV.6.3}
Consider the $N$-situation.  Then:
$$(f_i)_{\tT_k} = \cases
f\1_i \nu \qquad k=i\\
f_i \qquad\ \ \ \ T_i, T_k\qquad \text{weakly disjoint}\\
f_kf_i\qquad\  T_i, T_k\qquad \text{orderly adjacent}\\
f_if_k\1\quad\ \ T_i, T_k\qquad \text{are not orderly adjacent}.\endcases$$
$$(f_i)_{\tT_k\1} = \cases
f\1_i \nu \qquad k=i\\
f_i \qquad\ \ \ \ T_i, T_k\qquad \text{weakly disjoint}\\
f_if_k\qquad\  T_i, T_k\qquad \text{orderly adjacent}\\
f_k\1 f_i\quad\ \ T_i, T_k\qquad \text{are not orderly adjacent}.\endcases$$\ep
\sk
\demo{Proof}\ Recall that $f_i$ is a prime element with s.h.t. $T_i.$
Moreover,  $\nu\in\Center(G), \nu = \left[f_i^{\pm 1},
f_j^{\pm 1}\right]$ for $T_i, T_j$ adjacent half-twists (previous Lemma).

One can  see from Fig. IV.1.6 that if $T_i$ and $T_k$ are weakly disjoint
then they are disjoint unless $k=4$ and $i=3,7,$ in which case they are
transversal.
We shall treat separately the case where $T_i$ and $T_k$ are disjoint and the
case where $T_i$ and $T_k$ are transversal.

For $k=i$ and for $T_k$ and $T_i$ disjoint  we get from
the definition of prime element that:
$$\alignat 2(f_i)_{\tT_k\1}& = f_k\1\mu
\quad &&i=k\\
(f_i)_{\tT_k\1}& = f_i\quad &&\text{for}\ T_k\ \text{and}\ T_i\
\text{disjoint.}\endalignat$$
We conjugate the above formulas by $\tT_k$ to get the correct
formula  for $(f_i)_{\tT_k}$ where $k=i$ or when $T_i$ and $T_k$ are
weakly disjoint.

For $T_i$ and $T_k$ orderly adjacent we use Proposition III.4.1 to get:

\flushpar $(f_i)_{\tT_k\1} = f_i\ f_k.$

\flushpar Thus, $(f_i)_{\tT_k\2} = (f_i)_{\tT_k\1} (f_k)_{\tT_k\1} = f_i
f_kf_k\1\ \nu =  f_i\ \nu.$

\flushpar Since $\nu^2=1,\ \nu\in\Center( G),$ we get:
\ $f_i = \nu\1 (f_i)_{\tT_k\2} = (f_i)_{\tT_k\2} \nu.$

\flushpar Thus, $(f_i)_{T_k} = (f_i)_{T_k\1} \nu=
f_i\ f_k \nu = f_i f_k [f_k\1, f_i\1]  = f_i f_k f_k\1 f_i\1  f_k f_i =
f_k f_i.$

For $T_i$ and $T_k$ be non-orderly adjacent, we use III.4.1 again to get:

\flushpar $(f_i)_{\tT_k\1} = f_k\1\ f_i.$

\flushpar Thus, $(f_k)_{\tT_k\2} = (f_k)\1_{\tT_k\1} (f_i)_{\tT_k\1} = \nu\1
f_k f_k\1 f_i= \nu  f_i \Rightarrow$

\flushpar $(f_i)_{T_k} = \nu\ (f_i)_{T_k\1} = \nu\ f_k\1 f_i = [f_i,
f_k\1]f_k\1f_i =   f_i f_k\1 f_i\1 f_k f_k\1 f_i = f_i f_k\1.$

For $T_i$ and $T_k$ transversal we get from Lemma III.4.2 that
$$(f_i)_{\tT_k\1} = f_i/$$
We conjugate the above formula by $\tT_k$ to get the correct formula for
$(f_i)_{\tT_k}$ for $T_i$ and $T_k$ transversal. \qed\edm

\bk

\subheading{\S 7. New set of generators for $\bold{G:\ \{A_j,E_j\}}$}

Recall
$$\align &E_i=\left\{ \aligned   &\G_i\quad  i\ne 2,7\\
&\G_i\quad  i=2,7\endaligned\right.\\
&\quad\\
&E_i'=(E_i)\rho_i\\
&A_j = E'_j E_j\1.\endalign$$
\proclaim{Claim IV.7.1}\ $\{A_j, E_{j}\}_{j=1}^9$ generates $G.$\ep

\demo{Proof}\ By Lemma II.2, $\{E_j, E_{j'}\}_{j=1}^9$ generates $G.$
Since $E_{j'}=A_jE_j,$ \ $\{A_j,E_j\}_{j=1}^9$ generates $G.$ \qed\edm

\sk
\demo{Definition}\  $H_1 = \tB_n\text{-orbit of}\ A_1.$\edm
\proclaim{Proposition IV.7.2}
$A=  A_1 $ is a prime element of $H_1$ with s.h.t. $\tT_1.$\ep

\demo{Proof}\ Consider the following frame of $B_9$:\
$X_1 = T_1,\ X_2 = T_3,\ X_3 = (T_5)_{T_9},\ X_4 = T_9.$\
$X_5 = T_8,\ X_6 = T_7,\ X_7 = T_6,\ X_8 = (T_2)_{T_3\1T_5\1 T_8\1
T_7\1 T_6\1}$ (see Fig. IV.7.1).

\midspace{1.00in}\caption{Fig. IV.7.1}

In order to prove that $H_1$ is a prime element of $A_1,$ we shall prove that
all the necessary conditions of Proposition III.7.1 are fulfilled
using the above frame of
$B_9.$ Let $\nu = A\cd A_{\tT_1\1}.$
It is easy to see that $\nu = E_{1'}^2\ E_1\2.$
\newline ($\nu = A\cdot A_{\tT_1\1} = A\cdot A_{E_1\1} = E_{1'} E_1\1 E_1
E_{1'}
E_1\1 E_1\1 = E_{1'}^2 E_1\2.$
\roster
\item"(0)" By definition of $H_1,$ $H_1$ is the full orbit of $A_1=A.$
\item"(1)" We have to prove $A_{\tX_2\1 \tX_1\1} = A\1A_{\tX_2\1}.$
\item"(1a)" Since $AE_1 = E_1',$ we can use Proposition II.6 and Corollary II.3
to get $\la AE_1, E_3\ra  = \la E_1, E_3\ra = 1.$
Thus, by Lemma II.4(h) we get $A_{E_2\1 E_1\1} = A\1 \cd A_{E_3\1}.$
Since $A_{\tX} = A_{\a(\tX)}$ and $\a(\tX_1) = \a(\tT_1) = E_1,$
$\a(\tX_2) = \a(\tT_3) = E_3,$ we get (1a).
\item"(1b)" $A_{\tX_1} = A_{\tT_1} = A_{\a(\tT_1)} = A_{E_1}.$  Thus, \newline
$A_{\tX_1} E_1 = A_{E_1} E_1 = E_1\1 A E_1 E_1 = E_1\1 E_1' E_1\1 E_1 E_1 =
E_1\1 E_1' E_1 = (E_1) \rho_1\1.$
By Proposition II.6 and Corollary II.3 we get
$\la A_{\tX_1}E_1,E_3\ra=\la E_1,E_3\ra=1.$ By Lemma II.4.(h) we get $
A_{\tX_1E_3\1\ E_1\1} =  A_{\tX_1}\1 A_{\tX_1 E_3\1}.$
Since $\a(\tX_1) = E_1$ and $\a (\tX_2) = \a (\tT_3) = E_3,$ we get (1b).
\item"(2a)" We need to prove $\nu_{\tX_1^2} = \nu.$
Since $\a (\tX_1) = \a (\tT_1) = E_1,$ then
$\nu_{\tX_{1}^{2}}=\nu_{E_{1}^{2}}.$ To prove $\nu_{E_1^2} = \nu,$
it is enough to prove $[E_{1'}^{2}, E_1^2] = 1.$

Consider the following 4 half-twists in $B_9$:
$T_1,\ T_3, T_2\1 T_3 T_2,\ T_2 T_1 T_2\1.$\newline
Obviously, the above half-twist consists of a ``good triangle'' in $\tB_9$ (see
Fig. IV.7.2).

\midspace{1.00in}\caption{Fig. IV.7.2}

\flushpar By Lemma III.1.2,
$$\alignat 2
&\tT_1 \cd \tT_2\1 \tT_3 \tT_2  &&= \tT_{2}\1  \tT_3 \tT_2\cd \tT_1\\
\text{and} \quad &\tT_{1}^2 \cd \tT_2\1 \tT_3^2 \tT_2  &&= \tT_3^2\cd \tT_2
\tT_1^2 \tT_2\1. \endalignat$$
By Lemma IV.3.4, there exist $\a:\tB_9\ri G$ s.t. $\a(T_i)=E_i.$ Thus,
$$\alignat 2
&[E_1 , E_2\1 E_3 E_2] &&= 1  \\
and \quad &E_{1}^2 \cd E_2\1 E_3^2 E_2  &&= E_3^2E_2\cd E_1^2 E_2\1.
\endalignat$$
Thus, $E_1^2 = E_3^2 \cd E_2 E_1^2 E_2\1 \cd E_2\1 E_3\2 E_2.$

We shall find  the commutator of $E_{1'}^2$  with each of the 3
factors $\a=E_3,$\ $\be=E_2E_1^2E_2\1,$ \ $\g=E_2\1E_3\2E_2.$

To find $[E_1^2,\g]$ we apply on $[E_1,
E_2\1 E_3 E_2] = 1$ the Invariance Theorem (Corollary I.5) with $\rho_1$\
$(m_1 = 1$ and all other $m_i = 0)$ to get $[E'_1, E_2\1 E_3 E_2] = 1$ and
thus, we get $[E_{1'}^2, E_2\1 E_3^2 E_2] = 1.$ i.e., $[E_{1'}^2,\g]=1.$

 $[E_{1'}^2,\a]=[E_{1'}^2,
E_3^2] = [(E_1^2) \rho_1, E_3^2] = \left( [E_1^2, E_3^2]\right) \rho_1 =
(\mu)\ \rho_1 = \mu\ \ \text{(By IV.5.1)}.$
$$\allowdisplaybreaks\align [E_{1'}^2,\be]=
[E_{1'}^2, E_2 E_1^2 E_2\1]
&= [E_{1'}^2, E_1\1 E_2^2 E_1] = \ \text{(Claim II.4)}\\
&= [E_{1'}^2, E_1 \cd E_1\2 E_2^2 E_1^2 E_2\2 \cd E_2^2 E_1\1]\\
&= [E_{1'}^2, E_1 \cd \mu \cd E_2^2 E_1\1] = (\mu \in\Center (G))\\
&= [E_{1'}^2, E_1  E_2^2 E_1\1]\\
&= [E_1\1, E_{1'}^2 E_1,  E_2^2]_{E_1\1}\\
&= [(E_1^2)\ \rho_1\1, E_2^2]_{E_1\1}\\
&= \left([E_1^2, E_2^2]\ \rho_1\1\right)_{E_1\1}\\
&= \left((\mu)\ \rho_1\1\right)_{E_1\1} \quad \text{(Lemma IV.5.1)}\\
&= \mu_{E_1\1} = \ (\mu \in\Center (G))\\
&= \mu .\endalign $$
To find $[E_{1'}^2,\a\be\g]$ we use Claim II.4(d),
$$\align
[E_{1'}^2, E_{1}^2] = [E_{1'}^2,\a\be\g]
&= [E_{1'}^2, \a]\ [E_{1'}^2,\be]_{\a\1} \cd [E_{1'}^2,\g]_{\be\1\g\1}\\
&= \mu \cd (\mu)_\a \cd (1)_z = (\text{since}\ \mu \in \Center
(G))\\
&= \mu \cd \mu \cd 1 = \qquad \quad (\text{since}\ \mu^2 = 1)\\
&= 1.\endalign$$
\item"(2b)"
Let $B= (A)_{\tX_2\1} = (A)_{\tT_3\1} = (A)_{E_3\1}.$\newline
To prove (2b) we have to show $\nu_B = \nu\1_{E_1}.$\newline
Now:
$$\allowdisplaybreaks\align
\nu_B &= \nu_{(A)_{{}{E_3\1}}}\\
&= A_{E_3\1}\1 \cd \nu \cd A_{E_3\1}\\
&= E_3 E_1 \underbrace{E_{1'}\1 E_3\1 \cd E_{1'}^2} E_1\2 \cd E_3 E_{1'} E_1\1
E_3\1 = \text{(By Claim II.4(a))}\\
&= \underbrace{E_3 E_1 E_3^2}\ E_{1'}\1 \underbrace{E_3\1 E_1\2 E_3}\ E_{1'}
E_1\1 E_3\1 = \text{(By Claim II.4(a))}\\
&= E_1^2 E_3 E_1 E_{1'}\1 \overbrace{E_1 E_3\2 E_1\1}\ E_{1'} E_1\1 E_3\1\\
&= E_1^2 E_3 E_1^2 \underbrace{E_1\1 E_{1'}\1 E_1}\ E_3\2 \underbrace{E_1\1
E_{1'} E_1}\ E_1\2 E_3\1 \\
&= E_1^2 E_3 E_1^2 \underbrace{(E_1) \rho_1\1 \cd E_3\2 (E_1) \rho_1\1}\ E_1\2
E_3\1 = \binom{\text{By Claim II.4(a)}}{\text{and Prop. II.6 (1)}}\\
&= E_1^2 E_3 E_1\2 \overbrace{E_3 (E_1^2) \rho_1\1 E_3\1}\ E_1\2 E_3\1\\
&= E_1^2 \cd \underbrace{E_3 E_1^2 E_3\1} \cd \underbrace{E_3^2 (E_1^2)
\rho_1\1 E_3\2} \cd \underbrace{E_3 E_1\2 E_3\1} = \binom{\text{By
Claim II.4(a) and}}{\text{Lemmas IV.3.4, IV.4.1}}\\
&= E_1^2 \cd \overbrace{E_1\1 E_3^2 E_1}\ \overbrace{\mu \cdot(E_1^{2})
\rho_1\1} \cd \overbrace{E_1\1  E_3\2 E_1} = (\mu \in \Center (G))\\
&= \mu E_1^2\cd E_1\1 E_3^2 E_1 \cd E_1\1 E_{1'}\2 E_1 \cd E_1\1 E_3\2 E_1\\
&= \mu E_1 E_3^2 E_{1'}\2 E_3\2 E_1\\
&= \mu E_1 \mu E_{1'}\2 \cd E_1 = (\text{Since}\ \mu^2, \mu \in\Center
(G))\\
&= E_1 E_{1'}\2 E_1\\
\text{Thus,}\  \nu_B &= E_1 E_{1'}\2 E_1.\endalign$$
On the other hand, $\nu_{E_1}\1 = E_1\1 (E_{1'}^2 E_1\2)\1 E_1 = E_1 E_{1'}\2
E_1.$
\item"(3)" We have to prove $A_{\tX_j} = A\ \ \fa j \geq 3.$
We recall: $A=E_{1'} E_1\1$ and $A_{X_j}=A_{\a(X_j)}.$
Since $X_3 = (T_5)_{T_9},\ X_4 = T_9,\ X_5=T_8,\ X_6=T_7$ and $X_7=T_6$,
we get $\a(X_j)$ for $j\geq 3$ is a product of $E_i$ for $i=5,9,8,7,6.$ $[E_1,
E_i] = [E_{1'}, E_i] = 1$ for $i=5, 9, 8, 7, 6$ (Proposition II.6), we get
$A_{\tX_j} = A_{\a (\tX_j)} = A$ $\fa j = 3, 4, 5, 6, 7.$\newline Now: $\tX_8 =
(\tT_2)_{\tT_3\1 \tT_5\1 \tT_8\1 \tT_7\1 \tT_6\1}.$ Thus,\newline $\a(\tX_8) =
(E_2)_{E_3\1 E_5\1 E_8\1 E_7\1 E_6'{}\1}.$ Thus,\newline
$A_{\tX_8} = A_{\a (\tX_8)} = A_{(E_2)_{{}_{E_3\1 E_5\1 E_8\1 E_7\1 E_6\1}}}.$
\newline
We first prove that $(A)_{(E_2)_{{}_{E_3\1}}} = A.$

\flushpar Since $(T_2)_{T_3\1}$ and $T_1$ are disjoint $[T_1, (T_2)_{T_3\1}] =
1.$ Thus $\a \left( [\tT_1, (\tT_2)_{\tT_3\1}]\right) = 1.$
Thus, $[E_1, (E_2)_{E_3\1}] = 1.$
We apply on this relation the invariance Theorem (Corollary I.5) to get
$[E_{1'}, (E_2)_{E_3\1}] = 1.$
Since $A = E_{1'} E_1\1$ we get $[A, (E_2)_{E_3\1}] = 1.$
Thus, $A_{(E_2)_{{}_{E_3\1}}} = A.$
Thus, $A_{\tX_8} = A_{E_5\1 E_8\1 E_7\1 E_6\1}.$
Since $E$ for $i=5, 8, 7, 6$ commutes with $E_1, E'_1$ (Proposition II.6), we
get $A_{\tX_8} = A.$
\item"(4)" Let $c = [\tX_1^2, \tX_2^2].$
We have to show that $A_c = A.$\newline
Since $c=[\tX_1^2, \tX_2^2] = [\tT_1^2, \tT_3^2],\ \a (c) = [E_1^2, E_3^2] =
\mu,$ where $\mu \in\Center(G).$
Thus, $A_c = A_\mu = A.$ \endroster
Thus all the conditions of Proposition III.7.1 are fulfilled and $A$ is a
prime element of $H_1.$\quad $\square$\enddemo
\proclaim{Proposition IV.7.3}\ $A$ is a prime element of $G$ with s.h.t.
$\tT_1.$\ep
\demo{Proof}\ By the previous propositions, $A$ is a prime element of $H_1$
with s.h.t. $T_1$ and central element $\nu = A\cd A_{\tT_1\1},$\ $ \nu
\in\Center (H_1).$
By the definition of a prime element, in order to prove that $A$ is also a
prime element of $G,$ it is enough to show that $\nu = A\ A_{\tT_1\1} \in
\Center (G).$
Since $A$ is a prime element of $H_1$ with  central element $\nu,$ we get
$(\nu)_b = \nu\ \fa b\in \tB_9.$
In particular, $(\nu)_{\tT_i} = \nu\ \fa i=1\dots 9.$
Thus, $\nu_{E_i} = \nu\ \fa i = 1 \dots 9.$
For $i\ne 1,$ we apply  $\rho$ on $\nu_{E_i} = \nu,$ to get, using the
Invariance Theorem, the relation $\nu_{E_i'} = \nu$ (for $i\ne 1$\ $(\nu)\
\rho_i=\nu$).
For $i=1$ we use $[\nu, E_1] = 1$ (from above), and $[\nu, A] = 1$
(since $A \in H_1$ and $\nu \in\Center(H_1)$) to get $[E_{1'}, \nu] = [A
E_1, \nu] = 1.$
Thus $[E_i,\nu]=[E_{i'},\nu]=1\ \fa \ i=1\dots 9\Rightarrow
\nu\in\Center(G).$ \qed\enddemo \bk

\subheading{\S 8. New set of generators for $\bold{ G\
\{E_j,h_j,\eta_j\}_{j=1}^9}$}

We introduce here a new set of generators for $G.$
In \S3 we introduced a homomorphism $\a:\tB_9\ri G$ s.t. $\a(\tT_i)=E_i$ and
introduced $G$ as a $\tB_9$-group using $\a.$ In \S5 we proved that $\a
(\xi_1) = \eta_1 = (E_2^2)_{E_1} \cd E_2\2$ is a prime element of $G$ with
s.h.t. $\tT_1$ and central element $\mu = [E_2^2, E_1^2].$
 We proved in \S7 that $A = E_{1'} E_1$ is a prime element of $G$ with s.h.t.
$\tT_1$ and central element $\nu = A_{E_1\1} = E_{1'}^2 E_1\2.$
In \S6 we introduced the $N$-situation.
Here we consider the $N$-situation for $\a:\tB_9\ri G,$ and $h_1=A\in G.$

Consider the $\un{N\text{-situation}}$ with $\eta_i,h_i,N_1,N_2$ as
follows:

$h_i$ is the unique prime element with s.h.t. $\tT_i$
s.t. $(h_i, \tT_i)$ is coherent with $(A, \tT_1).$
The central element of $h_i$ is $\nu.$

$\eta_i$ is the unique prime element with s.h.t. $\tT_i$ s.t. $(\eta_i,
\tT_i)$ is coherent with $(\eta_{1'} \tT_1).$
The corresponding central element is $\mu.$

$
N_1 = \left\{ \eta_i \quad i = 1\dots 0\right\}.$

$N_2 = \left\{ h_i \quad i = 1\dots 0\right\}.$

\proclaim{Lemma IV.8.0}

\text{\rm(i)} Let $f$ be a prime element in $G$ with s.h.t. $\tT_i.$
Then $f$ commutes with $E_i^2.$

\text{\rm(ii)} If $T_j$ is transversal to $T_i$ then $f$ commutes with $E_j.$

\text{\rm(iii)} If $T_j$ is consecutive to $T_i$ then $\left[(E_i)_{E_j^{\pm
1}},E_i^2\right]=\mu.$\ep \demo{Proof}

(i) By Lemma III.2.1, $(f)_{\tT_{i}^2}=f.$
By definition of $G$ as a $\tB_9$-group: \ $(f)_{\tT_i^2}=f_{\a(\tT_{i'}^2)}.$
Since $\a(\tT_i)=E_i$ \ (Lemma IV.3.4) we get $(f)_{E_i^2}=f.$
Thus, $f$ commutes with $E_i^2.$

(ii) By Lemma III.4.2 $(f)_{\tT_j}=f.$
But $f_{\a(\tT_j)}=f_{E_{j}}.$ Thus it commutes with $E_j.$

(iii) Lemma IV.5.1.\quad $\square$\edm
\sk
\proclaim{Lemma IV.8.1}

Let $A_i = E'_i E_i\1.$
Let $h_i, \eta_i \mu, \nu$ be as above.
Then,
\roster
\item $A_1 = A = h_1.$
\item $A_2 = h_2\1 \eta_2.$
\item $A_3 = h_3 \eta_3\1 \mu \nu.$
\item $A_4 = h_4^2 \eta_4\2 \mu \nu.$
\item $A_5 = h_5^2 \eta_5\2 \mu \nu.$
\item $A_6 = \eta_6^3\ h_6^2 \mu \nu.$
\item $A_7 = h_7 \eta_7\1 \mu \nu.$
\item $A_8 = h_8\1 \eta_8.$
\item $A_9 = h_9^2 \eta_9^{-3} \nu.$\endroster\ep

\demo{Proof} We use the definition of prime element and Lemmas II.7,  IV.5.3,
IV.6.1
 IV.6.2, IV.6.3, and IV.8.0. Recall that $\nu^2 = \mu^2 = 1, \nu, \mu
\in\Center
(G).$ We use here often the following two facts: \ If $f$ is a prime element in
$G$ with s.h.t. $\tT_i$, then $f$ commutes with $E_i^2$ (Lemma IV.8.0(i));
$(E_i^2)_{E_j^{\pm 1}}$ and $E_i^2$ commute up to $\mu$ for $T_i$ and $T_j$
consecutive half-twists (Lemma IV.8.0(iii)). \roster \item By definition of
$h_1.$   \item By Corollary II.7, $A_2=E_1\2A_1\1(A_1)_{E_2\1}(E_1^2)_{E_2\1}.$
Since $A_1$ is a prime element with s.h.t. $\tT_1$ (Proposition IV.7.3), $A_1$
commutes with $E_1^2.$ Thus,  $A_2 = A_1\1 E_1\2
(E_1^2)_{E_2\1}(A_1)_{E_2\1}$\newline $ A_1\1 = h_1\1$ from (1).\newline $E_1\2
(E_1^2)_{E_2\1} = \eta_2$ by Lemma IV.5.3.\newline $(A_1)_{E_2\1} =
(h_1)_{\tT_2\1} = h_2\1 h_1$ since $\tT_1$ and $\tT_2$ are not orderly adjacent
(by Lemma IV.6.3).

Thus, $A_2 = h_1\1 \eta_2 h_2\1 h_1.$\newline
By Lemma IV.6.2,\ $\eta_2$ commutes with $h_2\1$ and $\eta_2 h_1 = h_1 \eta_2
\nu$ $(\nu \in\Center (G)).$
Thus, $A_2 = \nu h_1\1 h_2\1 h_1 \eta_2.$
By the same Lemma,\ $h_1\1 h_2\1 h_1 = \nu h_2\1.$
Thus, $A_2 = \nu^2 h_2\1 \eta_2 = h_2\1 \eta_2.$\item
By Corollary II.7,\
$A_3=E_1\2A_1\1(A_1)_{E_3\1}(E_1^2)_{E_3\1}.$ \newline Like in (2) we can
write,
$A_3 = A_1\1 E_1\2 (E_1^2)_{E_3\1} (A_1)_{E_3\1}.$

$A_1\1
=h_1\1$ from (1). \newline
$E_1\2 (E_1^2)_{E_3\1} = \eta_3\1\ \mu.$ (Lemmas IV.5.3 and IV.5.1(d)).\newline
$(A_1)_{E_3\1} = (h_1)_{\tT_3\1} = h_1 h_3.$ (Lemma IV.6.3).\newline
Thus, $A_3 = h_1\1 \mu \eta_3\1 h_1 h_3.$\newline
Since $\left[ h_1\1, \eta_3\1\right] = \nu,$ and $\mu\in\Center (G),$\
$A_3=\mu \eta_3\1 h_1\1 \nu h_1 h_3.$
\newline
Since $\mu,\nu\in\Center(G),$\ $A_3   = \mu\nu \eta_3\1 h_3\1$\newline
Since $\eta_3$ commutes with $h_3,$\
$A_3 = h_3 \eta_3\1 \mu \nu.$
\item By Corollary II.7,
$(A_4)_{E_2\1 E_4\1} = E_4^2 A_3 E_3^2 A_2 (E_3\2)_{E_2\1} (A_3\1)_{E_2\1}
(E_4\2)_{E_2\1}.$

Since $ A_3$ is a product of prime elements with s.h.t.
$\tT_3,$
 $A_3$ commutes with $E_3^2\ (=\a(\tT_3^2))$ (Lemma IV.8.0).
Thus, $(A_4)_{E_2\1 E_4\1} =$\newline $E_4^2 E_3^2 A_3 A_2 (E_3\2)_{E_2\1}
(A_3\1)_{E_2\1} (E_4\2)_{E_2\1}.$
Now, $E_4 = \a( \tT_4),\ \tT_4$ is a half-twist which is transversal to
$\tT_3.$ Thus, by Corollary IV.8.0(ii),
$E_4\2$ commutes with $A_3\1.$
Thus,
$(A_4)_{E_2\1, E_4\1} = E_1^2 E_3^2 A_3 A_2 (E_3\2)_{E_2\1}
(E_4\2)_{E_2\1} (A_3\1)_{E_2\1}.$
$Z_1 = (T_3)_{T_2\1}$ and $Z_2 = (T_4)_{T_2\1}$ are 2
half-twists which are adjacent to $T_3$ and to $T_2$ (Fig. IV.8). Thus, by
Lemma IV.6.1,\ $[\a(Z_i^2),f]=\nu$  where $f$ is a prime element with s.h.t.
$T_3$ or $T_2$ and central element $\nu$
and
$[\a(Z_i^2),\eta]=\mu$ where $\eta$ is a prime element with s.h.t. $T_3$ or
$T_2$ and central element $\mu$ \ $(i=1,2).$

\midspace{1.00in}\caption{Fig. IV.8}

\flushpar $A_3 A_2$ is a product of 4 prime elements with s.h.t. $\tT_2$ or
$\tT_3.$ Two of them have a central element $\mu,$ and 2 of them have a central
element $\nu.$
Thus, $A_3 A_2 \a (Z_1\2) \a (Z_2\2) = (\nu \mu)^2 \a (Z_1\2)
 \a (Z_2\2) A_3 A_2$.\ Since
  $\a(Z_1)=(E_3)_{E_2\1},$\ $\a(Z_2)=(E_4)_{E_2\1}$ and
$\mu^2=\nu^2=1,$
we have
$A_3A_2(E_3\2)_{E_2\1}(E_4\2)_{E_2\1}=(E_3\2)_{E_2\1}(E_4\2)_{E_2\1}\cdot
A_3A_2.$ Thus,\newline $(A_4)_{E_2\1 E_4\1} = E_4^2 E_3^2 (E_3\2)_{E_2\1}
(E_4\2)_{E_2\1} A_3 A_2 (A_3\1)_{E_2\1}.$

{}From Corollary IV.4.2,
$\tT_4^2(\tT_4^2)_{\tT_2\1}\2=\tT_3\2(\tT_3^2)_{\tT_2\1}.$ We use Lemma
IV.4.1(4) and $c^2=1$ to exchange factors and rewrite this equation as\newline
$\tT_4^2\tT_3^2(\tT_3\2)_{\tT_2\1}(\tT_4\2)_{\tT\1}=1.$
We apply $\a$ to it to get:\newline $E_4^2
E_3^2 (E_3\2)_{E_2\1} (E_4\2)_{E_2\1} = 1.$\
Thus, $(A_4)_{E_2\1 E_4\1} = A_3
A_2 (A_3\1)_{E_2\1}.$ By (2) (3),
$$\align
(A_4)_{E_2\1 E_4\1} &= h_3 \eta_3\1 \mu \nu h_2\1 \eta_2 (\eta_3
h_3\1)_{\tT_2\1}  \mu\nu\\
 &= h_3 \eta_3\1 h_2\1 \eta_2 \eta_3 \eta_2 h_2\1
h_3\1.\ \text{(By Lemma IV.6.3)}\endalign$$
Since $[h_2\1, \eta_3\1] = \nu$ and $[\eta_2, \eta_3\1] = \mu$\ (Lemma
IV.6.2), we get
$(A_4)_{E_2\1 E_4\1} = \mu\nu h_3 h_2\1 \eta_2^2 h_2\1 h_3\1.$
Since $h_2$ commutes with $\eta_2$ and $[h_3, h_2\1] = [h_3, \eta_2] = \nu,$
we get
$(A_4)_{E_2\1 E_4\1} = \nu^5 \mu \eta_2^2 h_2\2 = \nu \mu \eta_2^2 h_2\2.$
\newline
Thus,
$$
A_4  = (\nu \mu \eta_2^2 h_2\2)_{E_4 E_2}.$$
Since $$\align(\eta_2)_{E_2} &= \eta_2\1\ \mu\\
(h_2)_{E_2} &= h_2\1 \nu\\
(h_4)_{E_2} &= h_4 h_2\1\\
(\eta_4)_{E_2} &= \eta_4 \eta_2\1,\endalign$$
So,
$$\align
A_4 &= \nu \mu \left(\eta_2\1 \mu \eta_2 \eta_4\1\right)^2 \circ
\left(h_2\1 \nu h_2 h_4\1\right)\2\\
&= \nu \mu\eta_4\2  h_4^{+2}.\endalign$$
\item By Corollay II.7.
$$\align
A_5 &= (A_4)_{E_2\1 E_3 E_7\1 E_8}=\quad\text{(from  (4))}\\
 &= \nu \mu (\eta_4)\2_{E_2\1 E_3 E_7\1\ E_8}
(h_4)^2_{E_2\1 E_3 E_7\1 E_8}.\endalign$$

By Lemma IV.2.3 , $(h_4, \tT_4)$ is coherent with $\left((h_4)_{E_2\1
E_3 E_7\1 E_8}, \tT_5\right).$
But, $(h_5, \tT_5)$ is coherent with $(h_1, \tT_1).$
Thus, $\left( (h_4)_{E_2\1 E_3 E_7\1 E_8}, \tT_5\right)$ is coherent
with $(h_1, \tT_1).$
{}From uniqueness $(h_4)_{E_2\1 E_3 E_7\1 E_8} = h_5.$
Similarly, or from Claim IV.5.3, $(\eta_4)_{E_2\1 E_3 E_7\1 E_8} = \eta_5.$
Thus,
$$A_5 = \nu \mu \eta_5\2 h_5^{+2}.$$
\item By Corollary II.7,
$$A_6=E_4\2A_4\1(A_4)_{E_6\1}(E_4^2)_{E_6\1}.$$
By (4) above, $A_4$ is a product of prime elements, with s.h.t. $\tT_4.$
\newline Thus, $A_4$ commutes with $E_1^2$.  Thus:
 $$A_6 = A_4\1 E_4\2 (E_4^2)_{E_6\1} (A_4)_{E_6\1}.$$

\flushpar By Lemma IV.5.3, $E_4\2 (E_4^2)_{E_6\1} = \eta_6.$
Thus,
$$A_6 = \nu \mu\ h_4^{+2} \eta_4\2 \eta_6 \nu \mu
(\eta_4\2)_{E_6\1} (h_4^{+2})_{E_6\1}.$$
Since $\tT_4$ and $\tT_6$ are not orderly adjacent,
$$\alignat 2
&(\eta_4)_{E_6\1} = (\eta_4)_{\tT_6\1} &&= \eta_6\1 \eta_4\\
&(h_4)_{E_6\1} &&= h_6\1 h_4.\endalignat$$
Thus,
$$\align
A_6 &= h_4\2 \eta_4^{+2} \eta_6 (\eta_6\1 \eta_4)\2 (h_6\1 h_4)^2\\
&= h_4\2\eta_4^{+2} \eta_6 \eta_4\1 \eta_6 \eta_4\1 \eta_6 h_6\1
h_4 h_6\1 h_4.\endalign$$
Since $T_4$ and $T_6\1$ are adjacent, $[ h_4, \eta_6] = [h_4, h_6] = \nu$
(Lemma IV.6.2), and thus, $h_4 h_6\1 = h_6\1 h_4 \nu$ and $\eta_6 h_4 =
h_4 \eta_6 \nu.$
Moreover, $h_4$ and $\eta_4$ commute.
Thus,
$$A_6 = \nu^9 \eta_4^{+2} \eta_6 \eta_4\1 \eta_6 \eta_4\1 \eta_6h_6\2.$$
Since $[\eta_6, \eta_4] = \mu,$ $\eta_6 \eta_4\1 = \eta_4\1 \eta_6 \mu,$ and
thus, $A_6 = \nu^9 \mu^3 \eta_6^3 h_6\2=\eta^3_5h\2_6\nu\mu.$
\item and (8): Similar to the proof of (5).
\item"(9)" By Corollary II.7, $$A_9=E_5\2A_5\1(A_5)_{E_9\1}(E_5^2)_{E_9\1}.$$
Since $A_5$ is a product of prime elements in $G$ with s.h.t. $\tT_5,$\ it
commutes with $E_5^2.$
Thus, $$A_9 = A_5\1 E_5\2 (E_5^2)_{E_9\1} (A_5)_{E_9\1}.$$
By Lemma IV.5.3, $E_5\2 (E_5^{+2})_{E_9\1} = \mu \eta_9\1.$
Thus,
$$A_9 = \nu \mu h_5\2\eta_5^2 \cd \mu \eta_9\1
\cd\left(\nu\mu\eta_5\2 h_5^{+2}\right)_{E_9\1}.$$  Since $T_5$ and $T_9$ are
orderly adjacent, $$\alignat 2
&(\eta_5)_{E_9\1} = (\eta_5)_{\tT_9\1} &&= \eta_5 \eta_9\\
&(h_5)_{E_9\1} &&= h_5 h_9.\endalignat$$
Thus,
$$A_9 = h_5\2 \eta_5^2 \mu \eta_9\1 (\eta_5 \eta_9)\2 (h_5
h_9)^2.$$
Since $[\eta_5, \eta_9] = \mu,$ then $\eta_9\1 \eta_5\1 =
\eta_5\1 \eta_9\1 \mu.$\newline
Since $[h_5, h_9] = [h_5, \eta_9] = \nu,$ then $h_9 h_5 =
h_5 h_9\1 \nu$ and $\eta_9\1 h_5 = h_5\eta_9\1 \nu.$
Thus we can collect all $h_5$ and $\eta_5$ at the left  to
get
$$A_9 = h_5\2 \eta_5^2 \eta_5\2 h_5^2 \eta_9\1 \eta_9\2 h_9^2 \cd \mu^4\
\nu^3.$$ Thus, $A_9 = \eta_9^{-3} h_9^2 \nu.$ \qed\endroster\edm
\sk
\proclaim{Lemma IV.8.2}

\text{\rm(a)} $h_i^3 = \eta_i^3 \ \ \fa i = 1\dots 9.$

\text{\rm(b)} $\nu \mu = 1.$\ep
\sk
\demo{Proof}\ For $T_i$ and $T_j$ orderly adjacent
$$\alignat 3
&(h_i)_{E_j\1 E_i\1} &&= (h_i h_j)_{E_i\1}\\
& && = \nu h_i\1 h_j h_i \ \ &&\text{(By
Lemma IV.6.3)}\\
& &&=\nu h_i\1 \nu h_i h_j  &&\ \text{(By Lemma IV.6.2)}\\
& &&= h_j &&(\text{Since}\ \nu^2 =1, \nu \in\Center (G)).
\endalignat$$
Similarly,
$$(\eta_i)_{E_j\1 E_i\1}  = \eta_j$$
Thus
$$\gather(h_j)_{E_i E_j} = h_i\\
(\eta_j)_{E_i E_j} =
\eta_i.\endgather$$
\roster
\item"(a)" By Corollary II.7, $(A_9)_{E_8\1 E_9\1} = E_9^2\ A_8 (E_9\2)_{E_8\1}
= E_9^2\ h_8\1\ \eta_8 (E_9\2)_{E_8\1}.$ Since $T_9$ is a half-twist
 adjacent to
$T_8$ and $h_8, \eta_8$ are prime elements with s.h.t. $\tT_8,$ we get by Lemma
IV.6.1 that $E_9^2\ h_8\1 = \nu h_8\1\ E_9^2$ and $E_9^2\ \eta_8 = \mu \eta_8\
E_9^2.$ Thus,
$$(A_9)_{E_8\1 E_9\1}=E_9^2\ A_8 (E_9\2)_{E_8\1} = \nu\ \mu\ h_8\1\ \eta_8\
E_9^2 (E_9\2)_{E_8\1}.$$
By Lemmas IV.5.3 and IV.8.0 $E_9^2 (E_9\2)_{E_8\1} = \mu \eta_8\1.$
Thus, $E_9^2 A_8 (E_9\2)_{E_8\1} = \nu h_8\1.$
Thus, $(A_9)_{E_8\1 E_9\1} = \nu h_8\1$ and $A_9 = \nu (h_8\1)_{E_9 E_8} =
\nu h_9\1.$
We compare this with the previous Lemma to get $\eta_9^{-3} h_9^2\nu =
h_9\1 \nu.$
Thus, $\eta_9^{-3} h_9^3 = 1.$
 Since $(\eta_9, \tT_9)$ is coherent with $(\eta_i, \tT_i)$ and
$(h_9, \tT_9)$ is coherent with $(h_i, \tT_i),$ we can use Corollary III.3.5 to
conclude that $\fa i\ \exists \tB_i \in \tB_9$ s.t. $\eta _i =
(\eta_9)_{\tB_i}$ and $h_i = (h_9)_{\tB_i}.$
Thus, $\eta_i^{-3} h_i^3 = 1 \ \ \fa i.$
\item"(b)" By Corollary II.7,
$(A_9)_{E_8\1 E_9\1} = E_9^2 A_8 (E_9\2)_{E_8\1}.$
\flushpar Thus,
$$\allowdisplaybreaks\alignat2 A_8&=E_9\2(A_9)_{E_8\1E_9\1}(E_9^2)_{E_8\1}=\
&&\text{(See IV.8.1)}\\
&=E_9\2(h_9\2\eta_9^{-3}\nu)_{E_8\1E_9\1}(E_9\2)_{E_8\1}=\ &&\text{(See
above})\\ &=E_9\2 h_8^2\eta_8^{-3}\nu(E_9^2)_{E_8\1}=\ &&\text{(See IV.6.1)}\\
&= E_9\2(E_9^2)_{E_8\1}h_8^2\eta_8^{-3}\nu^3\mu^3=\ && \text{(See IV.5.3)}\\
&=\eta_8h_8^2\eta_8^{-3}\nu\mu=\ &&\text{(See IV.6.2)}\\
&=\eta_8\2h_8^2\nu\mu.\endalignat$$
We compare with the previous result to get
$\nu\mu=1.$\quad$\square$\endroster\edm
\proclaim{Proposition IV.8.3}
Let $A_i = E_{i'} E_i\1.$
Let $h_i, \eta_i \mu, \nu$ be as before.
Then:
\roster
\item $A_1 =h_1.$
\item $A_2 = h_2\1\eta_2.$
\item $A_3 = h_3 \eta_3\1.$
\item $A_4 = h_4\1 \eta_4.$
\item $A_5 = h_5\1 \eta_5.$
\item $A_6 = h_6.$
\item $A_7 = h_7 \eta_7\1.$
\item $A_8 = h_8\1 \eta_8.$
\item $A_9 = h_9\1 \nu.$\endroster\ep
\demo{Proof} Immediately from the previous 2 lemmas.\quad $\square$\edm
\proclaim{Lemma IV.8.4} $\{E_j,\eta_j,h_j\}_{j=1}^9$ generates $G.$
\ep
\demo{Proof}
By Lemma IV.7.1, $\{A_j,E_j\}$ generates $G.$
By the previous lemma, $A_j$ is a product of $\nu, h_j, \eta_j.$
However, $\nu$ is the central element of $(h_1,\tT_1).$
Thus, $\nu=(h_1)(h_1)_{\tT_1}.$\quad $\square$\edm

\end\magnification=1200
\parindent 20 pt
\NoBlackBoxes

\define\ora{\overset\rtimes\to\alpha}
\define \a{\alpha}
\define \be{\beta}
\define \Dl{\Delta}
\define \dl{\delta}
\define \g{\gamma}
\define \G{\Gamma}
\define \lm{\lambda}

\define \r{\rho}
\define \s{\sigma}\define \si{\sigma}

\define \ve{\varepsilon}
\define \vp{\varphi}

\define \cd{\cdot}
\define \df{\dsize\frac}

\define \fa{\forall}

\define \iy{\infty}
\define \la{\langle}
\define \ra{\rangle}

\define \ri{\rightarrow}
\define \Ri{\Rightarrow}

\define \ub{\underbar}
\define \un{\underline}
\define \ov{\overline}

\define \edm{\enddemo}
\define \ep{\endproclaim}

\define \sk{\smallskip}
\define \mk{\medskip}
\define \bk{\bigskip}

\define \1{^{-1}}
\define \2{^{-2}}

\define\hb{\hat\beta}
\define\ha{\hat\alpha}

\define \CP{\Bbb C\Bbb P}

\define \CPt{\Bbb C\Bbb P^2}

\define \BR{\Bbb R}
\define \BZ{\Bbb Z}

\define \ta{\tilde{\alpha}}
\define \tB{\tilde{B}}
\define \tE{\tilde{E}}
\define \tT{\tilde{T}}
\define \tX{\tilde{X}}
\define \tY{\tilde{Y}}
\define \tP{\tilde{\Pi}}
\define \tp{\tilde{P}}
\define \tv{\tilde{v}}
\define \tZ{\tilde{Z}}
\define \tz{\tilde{z}}

\define \CG{\Cal G}

\define \Aff{\operatorname{Aff}}
\define \Gal{\operatorname{Gal}}

\define \Cen{\operatorname{Center}}\define \Center{\operatorname{Center}}

\define \Ss{S^{(6)}}
\define \Xab{X_{ab}}

\heading{\bf CHAPTER V. \ \ Construction of $\bold{G_9,\quad\hat\alpha:
G_9\overset\sim\to\rightarrow G}$}\rm\endheading \baselineskip 20pt

\bk
\subheading{Construction of $\bold{G_0(9)}$}

Let $G_0(9)$ be the group generated by $g_i$\ $i=1\dots 9$ $i\neq 4$
with the following list of relations:

$[g_1,g_2]^2=1.$

$[g_1,g_2]\in\Cen(G_0(9).)$
\mk
$[g_i,g_j]=\cases [g_1,g_2]\quad & T_i, T_j\ \text{are adjacent}\\
1\quad & T_i, T_j\ \text{are disjoint}.\endcases$
\mk
\flushpar where $T_i$ are described as follows:

\midspace{1.50 in}\caption{Fig. V.1}

Denote $\nu=[g_1,g_2].$

Let us reformulate the relations of $G_0(9)$ as follows:
\mk
\flushpar $G_0(9) =\big\la g_1,\dots,\check g_4,\dots,
g_9\bigm|[g_i,g_j]=\cases  \tau\  & T_i, T_j\ \text{are adjacent},\\
&\qquad\qquad \qquad\qquad\  \tau^2=1\qquad \tau_{g_0}=\tau \\
1\quad & T_i, T_j\ \text{are disjoint}.\endcases\big\ra$

\demo{Remark}\ $G_0(9)$ or $G_0(n)$ in general can be described in a
different way:

\flushpar Take $A_{n-1},$ a free abelian group on $n-1$ generators $A_{n-1} =
\la w_1, \dots, w_{n-1}\ra.$
Define  a skew-symmetric form on $A_{n-1}$ as follows:
$$w_i\cdot w_j = \cases
1\quad &|i-j| =1\\
0 &|i-j| \neq 1.\endcases$$
Let $G_0(n)$ be the unique central extension that satisfies
$$1\ri\BZ/2\overset b\to\ri G_0(n)\overset a\to\ri A_{n-1}\ri 1$$
where $G_0(n)$ is generated by $u_1\dots u_{n-1},$ $a(u_i)=w_i$ and
$[u_i,u_j]=\cases 1\quad & |i-j|\neq 1\\ \tau\quad & |i-j|=1\endcases$
\mk
\flushpar and
$b(1)=\tau$.
We have: $\operatorname{Ab}(G_0(n))=A_{n-1}, \quad
G_0(n)'=\{\tau,1\}\simeq\BZ/2.$\edm
\mk
\subheading{$\bold{G_0(9)}$ as a $\bold{\tB_9}$-group}
$$(g_i)_{\tT_k}=\cases g_i\1\tau\quad & i=k\\
g_i\quad & T_i\ \text{and}\ T_k \ \text{are disjoint}\\
g_kg_i\quad & T_i\ \text{and}\ T_k\ \text{are orderly adjacent}\\
g_ig_k\1\quad & T_i\ \text{and}\ T_k \ \text{are not orderly
adjacent}.\endcases$$
\proclaim{Remark V.0} Let $G_0(9)$ and $g_i$ be as above.  Then
$g_i$ is a prime element of $G_0(9)$ with s.h.t. $\tT_i$ and central element
$\tau.$\ep
\demo{Proof} By the actions of $\tB_9$ on $g_i$ and the axioms of prime
element.\edm
\bk
Consider the semidirect product: $\tB_9\rtimes G_0(9).$

\subheading{Construction of $\bold{\overset\rtimes\to\a:\tB_9\rtimes G_0(9)\ri
G}$}\

$\overset\rtimes\to\a\bigm|_{\tB_9}=\a $

\qquad$\overset\rtimes\to\a(\tT_i)=E_i$

$\qquad \ora(\xi_i)=\a(\xi_i)=\eta_i$  (Lemma IV.5.3)

$\qquad \ora(c)=\a(c)=\mu$

$\ora\bigm|_{G_0(9)}$ defined by $\ora(g_i)=h_i.$

$\qquad \ora(\tau)=\ora[g_1,g_2]=[h_1,h_2]=\nu.$

 Since $\nu$ is the  central element of
$h_i$ it belong to $\Cen (G)$ and is of order 2.
Thus, by Proposition IV.4.2  all relations that
$g_i$ satisfy are also satisfied by $h_i,$ and thus $\ora\bigm|_{G_0(9)}$ is
well-defined.

By Lemma IV.4.3, $\ha\bigm|_{G_0(9)}$ is compatible with the action of
$\tB_9$ on $G_0(9)$ and thus we have $\ora:\tB_9\rtimes G_0(9)\ri G.$

\subheading{Construction of  $\bold{G_9 }$}

Let $g_i$ be the generators of $G_0(9)$ as above.

Let $\tau=[g_1,g_2].$

Let $\xi_i$ be the prime elements in $\tB_9$ defined in Chapter IV.4
($\xi$-situation).

Let $c$ be the central element of $\xi_i.$\ $c=[\tT_1^2,\tT_2^2].$

Let $N_9\subseteq \tB_a\rtimes G_0(9)$ be the normal subgroup generated by
$\tau c$ and $(g_i\xi_i)^3$:

$N_9=\la(\tau  c\1,\quad (g_i\xi_i\1)^3\ i=1\dots 9\ i\ne 4\ra.$

$G_9 = \tB_9\rtimes G_0(9)/N_9.$

\subheading{Construction of $\bold{\ha, \ \ha: G_9\ri G}$}

By Lemma IV.4.5, \ $\nu\mu=1,$\ $(\eta_i\eta_i\1)^3=1,$ so $\ora(N_9)=1.$
Thus $\ora$ induces a map on $G_9$ denoted by $\ha.$\ \
$\ha: G_9\ri G.$

\subheading{Construction of $\bold{\hb: G\ri G_9}$}

We start by using a set of generators for $G,\ \{\G_i,\G_{i'}\}_{i=1}^9.$

We then choose as a set of generators $\{E_i,E_{i'}\}_{i=1}^9$ where
$E_i = \cases
\G_i\  & i\neq 2,7\\
\G_{i'} \  & i=2,7\endcases$
and $E_{i'}=(E_i)\vp_i.$

 The third set of generators was $\{E_i,A_i\}_{i=1}^9$ where
$A_i=E_{i'}E_i^{-1}.$

Finally, we take the following set of generators for $G:$\ $E_i,h_i,\eta_i.$

We define $\hb: G\ri G_9$ on the third set of generators as follows:

$\hb(E_i)=\tT_i\quad i=1\dots 9.$

$\hb(A_1)=g_1.$

$\hb(A_2)=g_2\1\xi_2.$

$\hb(A_3)=g_3\xi_3\1.$

$\hb(A_4)=g_4\1\xi_4.$

$\hb(A_5)=g_5\1\xi_5.$

$\hb(A_6)=g_6.$

$\hb(A_7)=g_7\xi_7\1.$

$\hb(A_8)=g_8\1\xi_8.$

$\hb(A_9)=\tau g_9\1.$

\demo{Remark} By definition of $\ha$ and by the formula for expressing $A_i$ in
terms of $h_i$ and $\eta_i$ (Lemma IV.8.3), if $\hb$ is well-defined, then
$\ha\hb=Id.$

\proclaim{Theorem V.1} $\hb$ is well-defined.\ep

\demo{Proof} We recall that $G=F_{18}/G(\ve(18))$ where $F_{18}$ is the free
group generated on $\{\G_i,\G_{i'}\}_{i=1}^9$ and $G(\ve(18))$ is the subgroup
generated by the relations induced from the factors in the braid monodromy
factorization $\ve(18)$ (see Theorem I.1) by the Van Kampen method.
To prove that $\hb$ is well-defined, we have to prove that all relations
induced by the Van Kampen Theorem are valid when each generator in a
relation is replaced by its image under $\hb.$
In what follows we shall take every braid in the braid monodromy
factorization $\ve(18)$ and use the Van Kampen method to deduce from it a
relation on $\pi_1(\CPt-S)$ in terms of $\G_i$ and $\G_{i'}.$
Then we shall present the relation in terms of $E_i$ and $A_i.$
The next step is to substitute $\tT_i$ instead of $E_i$ and $\hb(A_i)$
instead of $A_i$ and confirm that the relation holds.

Denote:

$t_i=\hb(E_i) = \tT_i.$

$t_{i'}=\hb(E_i').$

$t_{i''}=\hb(E_iE_{i'}E_i\1).$

$a_i=\hb(A_i)=\be(E_{i'}E_i\1)=t_{i'}t_i\1.$

We have expressions for $a_i$ in terms of $g_i, \xi_i,\tau.$

$a_1=g_1.$

$a_2=g_2\1\xi_2.$

$a_3=g_3\xi_3\1.$

$a_4=g_4\1\xi_4.$

$a_5=g_5\1\xi_5.$

$a_6=g_6.$

$a_7=g_7\xi_7\1.$

$a_8=g_8\1\xi_8.$

$a_9=\tau g_9\1.$

Then:

$t_i=\tT_i.$

$t_{i'}=a_it_i.$

$t_{i''}=  (a_it_i)_{t_i}=(a_i)_{t_i}\cdot t_i=t_i(a_i)_{t_i^2}=t_ia_i$ (Lemma
II.3). \qed\edm
 \medskip
$\ve(18)$ implies relations of type $\la b_1,b_2\ra=1$ relations of type
$[b_1,b_2]=1$ and relations  of type $b_1=b_2.$ We shall treat all relations of
type $\la b_1,b_2\ra$ together. To do this we prove the following two lemmas:

\proclaim{Lemma 1}

$(a_i)_{t_i} = (a_i)_{t\1_i} = \cases
\tau a\1_i\quad &i=1,6,9\\
a_i\1 &\text{otherwise}.\endcases$\ep

\demo{Proof of Lemma 1}
The elements $g_i$ and $\xi_i$ are prime elements with s.h.t. $t_i.$
By Lemma II.5, for a prime element $g$ with s.h.t.    $t,$ we have:

$(g)_t = (g)_{t\1} = \tau g\1.$

For such $g$ we then also have:

$(g\1)_t = (g\1)_{t\1} = \tau\1 g.$

Since $\tau\1=\tau$ we have for every $f$ a prime element with s.h.t. $t$ or
an inverse of a prime element with s.h.t. $t:$

$(f)_t = (f)_{t\1} = \tau f.$

For $i=1,6,9,$ one can see from the above list that $a_i$ is a prime element,
or $\tau$ multiplied by an  inverse of a prime element, and since $\tau \in\Cen
(G),$ we get the Lemma for $i=1,6,9.$

For $i\neq 1,6,9,$\ $a_i$ is a product of a prime element with an inverse of a
prime element. Write $a_i = f_i\cd f'_i,$ thus, $(a_i)_{t_i} = (a_i)_{t\1_i}$
will be the product of the inverse of the 2 factors   times $\tau^2\
(\tau\in\Cen (G)).$ By Lemma IV.6.2 such two factors commute.
Thus,
$$(a_i)_{t_i} = (f_i)_{t_i} (f'_i)_{t_i} = \tau f\1_i \tau(f'_i)\1 =
\tau^2 f_i\1 f_i'{}\1 = f_i\1 f_i'{}\1 = f_i'{}\1 f_i\1 = a_i\1$$
The same is true for  $(a_i)_{t\1_i}.$ \quad $\square$ \quad for Lemma 1\edm

\proclaim{Lemma 2}
If\ $T_i$ and $T_j$ are adjacent and $d_i=a_i$ or $(a_i)_{t_i}$ then
$(d_i)_{t\1_{{\un{j}}} t\1_{i}}=d_i\1 (d_i)_{t\1_{\un{j}}}$ for
$t_{\un{j}}=t_j$ or $t_{j'}$ or $t_{j''}.$\ep

\demo{Proof of Lemma 2}
By the above list:

$a_i = g_i$ or $g_i\1 \xi_i$ or $g_i\xi_i\1$ or $\tau g_i\1.$

By the above Lemma:

$(a_i)_{t_i} = \tau g_i\1$ or $\xi_i\1 g_i$ or $\xi_ig_i\1$ or $g_i.$

Thus, $d_i = g_i$ or $\tau g_i\1$ or $(g_i\xi\1_i)^{\pm 1} (g\1_i\xi_i)^{\pm
1}.$\edm

If $d_i$ is of the form $g_i$ or $\tau g_i\1$ it satisfies the Lemma by the
definition of prime element (axiom 2).
For symmetry reasons we shall only treat the case $d_i=g_i\1\xi_i.$

\demo{Case 1} \ $t_{\un{j}}=t_j$
$$\align
(d_i)_{t_j\1 t_i\1} &=(g_i\1\xi_i)_{t_j\1t_i\1}=(g_i)_{t_j\1t_i\1}^{-1}\cdot
(\xi_i)_{t_j\1t_i\1} \\
&\overset\text{By axiom 2 of the prime element}\to=(g_i\1
g_{{_i}t_{j}\1})\1\cdot\xi_i\1\cdot (\xi_i)_{t_j\1}\\ &=(g\1_i)_{t_j\1}\ g_i\
\xi_i\1\ (\xi_i)_{t_j\1}\\ &\overset\text{By
Proposition  IV.6.2}\to=\xi_i\1 g_i(g_i\1)_{t_j\1}\cdot(\xi_i)_{t_j\1}\\
&=(g_i\1\xi_i)\1(g_i\1\xi)_{t_j\1}=d_i\1(d_i)_{t_j\1}.\endalign$$\edm

\demo{Case 2} $t_{\un{j}}=t_{j'}$

$(d_i)_{t_j'{}\1}=(d_i)_{(t_ja_j)\1}=([a_j,d_i]d_i)_{t_j\1}=[a_jd_i]_{t\1_j}
\cdot (d_i)_{t_j\1}.$

By Proposition IV.6.2, Claim II.4 and by the formulas for $a_i,$\  $[a_j,d_i]$
is a product of $\tau.$ Since $\tau$ is of order 2, $[a_j,d_i]=\tau^\ve$\quad
$\ve=0,1.$
Thus, $(d_i)_{t_{j'}\1}=\tau^\ve(d_i)_{t_j\1}.$
So we get the claim from case 1 when multiplying each side of the equation
there by $\tau^\ve$ to get the equation for $t_{j'}.$\edm

\demo{Case 3} $t_{\un{j}}=t_{j''}$

As in Case 2 we get $(d_i)_{t_j''}=(d_i)_{t_j}\tau^\ve$ and we use Case 1
to get  Case 3.
\flushpar \qed\ for Lemma 2 \quad \edm

Since $T_i$ and $T_j$ are adjacent, $\la t_i,t_j\ra=1.$
We use the above Lemma to deduce from Lemma IV.3.1 that $\la
d_it_i,t_{\un{j}}\ra=1$\quad $\un{j} =j,j',j''.$
Since $d_i= a_i$ or $(a_i)_{t_i}$ we get $\la t_{\un{i}},t_{\un{j}}\ra=1.$
This covers all the triple relations which are induced from $\ve(18).$

\proclaim{Lemma 3}  For $i,j$ s.t. $T_i$ and $T_j$ are disjoint or
transversal, we have $[a_i, t_j] = 1$ and $[t_{\un{i}},t_{\un{j}}]=1$\quad
where $\un{i}=i$ or $i',$\quad $\un{j}=j$ or $j'.$\ep

\demo{Proof of Lemma 3} It is enough to prove
$[a_i, t_j] = [t_i,t_j]=[t_{i'},t_j]=[t_{i'},t_{j'}]=1.$

If $T_i$ and $T_j$ are disjoint, then $[T_i,T_j]=1.$
If $T_i$ and $T_j$ are transversal, then $[\tT_i,\tT_j]=1.$
In any case, $[\tT_i,\tT_j]=1$ and thus, $[t_i,t_j]=1.$

Now, $a_i$ is a prime element with s.h.t. $t_i$ or a product of 2
prime elements.
If $t_j$ is disjoint from $t_i$ we know that each of the prime factors in
$a_i$ commutes with $t_j$ by the definition of prime element.
If $t_j$ is transversal to $t_i$ then each of the prime factors in $a_i$
commutes with $t_j$ by Lemma III.3.5.
In any case, $[a_i,t_j]=1.$
Now, $t_{i'}=a_it_i,\quad' [t_{i'}, t_j] = [a_i t_i, t_j].$
Since $t_i$ commutes with $t_j$ and $a_i,$ then $[t_{i'},t_j]=[a_i,t_j],$
which equal 1.
Now, $[t_{i'},t_{j'}]=[a_it_i,a_jt_j].$
Since $t_j$ commutes with $a_i t_i$ and $a_j$ commutes with $t_i$   to prove
that $t_{i'}$ commutes with $t_{j'}$ it is enough to prove $[a_i,a_j]=1.$
This follows from Proposition IV.6.2. \quad \qed

We use the following  Lemma to show that all commutation relations in $\ve(18)$
are satisfied when $t_i$ is replacing $E_i$ and $a_i$ is replacing $A_i.$

\proclaim{Lemma 4}
Let $\tZ_{ij}^2$ be a braid in $B_{18}$ s.t. $\tZ_{ij}$ connects $q_i$ or
$q_{i'}$ with $q_j$ or $q_{j'}$ where outside of $2$ small discs centered
at $q_i, q_{i'}$ and $q_j, q_{j'}$ respectively, the path goes below the
real line, except when it goes above some of the pairs $q_k,q_{k'}$ (for
$i<k<j\
k\in K).$
Assume $T_i$ and $T_j$ are disjoint or transversal.
Then the relation induced by $\tZ_{ij}^2$ via the Van Kampen-Zariski method is
mapped to $1$ under $\be.$ \ep

\demo{Proof}
We cut $\tZ_{ij}$ into 2 pieces, one connects $u$ with the disc around $q_i$
and  $q_{i'}$ from below and the other one connects $u$ with the disc around
$q_j$ and $q_{j'}$ above the pairs $q_kq_{k'}\ k\in K.$
Thus the relation induced from $\tZ_{ij}^2$ is
$$\left[\G_{\un{i}}, \left(\prod_{k\in K}\G_k\1\G_{k'}\1\right)
\G_{\un{j}}\left(\prod_ {k\in K}\G_{k'}\G_k\right)\right]=1$$
where $\G_{\un{i}}=(\G_i)\rho_i^m$ for some $m.$
Since $(\G_i)\rho_i=\G_{i'},$\ $(\G_{i'})\rho_i=\G_{i'}\G_i\G_{i'},$ we know
that $\G_{\un i}\in\la(\G_i,\G_{i'}\ra.$

It is enough to consider $\G_{\un{i}}=\G_i$ or $\G_{i'}$ and
$\G_{\un{j}}=\G_j$ or $\G_{j'},$
 since by proving that $\G_i$ and $\G_{i'}$ commutes with $X \G_j
X\1,$ we can conclude that $X \G_j X\1$ commutes with every element $g$ from
$\la \G_i, \G_{i'}\ra;$ in particular, $\G_{\un{i}}$ commutes with $X\G_j X\1.$
Similarly,
$\G_{\un{i}}$ commutes with $X\G_{j'} X\1.$
Thus, $[X\1\G_{\un{i}} X, \G_j] = [X\1\G_{\un{i}} X, \G_{j'}] = 1.$
So, $[X\1\G_{\un{i}} X, \G_{\un{j}}] = 1.$

Now, $\be(\G_{\un{i}})=t_i$ or $t_{i'},$ $\be(\G_{\underline j})=t_j$ or
$t_{j'},$ $\be(\G_{k'} \G_k) = t_{k'}\cd t_k = a_kt_k^2.$ So
$$\be\left(\left(\prod\limits_{k\in K}\G_k\1\G_{k'} \1\right)\G_{\underline
j}\left(\prod\limits_{k\in K}\G_{k'}\G_k\right)\right)=\left(\prod\limits_{k\in
K} t_k\2 a_k\1\right)a_j^\dl t_j\left(\prod\limits_{k\in K}a_kt_k^2\right),$$
which is a product of squares of half-twists and prime elements.
Thus, we can use Propositions IV.6.1 and IV.6.2 to rearrange the factors in the
above product while multiplying the product with the appropriate $\nu^\ve$
which is a central element to get
$$\tau^\ve a_j^\dl\prod_{k\in K}a\1_{k}(a_k)_{t_j\1}\prod_{k\in K}
t_k^{\2}(t_k^2)_{t_j\1}\cdot t_j.$$
Now,
$$a_k\1(a_k)_{t_j\1}=\cases (a_k)_{t_j\1t_k\1}\quad & t_j, t_k\
\text{adjacent\qquad\qquad (Lemma 2)}\\ 1\quad & t_j,t_k \ \text{disjoint or
transversal (Lemma 3)}.\endcases$$
Since $a_k$ is a product of prime elements with $\tau$ by Lemma III.2.2,
$(a_k)_{t_j\1t_h\1}$ is also such a product.
Thus, $a_k\1(a_k)_{t_j\1}$ is 1 or a prime element supported on $t_j.$

By Lemma IV.4.0  when $t_k$ and $t_j$ are adjacent, $(t_k^{-2})(t_k^2)_{t_j\1}$
is  a product of $c$ with an inverse of  a prime element
supported on $t_j$.
By Remark III.1.1, if $t_k$ and $t_j$ are disjoint or transveral then
$t_k^{-2}(t_k^2)_{t_j\1}=1.$

Thus,
$$\align
\be&\left(\left(\prod_{k\in
K}\G_k\1\G_{k'}\1\right)\G_{\un{j}}\left(\prod_{k\in
K}\G_{k'}\G_{k'}\1)\right)\right) = \tau^\ve a_j^\si
c^t
\left(\prod\eightpoint\binom{\text{prime elements supported}}  {\text{on}\
t_j\ \text{or their inverse}} \right)\cdot t_j
 \\
&
\qquad\qquad\qquad\qquad\qquad
=c^t\tau^\ve(t_{j'}t_j\1)^\dl\left(
\prod
\eightpoint
\binom{\text{prime elements supported}
} { \text{on}
\ t_j\ \text{or their
inverse}} \right)\cdot t_j.
\endalign$$
We have to check if this product commutes  with
$\be(\G_{\un{i}}).$ Since $\G_{\un i}$ belongs to the subgroup generated by
$\G_i$ and $\G_{i'}$ which is the same as the subgroup generated by $E_i$ and
$E_{i'},$\   $\be(\G_{\un{i}})$ is a product of $\be(E_i) \quad (=t_i)$ and of
$\be(E_{i'}) \ (= t'_i).$
Thus it is enough to check that the above product commutes with $t_i$ and
$t_{i'}.$

Now, $t_i = \tT_i$ and $t_j = \tT_j$ are transversal or disjoint, thus by
Remark III.1.1(d) $t_i, t_{i'}$ commutes with $t_j, t_{j'}.$
By the definition of prime element (Axiom 3) and Lemma III.4.2, a prime element
supported on $t_j$ commutes with $t_i$ and $t_i'.$
Finally, $c,\tau\in\Center(G_9),$ and thus all the factors in  the above
product
commute with $t_i$ and $t_{i'}$ and, therefore, the product commutes. \quad
\qed\quad for Lemma 4\edm

Lemma 4 covers all degree 2 elements that appear in $\ve (18).$

We still have to check the degree 1 elements that appear in $\ve (18).$

We recall from Lemma II.2 that for
$E_i=\cases
\G_i\ \ &i\neq 2, 7\\
\G_{i'} &i = 2,7\endcases
\ \ \text{and} \ \ A_i=E_{i'} E_i\1$ we have
\roster
\item $\G_{i'} \G_i = E_{i'} E_i = A_i E_i^2.$
\item $\G_i = \cases
E_i \qquad \qquad \qquad &i\neq 2, 7\\
E_i\1 E_{i'} E_i = E_i\1 A_i E_i^2 &i = 2, 7.\endcases$

\item $\G_{i'} = \cases
E_i = A_i E_i \qquad \qquad \quad &i\neq 2, 7\\
E_i \qquad \qquad \qquad &i= 2, 7\endcases$\endroster

We have to consider each of the degree 1 elements in $\ve (18)$ and to apply
$\be$ on both sides of the induced relation in order to prove that we get the
same element of $G_9.$\newline

We want to rewrite each side  of the induced relation as a product of the
following form:\quad $\tau^\ve$ $\cd$   prime elements $\cd$
half-twists, and then use formulas for the action of $\tB_9$ on $G_9$ and
Lemma 1.

We shall only consider here one braid which is a degree 1 element.
We shall consider $\tZ_{22'(1)}.$
\mk
$\tilde z_{22'(1)} = $
\mk

$\tZ_{22'(1)}$ implies by the RMS method the relation

$\G_{2'} = \G_{1\1} \G\1_{1'} \G_2 \G_{1'} \G_1.$

{}From (2) above, this relation is actually:

$E_2 = E_1\2 A_1\1 E_2\1 A_2 E_2^2 A_1 E_1^2.$

Recall that $t_{i} = \be (E_i).$ \
 $ a_{i} = \be_i (A_i).$ Thus,
$$\align
\be (\ell s) &= t_2\\
\be (r s) &= t_1\2 a\1_1 t_2\1 a_2 t_2^2 a_1 t_1^2=\\
&=(a_1)\1_{t_1^{+2}} (a_2)_{t_2 t_1^{+2}} (a_1)_{t_1\2 t_2 t_1^{+2}} \cd t_1\2
t_2 t_1^2=\  \text{(using Lemma 1)}\\
&=a_1\1 \cd (a_2\1)_{t_1^2} (a_1)_{t_2 t_1^2} t_1\2 t_2\ t_1^2=\ \
\text{(using Lemma III.2.1(2) applied on}\ a_2\ \text{and on}\ (a_1)_{t_2})\\
&= g_1\1 (g_2\1\xi_2)\1 (g_1)_{t_2} \cd t_1\2 t_2 t_1^2=\\
&= g_1\1 \xi_2\1 g_2 g_1 g_2\1 \cd t_1\2 t_2 t_1^2=\  \text{(using}\ [g_1\1,
g_2] = \tau\ \text{and}\ [g_1\1, \xi_2\1] = \tau)
\\ &= \tau^2 \cd \xi_2 g_2 g_1\1 g_1 g_2\1\cd t_1\2 t_2 t_1^2\\
&=  \xi_2\1 t_1\2 t_2 t_1^2 t_2\1 \cd t_2\\
&=  \xi_2\1 \cd t_1\2 (t_1^2)_{t_2\1} \cd t_2\\
&=  \xi_2\1 \cd \quad \xi_2\qquad\  \cd t_2 = t_2\ \ \text{(Lemma
IV.4.2)}\quad \endalign$$
\flushpar \qed\quad for Theorem V.1
\proclaim{Lemma V.2} $\ha\hb = I_d\qquad\ \hb\ta = I_d$\ep

\demo{Proof}  By definition of $\ha$ and by the formula for expressing $A_i$ in
terms of $h_i$ and $\eta_i$ (Lemma IV.8.3),
$\ha\hb=Id.$

For $\hat\be\ha$ recall that
$G_9$ is generated by $\tT_i$ and $g_i,$ i.e., by $t_i$ and $g_i,$ i.e., by
$\tT_i, a_i$ and $\xi_i.$

$\hb\ha (t_i) = \hb\ha(\tT_i) = \hb (\tE_i) = \tT_i = t_i.$

By Proposition IV.8.3 and by $\ha (\xi_i) = \eta_i,\quad\ha (g_i) = h_i,$ we
get $\ha(a_i)=A_i.$
Thus, $\hb\ha(a_i) = \hb (A_i)\overset{\text{by def.}}\to = a_i$

In Lemma IV.5.3 we have expressions for $\eta_i,$ in terms of $E_i.$
If we apply $\hb$ on these expressions, we get the same expressions where
$E_i$ is replaced by $\tT_i.$
These expressions are exactly the expressions for $\xi_i$ as products of
$\tT_i$'s from Lemma IV.4.2.
Thus, $\hb(\eta_i)=\xi_i.$

For $\eta_4$ we could also use the expression
$\eta_4=(\eta_5)_{\tT_8\1\tT_7\tT_3\1\tT_2}.$
Apply on it the $\tB_9$ homomorphism $\hb$  to get $\hb(\eta_4)=
(\hb(\eta_5))_{\tT_8\1\dots \tT_2}\overset\text{from above}\to=
(\xi_5)_{\tT_8\1\dots \tT_2}\overset\text{Lemma IV.4.2}\to = \xi_4.$

So $\hb\ha(\xi_i)=\hb(\eta_i) = \xi_i$ \ \ by IV.4.2 and IV.5.3.

Thus, $\hat\beta\hat \a=Id.$\quad $\square$\edm

\proclaim{Corollary V.3}

$G\simeq \df{\tB_9 \ltimes G_0(9)}{N_9} = G_9$\ep
\demo{Proof} $\hat\be: G\ri G_9$ is an isomorphism.\quad$\square$

\end\magnification=1200
\parindent 20 pt

\NoBlackBoxes

\define \a{\alpha}
\define \be{\beta}
\define \g{\gamma}
\define \lm{\lambda}

\define \fa{\forall}
\define\ti{\times}

\define \tg{\tilde{\gamma}}
\define \la{\langle}
\define \ra{\rangle}

\define \G{\Gamma}

\define\BR{\Bbb R}

\define\tB{\tilde B}
\define\tX{\tilde X}
\define \CP{\Bbb C\Bbb P}
\define \CPt{\Bbb C\Bbb P^2}
\define \ve{\varepsilon}
\define \vp{\varphi}

\define \Dl{\Delta}
\define \C{\Bbb C}

\define\BZ{\Bbb Z}

\define \1{^{-1}}

\define \ri{\rightarrow}

\define \Xab{X_{ab}}
\define \Aff{\operatorname{Aff}}
\define \Gal{\operatorname{Gal}}

\define \Ab{\operatorname{Ab}}
\define\edm{\enddemo}
\define\ep{\endproclaim}
\define\un{\underline}
\define\bk{\bigskip}

\heading{\bf CHAPTER VI. \ Main Results and Formulations of Additional
Results}\rm\endheading \baselineskip 20pt

In this chapter we shall state the main results concerning the fundamental
group of the complement of a branch curve of a Veronese surface of order 3
proven in the previous chapters. The theorem is formulated in Theorem VI.1. We
shall also formulate additional results on fundamental groups that were proven
in earlier works as well as future results.
\subheading{\S1. Main results and ``forthcoming'' results}

 In order to
phrase the main results we recall a few definitions.

\demo{Definition} $\un{\text{Braid group}\ B_n=B_n[D,K]}$

Let $D$ be
a closed disc in $\Bbb R^2,$ \ $K\subset D,$ $K$ finite.
Let $B$ be the group of
all diffeomorphisms $\beta$ of $D$ such that $\beta(K) = K\,,\, \beta
|_{\partial D} = \text{Id}_{\partial D}$\,.
For $\beta_1 ,\beta_2\in B$\,, we
say that $\beta_1$ is equivalent to $\beta_2$ if $\beta_1$ and $\beta_2$ induce
the same automorphism of $\pi_1(D-K,u)$\,.
The quotient of $B$ by this
equivalence relation is called the braid group $B_n[D,K]$ ($n= \#K$).
We sometimes denote by $\overline\beta$ the braid represented by $\beta.$
The elements of $B_n[D,K]$ are called braids.
\enddemo
\demo{Definition}\ \underbar{$H(\sigma)$, half-twist defined by $\sigma$}

Let
$D,K$ be as above.
Let $a,b\in K\,,\, K_{a,b}=K-a-b$ and $\sigma$ be a simple
path in $D-\partial D$ connecting $a$ with $b$ s.t. $\sigma\cap
K=\{a,b\}.$
Choose a small regular neighborhood $U$ of $\sigma$ and an
orientation preserving diffeomorphism $f:{\Bbb R}^2 \rightarrow {\Bbb C}^1$
(${\Bbb C}^1$ is taken with usual ``complex'' orientation) such that
$f(\sigma)=[-1,1]\,,\,$ \ $ f(U)=\{z\in{\Bbb C}^1 \,|\,|z|<2\}$\,.
Let $\alpha(r),r\geq 0$\,,
be a real smooth monotone function such that $
\alpha(r) = 1$ for $r\in [0,\tsize{3\over 2}]$ and
                $\alpha(r) =   0$ for $ r\geq 2.$

Define a diffeomorphism $h:{\Bbb C}^1 \rightarrow {\Bbb C}^1$ as follows.
For $z\in {\Bbb C}^1\,,\, z= re^{i\varphi}$ let
$h(z) = re^{i(\varphi +\alpha(r))}$\,.
It is clear that on
$\{z\in{\Bbb C}^1\,|\,|z|\leq\tsize{3\over 2}\}$\,,\,$h(z)$ is the
positive rotation on $180^{\tsize{\circ}}$ and that
$h(z)=\text{Identity on }\{z\in{\Bbb C}^1\,|\,|z|\geq 2\}$\,, in
particular on ${\Bbb C}^1 -f(U)$\,.
Considering $(f\circ h\circ f^{-1})|_{D}$ (we always take composition from
left to right) we get a diffeomorphism of $D$ which switches $a$ and $b$ and is
the identity on $D-U$\,.
Thus it defines an element of $B_n[D,K],$ called the
half-twist defined by $\sigma$ and denoted $H(\sigma).$ \edm
 \demo{Definition}\  $\underline{\text{Frame of}\ B_n[D,K]}$

Let $D$ be a disc in $\BR^2.$
Let $K=\{a_1,\ldots ,a_n\}\ K\subset D.$
 Let $\sigma_1,\ldots ,\sigma_{n-1}$
be a system of simple paths in $D-\partial D$ such that each $\sigma_i$
connects
$a_i$ with $a_{i+1}$ and for
$$
i,j\in\{1,\ldots ,n-1\}\ ,\ i<j\quad ,\quad
\sigma_i\cap\sigma_j =
    \cases \emptyset \ \ &\text{if } |i-j|\geq 2\\
           a_{i+1} \ \ &\text{if } j=i+1\,.
    \endcases
$$
Let $H_i = H(\sigma_i)$\,.
We call the ordered system of (positive) half-twists $(H_1,\ldots ,H_{n-1})$ a
frame of $B_n[D,K]$ defined by $(\sigma_1,\ldots ,\sigma_{n-1})$\,, or a frame
of $B_n[D,K]$ for short. \enddemo

 \demo{Definition}
\underbar{Transversal half-twists}

The half-twists $H(\sigma_1)$ and $H(\sigma_2)$ will be called {\it
transversal} if $C_1$ and $C_2$ intersect transversally in one point which
is not an end point of either of the $\sigma_i$'s.\edm
\demo{Definition}\ $\underline{\tilde B_n}$

Let $T_n$ be the subgroup of $B_n$ normally generated by $[X,Y]$ for $X,Y$
transversal half-twists.
$\tilde B_n$ is the quotient of $B_n$ modulo $T_n.$
We choose a frame $X_i$ of $B_n.$
We denote their images in $\tilde B_n$ by $\tilde X_i.$\enddemo

\proclaim{Definition-Proposition}\ $\underline{G_0(n),\tau,u_1}$

Let $A_{n-1}$ be the free abelian group on $w_1,\dots,w_{n-1}d.$
Let us define a $\Bbb Z/2$ skew-symmetric form on $A_{n-1}$ as follows:
$$w_i\cdot w_j=\cases 1 \quad & (i-j)=1\\
0 \quad & \text{otherwise}.\endcases$$
There exists a unique central extension $G_0(n),$ of $\Bbb Z/2$ by
$A_{n-1},$ with generators $u_1\dots u_{n-1}$ that satisfies
$$\align &1\rightarrow \Bbb Z/2\overset b\to\rightarrow G_0(n)\overset
a\to\rightarrow A_{n-1}\rightarrow 1\\
&a(u_i)=w_i\\
&[u_i,u_j]=b(w_i\cdot w_j)=\cases \tau \quad &|i-j|=1\\
0\quad &\text{otherwise}.\endcases\endalign$$
We always consider $G_0(n)$ with the standard $\tilde B_n$-action as follows:
$$(u_i)_{\tilde X_{k}}=\cases u_i^{-1}\tau&\quad k=i\\
u_ku_i&\quad |i-k|=1\\
u_i&\quad |i-k|\geq 2\endcases.$$
\endproclaim
\proclaim{Claim}
$\Ab(G_0(n))=A_{n-1}$ (free abelian group on $n-1$
generators),\quad $G_0(n)'=\{\tau,1\}(\simeq \Bbb Z/2).$\ep \demo{Proof} Claim
III.6.4.\quad $\square$\edm \bk
Let $\psi_n$ be the standard homomorphisms $B_n\overset \psi_n\to\rightarrow
S_n (=$ symmetric group).

Let $\Ab$ be the standard homomorphism $B_n\overset
\Ab\to \rightarrow\Bbb Z.$

Since $\psi_n([X,Y])=1,$ and $\Ab([X,Y])=1,$ $\psi_n$
 and $\Ab$ induce  homomorphisms on $\tB_n.$

\demo{Definition}\ $\underline{\tilde\psi_n,\tilde P_n,\tilde P_{n,0},c}$

$\tilde \psi_n:\tilde B_n\rightarrow S_n,$  the induced homomorphism
from $\psi_n$.

 $\widetilde\Ab:\tilde B_n\overset \widetilde\Ab\to \rightarrow\Bbb Z,$ the
induced
homomorphism from $B_n \overset \Ab\to \rightarrow\Bbb Z.$

 $\tilde P_n=\ker
\tilde\psi_n.$

$\tilde P_{n,0}=\ker \tilde\psi_n\cap \ker \widetilde \Ab= \ker\tilde
P_n\ri\Ab(\tB_n)=\Bbb Z$.

$c=[\tilde X_1^2,\tilde X_2^2]$\quad for 2 consecutive half-twists.
\edm
\bk
Consider $\tilde B_9\ltimes G_0(9)$ with respect to the standard $\tilde
B_9$ action on $G_0(9).$
\demo{Definition}
$\underline{v_1, N_9, G_9, \tilde\psi_9: G_9\rightarrow S_9}$

$v_1=(\tilde X_2\tilde X_1\tilde X_2^{-1})^2\tilde X_2^{-2}$ for a frame
$X_1,\dots,X_8$ of $B_9.$

$N_9=$ normal subgroup generated by $c\tau^{-1},\ (u_1v_1^{-1})^3.$

$G_9=\tilde B_9\rtimes G_0(9)/{N_9}.$

$\tilde \psi_9: G_9\rightarrow S_9\quad \tilde\psi_9(\alpha,\beta)=\tilde
\psi_9(\alpha).$\edm

\demo{Definition}
$\underline{\Ab_9,H_9,H_{9,0}}$

$\Ab_9:G_9\rightarrow\Bbb Z\quad \Ab_9(\alpha,\beta)=\widetilde{\Ab}(\alpha).$

$H_9=\ker\tilde\psi_9.$

$H_{9,0}=\ker\tilde\psi_9\cap \ker \Ab_9.$\enddemo

\proclaim{Theorem VI.1} Let $V_3$ be the Veronese surface of order $3.
$ Let $S_3$ be the branch curve of a generic projection $V_3\ri \CP^2.$
Let $\Bbb C^2$ be an ``affine piece'' of $\Bbb C\Bbb P^2.$
Let $S=S_3\cap \Bbb C^2.$
 Let $G=\pi_1(\Bbb C^2-S).$
Then $G\cong G_9$ s.t. $\psi: G\rightarrow S_9$ is compatible with
$\tilde\psi_9:G_9\rightarrow S_9.$
\endproclaim
\demo{Proof} Corollary V.3.\quad $\square$\edm
\proclaim{Proposition VI.2}  Let $V_3$ be the Veronese surface of order $3.$
Let $S_3$ be the branch curve of a generic projection $V_3\ri \CP^2.$
Let $\bar G=\pi_1(\Bbb C\Bbb P^2-S_3$).
Consider $v_1$ as an element of $G_9.$
Then there exist $w_0\in H_{9,0}$ s.t. $\overline G\simeq \overline G_9$
where $\overline G_9=G_9/\langle v_1^{18}w_0\rangle.$\ep
\demo{Proof} To appear in \cite{MoTe9}.\quad $\square$\edm
\proclaim{Proposition VI.3} Let $X_i$ be a frame in $B_n.$
Let $c=[\tX_1^2,\tX_2^2].$ Then

$c=[\tilde X_1^2,\tilde
X_2^2]=[\tilde X_i^2,\tilde X_{i+1}^2]=\dots =[\tilde X_{n-2}^2,\tilde
X_{n-1}^2].$

Moreover, $(\tilde P_n')=(\tilde P_{n,0})'=\{1,c\}\simeq \Bbb
Z_2.$

$\operatorname{Ab}(\tilde P_n)=$ free abelian group on $n$ generators.

$\tilde B_n$ acts on $\tilde P_{n,0}$ by conjugation.

$\tilde P_{n,0}$ with this action is isomorphic to $G_0(n)$ with the
standard $\tilde B_n$-action as defined previously.

There exists a series:
$1\subseteq(\tilde P_{n,0})'\subseteq\tilde P_{n,0}\subseteq \tilde
P_n\subset \tilde B_n$
 s.t.
$\tilde B_n/\tilde P_n=S_n,$\quad
$\tilde P_n/\tilde P_{n,0}\simeq \Bbb Z,$\quad
$\tilde P_{n,0}/(\tilde P_{n,0})'\simeq A_{n-1}=\Ab(G_0(n)),$\quad
$(\tilde P_{n,0})'\simeq\Bbb Z_2.$

 $\operatorname{Ab}(B_n)=\BZ.$\ep
\demo{Proof} See in \cite{MoTe4}, Chapters 4, 5, definitions of $\tilde P_n$
and $\tilde P_{n,0}$ and Theorem III.6.4.\quad $\square$\edm
\proclaim{Proposition VI.4}
There exists a series
$1\subseteq H_{9,0}'\subset H_{9,0}\subset H_9\subset G_9,$
where $G_9/H_9\simeq S_9,$ \quad\qquad $H_9/H_{9,0}\simeq\Bbb Z,$\quad
\qquad  $H_{9,0}/H'_{9,0}\simeq
(\Bbb Z+\Bbb Z/3)^8,$\linebreak $H_{9,0}'=H_9'=\{1,c\}\cong \Bbb Z/2.$

\endproclaim
\demo{Proof} To appear in \cite{MoTe9}.\quad $\square$\edm
\proclaim{Proposition VI.5}
Let $\overline H_9$ and $\overline H_{9,0}$ be the images of $H_{9,0}$ and
$H_q$ in $\overline G_9.$
Then $\overline H_{9,0}'=\overline H_q'$ and
$$1\subseteq \overline H_{9,0}'\subseteq \overline H_{9,0}\subseteq
\overline H_g\subset \overline G_g$$
where
$
\overline G_9/\overline H_9\simeq S_9,$\quad
$\overline H_9/\overline H_{9,0}\simeq \Bbb Z_9,$\quad
$\overline H_{9,0}/\overline H_{9,0}'\simeq (\Bbb Z+\Bbb Z/3)^8,$\quad
$\overline H_{9,0}'\cong \Bbb Z/2.$\endproclaim
\demo{Proof} To appear in \cite{MoTe9}.\quad $\square$\edm
\bk
\subheading{\S2. The Galois cover of $\bold{V_n}$}

We shall quote here other results on fundamental groups related to Veronese
surfaces proven in earlier works.
\proclaim{Theorem} \  Let $V_n$ be a Veronese surface of order $n.$
Let $(V_n)_{\Gal}$ be its Galois cover with respect to $f,$ a generic
projection.
Then $\pi_1 ((V_n)_{\Gal}^{\Aff})$ the fundamental group of the part of
$(V_n)_{\Gal}$ that lies over a generic affine part of $\CPt$ is a
direct sum of $n^2-1$ cyclic groups of order $n.$\endproclaim
\demo{Proof} See
\cite{MoTe3}. \qed\edm
\proclaim{Theorem} \ Let $V_n$ be a Veronese surface of
order $n.$ Let $(V_n)_{\Gal}$ be its Galois cover with respect to $f,$ a
generic
projection to $\CPt.$
Then $\pi_1 (V_n)_{\Gal}$ is a direct sum of $n^2-2$ cyclic groups of order
$n.$\endproclaim
\demo{Proof}  See \cite{MoTe 3}. \qed\edm
\bk
The above results concern the computation of $\ker \ \psi$ for $\psi:
\dsize\frac{\pi_1 (\C^2-S)}{\la\G^2_j\ra} \ri S_n$ for $S,$ the branch curve of
a generic projection to $\CPt$ of $V_n,$ and $\{\G_j\}$ a $g$-base
for $\pi_1 (\C-S, *).$ To carry out the computation we used the relations
induced from the Van Kampen method (Chapter II), $\G_j^2=\G_{j'}^2=1$ and the
$RMS$ method without using the computations of Chapter IV.  These results are
easier since we assume there that all generators of fundamental groups are of
order 2. \bk \subheading{\S3.\ The Galois cover of $\bold{X_{ab}}$}

We shall also quote here a few results concerning the fundamental group of the
Galois cover of $X_{ab}.$
\demo{Definition}\ $\un{X_{ab}}$

Let $X=\CP^1\times \CP^1.$
Let $\ell_1$ be the first $\CP^1$ and $\ell_2$ the second one.
Let $a,b\in N^+.$
Look at $E=a\ell_1+b\ell_2.$
Embed $X$ in some $\Bbb P^N$ with respect to the linear system $|E|.$
Denote the image of the embedding by $X_{ab}.$\edm

\proclaim{Theorem } \ The fundamental group $\pi_1 ((\Xab)_{\Gal}^{\Aff})$
is a finite commutative group on $n-1$ generators $(n=2ab),$ each of order $C$
($c= g.c.d. (a,b))$ and there are no further relations.\endproclaim
\smallskip
\demo{Proof} Theorem 10.1 from \cite{MoTe5}. \qed \edm
\proclaim{Theorem} \ The fundamental group $\pi_1 ((\Xab)_{\Gal})$ is a
finite abelian group with $n-2$ generators, each of order $c$ ($c=g.c.d.
(a,b))$ and there are no further relations.\endproclaim

\demo{Proof} Theorem 10.2 from \cite{MoTe5}. \qed\edm
\proclaim{Corollary}
If $a, b$ are relatively  prime, then $(\Xab)_{\Gal}$ is simply
connected.\endproclaim

These results give us very interesting examples of surfaces of general type.
The Galois covers are minimal surfaces of general type.
Their index is zero for $a=b=5$ or for $a=4,$ and  $b=7$
and positive for $a\geq 5, \ b\geq 6.$
By the above Corollary, they are simply connected for $a,b$ relatively prime.
Thus we get a series of simply connected surfaces of general type with
positive index unlike the Bogomolov-watershed Conjecture (see \cite{FH}).
Moreover, $X_{55}$ is an example of a surface of general type with zero index
and even type with finite commutative fundamental group whose universal cover
is homomorphic to a connected sum of $S^2\ti S^2.$
\ $X_{55}$ gives also an exotic differential structure on a connected sum of
several copies of $S^2\ti S^2.$
There are only a few other such examples (one of them is $X_{4,7}).$
The other 3 examples will appear in \cite{MoTe10}.

In this work we have computed fundamental groups of complements of branch
curves as part of our research on algebraic surfaces.
This work also has implications to the topology of complements of curves in
general.
For general singular curves see, for example, \cite{L1} and \cite{L2}.

\end